%% file: main.tex
\documentclass[12pt,reqno]{amsart}
\input{src/preamble.tex}

\title[Flag correspondences: Clifford relations and flops]{Flag correspondences between derived Grassmannians: Clifford relations and flops}
\author{Qingyuan Jiang}
\address{Department of Mathematics, The Hong Kong University of Science and Technology, Clear Water Bay, Kowloon, Hong Kong}
\email{jiangqy@ust.hk}
\author{Shuai Zong}
\address{Department of Mathematics, The Hong Kong University of Science and Technology, Clear Water Bay, Kowloon, Hong Kong}
\email{szong@connect.ust.hk}
\date{}

\begin{document}
\begin{abstract}
    Flag correspondences between Grassmannians give rise to natural functors between their derived categories.
    For derived Grassmannians of perfect complexes, we show that composites of these functors and their adjoints fit into canonical fiber sequences, which provide a categorical refinement of the Clifford relations.
    In particular, for perfect complexes of positive rank, these functorial relations recover the spin representation constructed by Jiang, Li, and Zhao on \(K_{0}\) and realize its Fock decomposition as semiorthogonal decompositions of the derived categories.
    Moreover, when the perfect complex has rank zero, we show that the functors are spherical and obtain a ``flop--flop = twist" formula, extending results of Jiang--Leung for projectivizations and of Donovan--Segal for Grassmannian flops to derived Grassmannians over arbitrary bases.
\end{abstract}
\maketitle

\tableofcontents
\clearpage

\input{src/introduction.tex}
\input{src/flag-geometry-rees.tex}
\input{src/clifford-relations.tex}
\input{src/spherical-flags.tex}
\input{src/abel-applications.tex}

\appendix
\input{src/incidence-and-extended-clifford-relations.tex}
\clearpage
\printbibliography[heading=bibliography,title={References}]
\end{document}

%% file: src/preamble.tex
\usepackage{amsmath,amssymb,amsthm,mathtools}
\usepackage[dvipsnames]{xcolor}
\usepackage{microtype}
\usepackage{array}
\usepackage{tikz-cd}
\usepackage{enumitem}
\usepackage[a4paper,top=25mm,bottom=25mm,left=25mm,right=25mm,
    includefoot,footskip=30pt]{geometry}
\usepackage{indentfirst}

\usepackage{xurl}
\usepackage[backend=biber,style=alphabetic,giveninits=true,maxbibnames=99,
    maxalphanames=99,doi=true,url=false,isbn=false]{biblatex}
\definecolor{paperlink}{HTML}{8B1E3F}
\usepackage[unicode,colorlinks=true,
    linkcolor=black,citecolor=paperlink,
    urlcolor=paperlink]{hyperref}
\usepackage{bookmark}
\usepackage{zref-clever}
\zcsetup{cap,noabbrev,countertype={maintheorem=theorem,subsection=subsection}}
\zcLanguageSetup{english}{
    type=subsection,
    Name-sg=Subsection,Name-pl=Subsections,
    name-sg=subsection,name-pl=subsections
}
\AddToHook{cmd/appendix/after}{\zcsetup{countertype={subsection=subsection}}}
\hypersetup{
    pdftitle={Flag correspondences between derived Grassmannians: Clifford relations and flops},
    pdfauthor={Qingyuan Jiang and Shuai Zong}
}

\makeatletter
\renewcommand{\paragraph}{\@startsection{paragraph}{4}%
    {\z@}{\z@}{-\fontdimen2\font}%
    {\normalfont\bfseries}}
\makeatother

\newtheoremstyle{paperplain}
    {6pt plus 2pt minus 1pt}{6pt plus 2pt minus 1pt}
    {\itshape}{}{\normalfont\small\scshape}{.}{.5em}
    {\thmname{#1}\thmnumber{ #2}\thmnote{\normalfont\small\ (#3)}}
\newtheoremstyle{paperdefinition}
    {6pt plus 2pt minus 1pt}{6pt plus 2pt minus 1pt}
    {\normalfont}{}{\normalfont\small\scshape}{.}{.5em}
    {\thmname{#1}\thmnumber{ #2}\thmnote{\normalfont\small\ (#3)}}

\theoremstyle{paperplain}
\newtheorem{theorem}{Theorem}[section]
\newtheorem{lemma}[theorem]{Lemma}
\newtheorem{proposition}[theorem]{Proposition}
\newtheorem{corollary}[theorem]{Corollary}

\newtheorem{maintheorem}{Theorem}

\theoremstyle{paperdefinition}
\newtheorem{definition}[theorem]{Definition}
\newtheorem{example}[theorem]{Example}
\newtheorem{remark}[theorem]{Remark}
\numberwithin{equation}{section}

\AddToHook{env/lemma/begin}{\zcsetup{countertype={theorem=lemma}}}
\AddToHook{env/proposition/begin}{\zcsetup{countertype={theorem=proposition}}}
\AddToHook{env/corollary/begin}{\zcsetup{countertype={theorem=corollary}}}
\AddToHook{env/definition/begin}{\zcsetup{countertype={theorem=definition}}}
\AddToHook{env/example/begin}{\zcsetup{countertype={theorem=example}}}
\AddToHook{env/remark/begin}{\zcsetup{countertype={theorem=remark}}}

\DeclareMathOperator{\id}{id}
\DeclareMathOperator{\Map}{Map}
\DeclareMathOperator{\Hom}{Hom}
\DeclareMathOperator{\Aut}{Aut}
\DeclareMathOperator{\Perf}{Perf}
\DeclareMathOperator{\QCoh}{QCoh}
\DeclareMathOperator{\fib}{fib}
\DeclareMathOperator{\cofib}{cofib}
\DeclareMathOperator{\im}{im}

\DeclareMathOperator{\Gr}{Gr}
\DeclareMathOperator{\Flag}{Flag}
\DeclareMathOperator{\Incid}{Incid}
\DeclareMathOperator{\Spec}{Spec}
\DeclareMathOperator{\Pic}{Pic}

\DeclareMathOperator{\gr}{gr}
\DeclareMathOperator{\Sym}{Sym}

%% file: src/introduction.tex
\section{Introduction}
\label{sec:introduction}

A basic approach in geometric representation theory is to realize algebraic operators by correspondences between varieties and to study their relations through convolution.
In the classical work of Beilinson, Lusztig, and MacPherson \cite{BeilinsonEtAl1990QuantumDeformationGL}, convolution on varieties of pairs of partial flags gives a geometric construction of the quantum group \(U_{q}(\mathfrak{gl}_{n})\).
Grassmannians are among the simplest partial flag varieties, and two-step flags provide natural correspondences between Grassmannians of different ranks.

Schubert classes give another natural connection between Grassmannians and representation theory through exterior powers.
Let \(V\) be an \(n\)-dimensional complex vector space.
After choosing a complete flag in \(V\), the Schubert classes of \(\Gr(d,V)\) form a basis of its cohomology indexed by \(d\)-element subsets \(I = \{i_{1} < \cdots < i_{d}\}\) of \(\{1,\ldots,n\}\).
Identifying the class indexed by \(I\) with \(e_{i_{1}} \wedge \cdots \wedge e_{i_{d}}\), where \(e_{1},\ldots,e_{n}\) is the standard basis of \(\mathbb{Q}^{n}\), gives an isomorphism of vector spaces
\[
    \bigoplus_{d = 0}^{n}H^{*}(\Gr(d,V),\mathbb{Q}) \simeq \bigwedge\nolimits^{\bullet}\mathbb{Q}^{n},
\]
after forgetting the cohomological grading.
The exterior algebra carries the spin representation of the Clifford algebra, with exterior multiplication and contraction acting as creation and annihilation operators.

The Clifford action admits a natural geometric realization through flag correspondences between Grassmannians.
Xiao \cite{Xiao2015CohomologicalHallAlgebraA1} identified Grassmannians with moduli spaces of stable framed representations of the \(A_{1}\)-quiver and used Hall algebra techniques to establish the Clifford relations at the level of cohomology.
His construction thus connects the spin representation arising from Schubert calculus with the convolution of geometric correspondences.

A more general setting for these constructions is provided by derived Grassmannians of perfect complexes \cite{Jiang2023DerivedGrassmanniansSchurFunctors}.
Such Grassmannians arise naturally in the study of linear series on curves \cite[Chapter~IV, \S3]{ArbarelloEtAl1985GeometryAlgebraicCurvesI} and moduli spaces of sheaves on K3 surfaces \cite{Markman2001BrillNoetherDualityK3,AddingtonTakahashi2022CategoricalSl2Action}, providing a broader geometric setting for the study of Clifford relations.
The first-named author, Li, and Zhao \cite{JiangEtAl2025InstantonsAffineVertexAlgebra} established the Clifford relations on \(K_{0}\) and constructed spin representations for derived Grassmannians.
They further applied these results to obtain \(\widehat{\mathfrak{gl}}_{r}\)-actions on various cohomology theories of instanton moduli spaces on blown-up surfaces, answering a question of Vafa and Witten \cite{VafaWitten1994StrongCouplingTest}.

While these results concern cohomology or \(K_{0}\), flag correspondences also induce Fourier--Mukai functors between derived categories.
We study the relations among these functors and their adjoints for derived Grassmannians by comparing their compositions directly.
The resulting relations recover the Clifford action on \(K_{0}\), give semiorthogonal decompositions for perfect complexes of positive rank, and lead to spherical functors and flop--flop formulas in rank zero.

\subsection{Main results}

Let \(X\) be a prestack and let \(\mathcal{E}\) be a perfect complex of rank \(r\) and Tor-amplitude in \([0,1]\) over \(X\); locally, \(\mathcal{E}\) is a two-term complex \([\mathcal{W} \to \mathcal{V}]\) of vector bundles with \(\mathcal{V}\) in degree \(0\).
(Readers less familiar with derived algebraic geometry may assume that \(X\) is a scheme and \(\mathcal{E}\) is a vector bundle, and interpret the constructions below classically.)
For \(d \geq 0\), let
\[
    \operatorname{pr}_{d} \colon \Gr_{X}(\mathcal{E},d) \to X
\]
denote the derived Grassmannian of rank \(d\) quotients \cite{Jiang2023DerivedGrassmanniansSchurFunctors}, with universal fiber sequence
\[
    R_{d} \longrightarrow \operatorname{pr}_{d}^{*}\mathcal{E} \longrightarrow Q_{d}.
\]

The natural forgetful maps from the flag variety give the correspondence
\[
    \Gr_{X}(\mathcal{E},d) \xleftarrow{j_{d}} \Flag_{X}(\mathcal{E};d,d + 1) \xrightarrow{\pi_{d}} \Gr_{X}(\mathcal{E},d + 1).
\]
For each \(a \in \mathbb{Z}\), we consider the Fourier--Mukai functor
\[
    \Psi_{d,a} \colon M \longmapsto j_{d,*}\bigl(L_{d + 1}^{a} \otimes \pi_{d}^{*}M\bigr),
\]
where \(L_{d + 1} = \ker(Q_{d + 1} \to Q_{d})\) is a line bundle on the flag variety.
Its left adjoint is
\[
    \Psi_{d,a}^{L} \colon M \longmapsto \pi_{d,!}\bigl(L_{d + 1}^{-a} \otimes j_{d}^{*}M\bigr),
\]
where \(\pi_{d,!}(-) = \pi_{d,*}(- \otimes \omega_{\pi_{d}})\) and \(\omega_{\pi_{d}}\) is the relative dualizing complex.

For \(d \geq 1\), we want to compare the compositions of the following forms
\[
    \Psi_{d - 1,b}^{L}\Psi_{d - 1,a}^{\vphantom{L}}, \qquad \Psi_{d,a^{\prime}}^{\vphantom{L}}\Psi_{d,b^{\prime}}^{L}.
\]
Both compositions are Fourier--Mukai transforms, whose kernels are obtained by convolution along the following fiber products of flag varieties:
\[
    \begin{aligned}
        Z^{-} = {} & \Flag_{X}(\mathcal{E};d - 1,d) \times_{\Gr_{X}(\mathcal{E},d - 1)} \Flag_{X}(\mathcal{E};d - 1,d), \\
        Z^{+} = {} & \Flag_{X}(\mathcal{E};d,d + 1) \times_{\Gr_{X}(\mathcal{E},d + 1)} \Flag_{X}(\mathcal{E};d,d + 1).
    \end{aligned}
\]

There is a natural correspondence between \(Z^{-}\) and \(Z^{+}\) away from their diagonals.
A point of \(Z^{-}\) consists of two rank \(d\) quotients \(V_{1},V_{2}\) of \(\mathcal{E}\), together with a common rank \(d - 1\) quotient \(W\).
When \(V_{1}\) and \(V_{2}\) are distinct, the induced map
\[
    \mathcal{E} \longrightarrow V_{1} \times_{W}V_{2}
\]
is a rank \(d + 1\) quotient, giving a point of \(Z^{+}\).
Conversely, two distinct rank \(d\) quotients \(V_{1},V_{2}\) of a common rank \(d + 1\) quotient \(V\) have kernel lines spanning a rank two subbundle of \(V\).
Quotienting by this subbundle gives a common rank \(d - 1\) quotient.
Thus we obtain a natural identification
\[
    Z^{-} \setminus \Delta_{-} \simeq Z^{+} \setminus \Delta_{+},
\]
where \(\Delta_{-}\) and \(\Delta_{+}\) are the respective diagonals.

This identification suggests considering the blow-ups of \(Z^{-}\) and \(Z^{+}\) along their diagonals.
We show that these blow-ups are naturally identified with the space \(\widetilde{Z}\) parametrizing two rank \(d\) quotients sharing both a rank \(d - 1\) quotient and a rank \(d + 1\) extension.
More precisely, there is a diagram
\[
    Z^{-}\xleftarrow{b_{-}}\widetilde{Z} \xrightarrow{b_{+}}Z^{+},
\]
where
\[
    \widetilde{Z} \simeq \operatorname{Bl}_{\Delta_{-}}(Z^{-}) \simeq \operatorname{Bl}_{\Delta_{+}}(Z^{+}).
\]
Here, the blow-ups are taken in the derived sense of \cite{HekkingEtAl2025BlowupsDerivedWeilRestrictions}.
Using this common blow-up, we compare the two compositions of flag functors and obtain the following result.

\begin{maintheorem}[\zcref{thm:mixed-flag-rees-comparison}]
    \label{thm:intro-mixed-comparison}
    For \(d \geq 1\) and \(a,b \in \mathbb{Z}\), put \(s = a - b\) and \(u = -r - s\).
    There are endofunctors \(\mathcal{M}_{a,b}\), \(\mathcal{H}^{+}_{a,b}\), and \(\mathcal{H}^{-}_{a,b}\) of \(\Perf(\Gr_{X}(\mathcal{E},d))\) and canonical fiber sequences
    \[
        \begin{aligned}
            \mathcal{M}_{a,b} & \longrightarrow \Psi_{d - 1,b}^{L}\Psi_{d - 1,a} \longrightarrow \mathcal{H}^{+}_{a,b}, \\
            \mathcal{M}_{a,b} & \longrightarrow \Psi_{d,a + 1}\Psi_{d,b + 1}^{L}[-1] \longrightarrow \mathcal{H}^{-}_{a,b}.
        \end{aligned}
    \]

    The functors \(\mathcal{H}^{+}_{a,b}\) and \(\mathcal{H}^{-}_{a,b}\) vanish for \(s < 0\) and \(u < 0\), respectively.
    For \(s \geq 0\) and \(u \geq 0\), respectively, they admit finite filtrations with graded pieces
    \[
        \begin{aligned}
            \gr_{t}\mathcal{H}^{+}_{a,b} & \simeq - \otimes \Sym^{t}R_{d} \otimes \Gamma^{s - t}Q_{d},                                                           &  & 0 \leq t \leq s, \\
            \gr_{t}\mathcal{H}^{-}_{a,b} & \simeq - \otimes \bigl(\Sym^{t}R_{d} \otimes \Gamma^{u - t}Q_{d}\bigr)^{\vee} \otimes \det \mathcal{E}^{\vee}[2d - r], &  & 0 \leq t \leq u.
        \end{aligned}
    \]
\end{maintheorem}

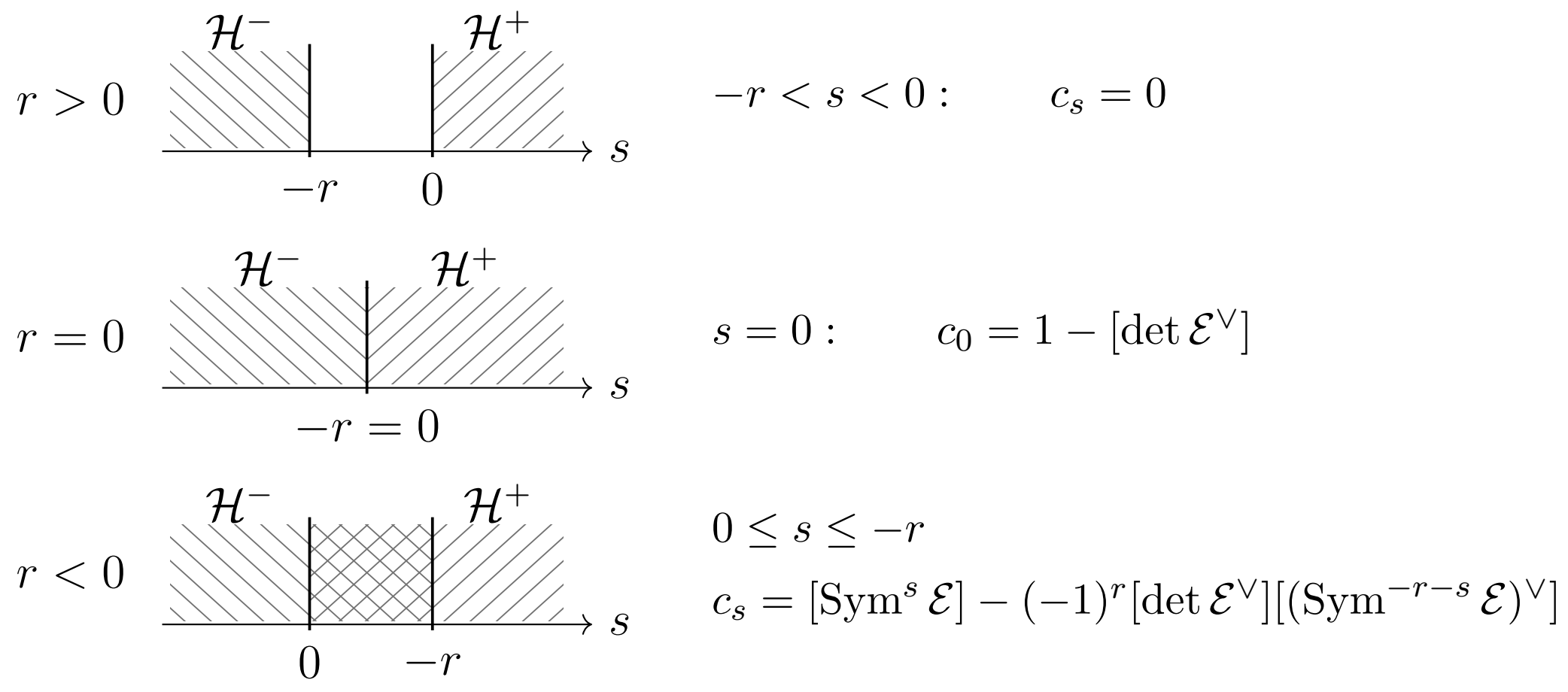
\begin{figure}[htbp]
    \centering
    \begin{tikzpicture}[x=0.72cm,y=0.65cm,
            hatch/.style={draw=black!55,line width=0.35pt},
            wall/.style={line width=0.7pt},
            pairing/.style={anchor=west,text width=8.8cm,align=left,
                    inner sep=0pt,font=\small}]
        \begin{scope}[shift={(0,6.4)}]
            \node[pairing] at (6.7,0.75) {
                \mbox{\(-r < s < 0:\qquad c_{s} = 0\)}
            };
            \node[anchor=east] at (-0.25,0.7) {$r > 0$};
            \begin{scope}
                \clip (0.1,0.05) rectangle (1.8,1.35);
                \foreach \x in {-1.4,-1.08,...,5.6}
                \draw[hatch] (\x,1.4) -- ++(1.4,-1.4);
            \end{scope}
            \begin{scope}
                \clip (3.3,0.05) rectangle (4.9,1.35);
                \foreach \x in {-1.4,-1.08,...,5.6}
                \draw[hatch] (\x,0) -- ++(1.4,1.4);
            \end{scope}
            \draw[->] (0,0) -- (5.25,0) node[right] {$s$};
            \draw[wall] (1.8,-0.08) -- (1.8,1.45);
            \draw[wall] (3.3,-0.08) -- (3.3,1.45);
            \node[below] at (1.8,-0.08) {$-r$};
            \node[below] at (3.3,-0.08) {$0$};
            \node at (0.95,1.65) {$\mathcal{H}^{-}$};
            \node at (4.1,1.65) {$\mathcal{H}^{+}$};
        \end{scope}
        \begin{scope}[shift={(0,3.2)}]
            \node[pairing] at (6.7,0.75) {
                \mbox{\(s = 0:\qquad c_{0} = 1 - [\det \mathcal{E}^{\vee}]\)}
            };
            \node[anchor=east] at (-0.25,0.7) {$r = 0$};
            \begin{scope}
                \clip (0.1,0.05) rectangle (2.5,1.35);
                \foreach \x in {-1.4,-1.08,...,5.6}
                \draw[hatch] (\x,1.4) -- ++(1.4,-1.4);
            \end{scope}
            \begin{scope}
                \clip (2.5,0.05) rectangle (4.9,1.35);
                \foreach \x in {-1.4,-1.08,...,5.6}
                \draw[hatch] (\x,0) -- ++(1.4,1.4);
            \end{scope}
            \draw[->] (0,0) -- (5.25,0) node[right] {$s$};
            \draw[wall] (2.5,-0.08) -- (2.5,1.45);
            \node[below] at (2.5,-0.08) {$-r = 0$};
            \node at (1.3,1.65) {$\mathcal{H}^{-}$};
            \node at (3.7,1.65) {$\mathcal{H}^{+}$};
        \end{scope}
        \begin{scope}
            \node[pairing] at (6.7,0.75) {
                \(0 \leq s \leq -r\)\\[5pt]
                \mbox{\(c_{s} = [\Sym^{s}\mathcal{E}]
                    - (-1)^{r}[\det \mathcal{E}^{\vee}]
                    [(\Sym^{-r - s}\mathcal{E})^{\vee}]\)}
            };
            \node[anchor=east] at (-0.25,0.7) {$r < 0$};
            \begin{scope}
                \clip (0.1,0.05) rectangle (3.3,1.35);
                \foreach \x in {-1.4,-1.08,...,5.6}
                \draw[hatch] (\x,1.4) -- ++(1.4,-1.4);
            \end{scope}
            \begin{scope}
                \clip (1.8,0.05) rectangle (4.9,1.35);
                \foreach \x in {-1.4,-1.08,...,5.6}
                \draw[hatch] (\x,0) -- ++(1.4,1.4);
            \end{scope}
            \draw[->] (0,0) -- (5.25,0) node[right] {$s$};
            \draw[wall] (1.8,-0.08) -- (1.8,1.45);
            \draw[wall] (3.3,-0.08) -- (3.3,1.45);
            \node[below] at (1.8,-0.08) {$0$};
            \node[below] at (3.3,-0.08) {$-r$};
            \node at (0.95,1.65) {$\mathcal{H}^{-}$};
            \node at (4.1,1.65) {$\mathcal{H}^{+}$};
        \end{scope}
    \end{tikzpicture}
    \caption{The two walls in \zcref{thm:intro-mixed-comparison} and the corresponding coefficients \(c_{s}\) defined in \eqref{eq:clifford-coefficient-definition}.
        Hatching indicates the ranges of the two filtrations.}
    \label{fig:comparison-walls}
\end{figure}

The two terms \(\mathcal{H}^{+}_{a,b}\) and \(\mathcal{H}^{-}_{a,b}\) are governed by the walls \(s = 0\) and \(s = -r\), respectively.
For \(r > 0\), both terms vanish in the range \(-r < s < 0\), giving an equivalence
\[
    \Psi_{d - 1,b}^{L}\Psi_{d - 1,a} \simeq \Psi_{d,a + 1}\Psi_{d,b + 1}^{L}[-1].
\]
For \(r = 0\), the two walls coincide, while for \(r < 0\), both filtrations occur throughout \(0 \leq s \leq -r\).
These three cases are illustrated in \zcref{fig:comparison-walls}.

Passing to \(K_{0}\), the common term \(\mathcal{M}_{a,b}\) in the two fiber sequences cancels, while the shift \([-1]\) in the second sequence introduces a sign.
Thus, the sum of the operators induced by the two compositions is given by \([\mathcal{H}^{+}_{a,b}] - [\mathcal{H}^{-}_{a,b}]\).
By the above filtrations and the generalized Serre theorem \cite{Jiang2025DerivedProjectivizations}, this difference acts, via the natural \(K_{0}(\Perf(X))\)-action, as multiplication by the class
\begin{equation}
    \label{eq:clifford-coefficient-definition}
    c_{s} = [p_{*}\mathcal{O}_{\mathbb{P}_{X}(\mathcal{E})}(s)] \in K_{0}(\Perf(X)), \qquad s \in \mathbb{Z},
\end{equation}
where \(p \colon \mathbb{P}_{X}(\mathcal{E}) \to X\) is the projection.
After normalizing the flag functors by determinant line bundles, this identity becomes the mixed Clifford anticommutation relation for the operators introduced below.

Put
\[
    \mathcal{C} = \bigoplus_{d \geq 0}\Perf(\Gr_{X}(\mathcal{E},d)).
\]
For each \(i \in \mathbb{Z}\) and \(d \geq 0\), define the functors
\[
    \begin{aligned}
        F_{i} \colon \Perf(\Gr_{X}(\mathcal{E},d + 1)) & \longrightarrow \Perf(\Gr_{X}(\mathcal{E},d)), \\
        M                                              & \longmapsto \Psi_{d,i - r}(M) \otimes (\det Q_{d})^{-1},
    \end{aligned}
\]
and
\[
    \begin{aligned}
        F_{i}^{L} \colon \Perf(\Gr_{X}(\mathcal{E},d)) & \longrightarrow \Perf(\Gr_{X}(\mathcal{E},d + 1)), \\
        N                                              & \longmapsto \Psi_{d,i - r}^{L}(N \otimes \det Q_{d}).
    \end{aligned}
\]
Together, these give an adjoint pair \(F_{i}^{L} \dashv F_{i}\) of endofunctors of \(\mathcal{C}\).

\begin{maintheorem}[\zcref{lem:flag-squares} and \zcref{thm:k-theory-clifford-relations}]
    \label{thm:intro-clifford}
    On \(K_{0}(\mathcal{C})\), the induced \(K_{0}(\Perf(X))\)-linear operators \(f_{i} = [F_{i}]\) and \(f_{i}^{*} = [F_{i}^{L}]\) are square-zero and satisfy, for all \(i,j \in \mathbb{Z}\),
    \[
        \begin{aligned}
            f_{i}f_{j} + f_{j}f_{i} & = f_{i}^{*}f_{j}^{*} + f_{j}^{*}f_{i}^{*} = 0, \\
            f_{i}^{*}f_{j} + f_{j}f_{i}^{*} & = c_{j - i}\,\id.
        \end{aligned}
    \]
\end{maintheorem}

For \(r > 0\), after a change of generators, these relations give an action of the standard Clifford algebra with \(r\) pairs of generators on \(K_{0}(\mathcal{C})\).
This Clifford module structure was constructed by Jiang, Li, and Zhao \cite{JiangEtAl2025InstantonsAffineVertexAlgebra}, who obtained it from the semiorthogonal decompositions of \cite{Jiang2023DerivedCategoriesGrassmannians} by viewing a semiorthogonal decomposition as a categorical representation.
\zcref{thm:intro-clifford} instead proves the relations directly from the flag functors, for every rank \(r\).
Reversing the direction of \cite{JiangEtAl2025InstantonsAffineVertexAlgebra}, we then use these categorical relations to give a new proof of semiorthogonal decompositions of derived Grassmannians over arbitrary bases (\zcref{thm:categorical-fock-decomposition} and \zcref{prop:subset-flag-words}).

When \(r = 0\), the two walls in \zcref{thm:intro-mixed-comparison} coincide.
For each flag functor \(\Psi_{d,a}\), let
\[
    \begin{aligned}
        C_{d,a} & = \fib\bigl(\Psi_{d,a}^{L}\Psi_{d,a} \to \id\bigr), \\
        T_{d,a} & = \cofib\bigl(\Psi_{d,a}\Psi_{d,a}^{R} \to \id\bigr)
    \end{aligned}
\]
denote its cotwist and twist, respectively.

Taking \(a = b\) in \zcref{thm:intro-mixed-comparison}, for \(d \geq 1\) we obtain the fiber sequences
\[
    C_{d - 1,a} \longrightarrow \Psi_{d - 1,a}^{L}\Psi_{d - 1,a} \longrightarrow \id,
\]
and
\[
    \bigl(- \otimes \det \mathcal{E}^{\vee}\bigr) T_{d,a + 1}[2d - 1] \longrightarrow \Psi_{d,a + 1}\Psi_{d,a + 1}^{L}[-1] \longrightarrow \bigl(- \otimes \det \mathcal{E}^{\vee}\bigr)[2d].
\]
The last arrow in the first sequence is an adjunction counit, while the last arrow in the second is identified with a right adjunction counit under relative Serre duality.
Since \zcref{thm:intro-mixed-comparison} identifies the fibers of these two maps, we obtain
\[
    C_{d - 1,a} \simeq \bigl(- \otimes \det \mathcal{E}^{\vee}\bigr) T_{d,a + 1}[2d - 1].
\]
This relation allows us to establish sphericality of the flag functors.
Comparing their spherical twists with the incidence equivalences between the Grassmannians of \(\mathcal{E}\) and \(\mathcal{E}^{\vee}[1]\) then gives the following flop--flop formulas.

\begin{maintheorem}[\zcref{thm:spherical-flags,thm:flop-product}]
    \label{thm:intro-rank-zero}
    Suppose \(r = 0\).
    For every \(d \geq 0\) and \(a \in \mathbb{Z}\), the flag functor \(\Psi_{d,a}\) is spherical.

    For \(k \in \mathbb{Z}\), denote the functors associated with the incidence correspondences in \eqref{eq:rank-zero-incidence-functors} by
    \[
        \begin{aligned}
            \Phi_{d,k} \colon \Perf(\Gr_{X}(\mathcal{E}^{\vee}[1],d)) & \xrightarrow{\sim} \Perf(\Gr_{X}(\mathcal{E},d)), \\
            \Phi_{d,k}^{\vee} \colon \Perf(\Gr_{X}(\mathcal{E},d))    & \xrightarrow{\sim} \Perf(\Gr_{X}(\mathcal{E}^{\vee}[1],d)).
        \end{aligned}
    \]
    There are natural equivalences of \(X\)-linear endofunctors:
    \[
        \begin{aligned}
            \Phi_{d,k}\Phi_{d,k + 1}^{-1}                & \simeq \bigl(- \otimes \det \mathcal{E}^{\vee}\bigr)T_{d,-k}, \\
            \bigl(\Phi_{d,0}\Phi_{d,0}^{\vee}\bigr)^{-1} & \simeq T_{d,d}\cdots T_{d,2}T_{d,1}.
        \end{aligned}
    \]
\end{maintheorem}

\zcref{thm:intro-rank-zero} provides a common generalization of two earlier results: for \(d = 1\), it recovers the flop--flop formula of Jiang and Leung for projectivizations \cite{JiangLeung2022ProjectivizationsFlops} (\zcref{rem:d=1}), and for the local Grassmannian flop, it recovers the Grassmannian twists of Donovan and Segal \cite{DonovanSegal2014WindowShiftsGrassmannianTwists} (\zcref{ex:DS}).

In \zcref{sec:abel-applications}, we apply our semiorthogonal decompositions and flop--flop formulas to spaces of linear series on curves, realized as derived Grassmannians over compactified Picard schemes.
We then investigate the interaction between spherical twists and Brill--Noether loci, showing that the twists are supported over deeper Brill--Noether loci and obtaining autoequivalences of infinite order.
For a smooth Brill--Noether--Petri general curve of genus \(g = n^{2}\) with \(n \geq 2\), we obtain a free abelian group of spherical twists of rank \(N_{n}\), where \(N_{n}\) is the Castelnuovo number.
These twists give an explicit flop--flop factorization, while a subgroup of rank \(N_{n} - 1\) acts trivially on rational cohomology and numerical \(K\)-theory (\zcref{cor:square-genus-flop-product,prop:square-genus-cohomological-action}).

In \zcref{sec:integral-grassmannian-decompositions}, we prove the incidence--flag comparison used for the flop--flop formula and categorical analogues of the extended Clifford relations of Jiang, Li, and Zhao \cite{JiangEtAl2025InstantonsAffineVertexAlgebra}.
We also identify the vacuum category, compare our semiorthogonal decompositions with those of \cite{Jiang2023DerivedCategoriesGrassmannians}, and recover the spin representation on \(K_{0}\).

\subsection{Related work}

The study of derived categories of Grassmannians goes back to Kapranov's exceptional collections \cite{Kapranov1988DerivedCategoriesHomogeneousSpaces}.
In the derived setting, the first-named author developed the theory of derived projectivizations, Grassmannians and flag schemes of perfect complexes, together with derived Schur and Weyl functors and a generalization of the Borel--Weil--Bott theorem \cite{Jiang2025DerivedProjectivizations,Jiang2023DerivedGrassmanniansSchurFunctors}.
His results include a generalized Serre theorem for derived projectivizations and semiorthogonal decompositions of derived Grassmannians, extending earlier work on Quot schemes by Toda \cite{Jiang2025DerivedProjectivizations,Jiang2023DerivedCategoriesGrassmannians,Toda2023QuotSchemesHallProducts}.
These constructions provide much of the geometric and categorical framework used in this paper.

A key ingredient in our proof is a geometric comparison of compositions of flag functors appearing in their commutation relations.
Such comparisons have been extensively studied for Hecke correspondences.
Negu{\c{t}} used quadruple moduli spaces to construct natural transformations between compositions of Hecke functors and describe their cones \cite{Negut2022HeckeCorrespondences}.
Zhao developed this approach to mixed Hecke commutators on moduli spaces of sheaves \cite{Zhao2024CategoricalQuantumToroidal,Zhao2021CategorifiedCommutator}, and subsequently established common derived blow-up descriptions for nested quiver varieties \cite{Zhao2023DerivedBlowupsNestedQuivers} and a generalized vanishing theorem for derived blow-ups \cite{Zhao2023GeneralizedVanishingBlowups}.
Our comparison of flag functors follows a similar geometric strategy, using a common derived blow-up to relate the two convolution spaces.

Geometric realizations of Clifford algebras and spin representations were studied by Savage using quiver varieties and by Xiao using the cohomological Hall algebra of the \(A_{1}\)-quiver acting on Grassmannian cohomology \cite{Savage2006GeometricSpinRepresentations,Xiao2015CohomologicalHallAlgebraA1}.
Closely related flag functors were studied by Hsu in the context of categorical actions of shifted \(q = 0\) affine algebras and their relation to semiorthogonal decompositions \cite{Hsu2025ShiftedZeroAffineAlgebra,Hsu2024SemiorthogonalDecompositionCategoricalAction}.
More recently, Jiang, Li, and Zhao constructed spin representations of finite-dimensional Clifford algebras on \(K_{0}\), Hochschild homology, and other homological theories from semiorthogonal decompositions of derived Grassmannians, with applications to affine vertex algebras \cite{JiangEtAl2025InstantonsAffineVertexAlgebra}.
Our approach is complementary: we establish explicit relations between flag functors and their adjoints at the categorical level.
These relations yield the Clifford relations on \(K_{0}\) for perfect complexes of arbitrary rank and, in positive rank, recover the semiorthogonal decompositions.

Other approaches to categorical Clifford algebras include Tian's diagrammatic categorification \cite{Tian2015DiagrammaticCliffordCategorification} and the categorical boson--fermion correspondence proposed by Cautis and Sussan \cite{CautisSussan2015CategoricalBosonFermionCorrespondence}.
In the latter framework, Gonz{\'a}lez proved the conjectured categorical Clifford relations \cite{Gonzalez2020BernsteinOperatorsBosonFermion}.
These constructions provide categorical realizations of Clifford relations through different representations and categories, complementary to the geometric flag correspondences considered here.

Derived equivalences associated with flops originate in the work of Bondal--Orlov and Bridgeland \cite{BondalOrlov1995SemiorthogonalDecomposition,Bridgeland2002FlopsDerivedCategories}.
An important phenomenon is the flop--flop = twist relation, which identifies compositions of flop equivalences with spherical or \(\mathbb{P}\)-twists, as studied by Addington, Donovan, and Meachan for standard flops and Mukai flops \cite{AddingtonEtAl2019MukaiFlopsPTwists}.
For Grassmannian flops, Donovan and Segal identified window shifts with twists of spherical functors \cite{DonovanSegal2014WindowShiftsGrassmannianTwists}, while Halpern-Leistner and Shipman showed that, more generally, the window-shift autoequivalences arising from variation of GIT quotient are spherical twists induced by mutations of semiorthogonal decompositions \cite{HalpernLeistnerShipman2016AutoequivalencesGIT}.
Ballard et al. further studied the Fourier--Mukai kernels realizing Grassmann flop equivalences \cite{BallardEtAl2021KernelsGrassmannFlops}, while Cautis, Kamnitzer, and Licata constructed related equivalences for cotangent bundles of Grassmannians via categorical \(\mathfrak{sl}_{2}\) actions \cite{CautisEtAl2013CotangentGrassmanniansDerivedEquivalences}.
Jiang and Leung established ``flop--flop = twist'' formulas for rank-zero projectivizations, with applications to symmetric powers of curves and \(\Theta\)-flops \cite{JiangLeung2022ProjectivizationsFlops}.
\zcref{thm:intro-rank-zero} provides a common generalization of the results of Jiang--Leung and Donovan--Segal; see the discussion after \zcref{thm:intro-rank-zero}.

\subsection*{Acknowledgments}
The authors would like to thank Yu Zhao for helpful discussions.
Q.J. is especially grateful to Wei-Ping Li and Yu Zhao for their collaboration on \cite{JiangEtAl2025InstantonsAffineVertexAlgebra}, through which he first came to realize the connection between Clifford algebras and semiorthogonal decompositions.
S.Z. is grateful to Yifan Jin for valuable conversations and comments on the manuscript.
This work was supported by the Hong Kong Research Grants Council ECS grant (Grant No. 26311724) and GRF grant (Grant No. 16311625).

\subsection*{Declaration of AI use}
Generative AI tools were used to assist with checking mathematical calculations, proofreading, and improving the exposition of this paper.
The authors take full responsibility for the mathematical content and any remaining errors.

\begin{samepage}
    \subsection{Conventions and notation}
    \paragraph{Categories.}
    We work with \(\infty\)-categories.
    All constructions and functors are derived.
    We write \(\operatorname{Ani}\) for the category of anima (spaces), and \(\operatorname{Ani}(\operatorname{Ring})\) for the category of animated commutative rings.
    \smallskip
    \paragraph{Geometric Setting}
    A prestack is a functor \(X \colon \operatorname{Ani}(\operatorname{Ring}) \to \operatorname{Ani}\).
    We write \(\QCoh(X)\) for the stable category of quasi-coherent complexes on \(X\), and \(\Perf(X)\) for its full subcategory of perfect complexes.
    \smallskip
    \paragraph{Indexing and Tor-amplitude.}
    We use homotopy indexing, so \(\pi_{i}(M) = H^{-i}(M)\) and \(\pi_{i}(M[n]) = \pi_{i - n}(M)\).
    A complex \(M\) over an animated commutative ring \(R\) has Tor-amplitude in \([0,1]\) if, for every discrete \(R\)-module \(N\),
    \[
        \pi_{i}(M \otimes_{R} N) = 0 \qquad \text{for }i \notin [0,1].
    \]
\end{samepage}

\smallskip

\paragraph{Bases.}
Throughout the paper, we fix a prestack \(X\) as the base.
Unless explicitly stated otherwise, \(\mathcal{E}\) denotes a perfect complex on \(X\) of Tor-amplitude in \([0,1]\) and rank \(r\); in particular, \(\mathcal{E}^{\vee}[1]\) is again a perfect complex of Tor-amplitude in \([0,1]\), of rank \(-r\).

%% file: src/flag-geometry-rees.tex
\section{Flag correspondence of derived Grassmannians}
\label{sec:flag-geometry-rees}

In this section, we study the geometry of flag correspondences on derived Grassmannians and use derived blow-ups and Rees algebras to compare their Fourier--Mukai functors.

\subsection{Derived Grassmannians and flag functors}
\label{subsec:grassmannians-flags}
We recall the basic geometry of derived Grassmannians and flag spaces from \cite{Jiang2023DerivedGrassmanniansSchurFunctors}.

\begin{definition}
    \label{def:derived-flags}
    Let \(X\) be a prestack and let \(\mathcal{E} \in \QCoh(X)^{\leq 0}\) be connective.
    For \(d \geq 0\), the derived Grassmannian of rank \(d\) quotients of \(\mathcal{E}\) over \(X\) is the relative prestack
    \[
        \begin{aligned}
            \Gr_{X}(\mathcal{E},d) \colon \operatorname{Ani}(\operatorname{Ring})_{/X} & \longrightarrow \operatorname{Ani}, \\
            (u \colon T \to X)                                                         & \longmapsto \{u^{*}\mathcal{E} \twoheadrightarrow U\},
        \end{aligned}
    \]
    where the right-hand side is the groupoid of quotient maps with \(U\) a vector bundle of rank \(d\) on \(T\), surjective on \(\pi_{0}\).
    We write
    \[
        \operatorname{pr}_{d} \colon \Gr_{X}(\mathcal{E},d) \longrightarrow X
    \]
    for the structure morphism and denote the universal fiber sequence by
    \[
        R_{d} \longrightarrow \operatorname{pr}_{d}^{*}\mathcal{E} \longrightarrow Q_{d}.
    \]

    For \(d = 1\), we also write
    \[
        \mathbb{P}_{X}(\mathcal{E}) \coloneqq \Gr_{X}(\mathcal{E},1)
    \]
    for the derived projectivization of \(\mathcal{E}\), with universal quotient line bundle \(\mathcal{O}(1)\).

    More generally, for an increasing sequence \(\mathbf{d} = (d_{1},\ldots,d_{k})\) of nonnegative integers, the derived flag prestack of type \(\mathbf{d}\) is
    \[
        \begin{aligned}
            \Flag_{X}(\mathcal{E};\mathbf{d}) \colon \operatorname{Ani}(\operatorname{Ring})_{/X} & \longrightarrow \operatorname{Ani}, \\
            (u \colon T \to X)                                                                    & \longmapsto \{u^{*}\mathcal{E} \twoheadrightarrow U_{k} \twoheadrightarrow \cdots \twoheadrightarrow U_{1}\},
        \end{aligned}
    \]
    where the right-hand side is the groupoid of chains of quotient maps, with each \(U_{i}\) a vector bundle of rank \(d_{i}\) on \(T\) and each map surjective on \(\pi_{0}\).
    We write
    \[
        \operatorname{pr}_{\mathbf{d}} \colon \Flag_{X}(\mathcal{E};\mathbf{d}) \longrightarrow X
    \]
    for the structure morphism and \(Q_{d_{i}}\) for the universal quotient bundles.
    If \(i \geq 2\) and \(d_{i} = d_{i - 1} + 1\), we write
    \[
        L_{d_{i}} \coloneqq \ker(Q_{d_{i}} \to Q_{d_{i - 1}})
    \]
    for the corresponding universal line bundle.
    When \(d_{1} = 1\), we set \(L_{1} = Q_{1}\).
\end{definition}

By abuse of notation, we use the same symbols for complexes and their pullbacks throughout the paper whenever the relevant morphism is clear.

\begin{theorem}[{\cite[Theorem~4.18 and Section~4.3]{Jiang2023DerivedGrassmanniansSchurFunctors}}]
    \label{thm:flag-geometry}
    \leavevmode
    \begin{enumerate}
        \item The flag prestack is represented by a relative derived scheme.
              For every \(u \colon X^{\prime} \to X\), there is a canonical equivalence
              \[
                  \Flag_{X}(\mathcal{E};\mathbf{d}) \times_{X} X^{\prime} \simeq \Flag_{X^{\prime}}(u^{*}\mathcal{E};\mathbf{d}).
              \]
        \item If \(\mathcal{E}\) is perfect, then \(\operatorname{pr}_{\mathbf{d}}\) is proper and locally of finite presentation.
        \item If \(\mathcal{E}\) is perfect of Tor-amplitude in \([0,a]\), with \(a \geq 0\), then the relative cotangent complex \(\mathbb{L}_{\operatorname{pr}_{\mathbf{d}}}\) is perfect of Tor-amplitude in \([0,a]\).
        \item If \(\mathcal{E}\) is perfect of Tor-amplitude in \([0,1]\), then \(\operatorname{pr}_{\mathbf{d}}\) is quasi-smooth and locally of finite Tor-amplitude.
              The pushforward \(\operatorname{pr}_{\mathbf{d},*}\) preserves perfect complexes.
    \end{enumerate}
\end{theorem}

By the definitions of the derived Grassmannian and flag space, for every \(d \geq 0\) there is a natural commutative diagram
\[
    \begin{tikzcd}
        \Flag_{X}(\mathcal{E};d,d + 1) \arrow[r,"\pi_{d}"] \arrow[d,"j_{d}"'] & \Gr_{X}(\mathcal{E},d + 1) \arrow[d,"\operatorname{pr}_{d + 1}"] \\
        \Gr_{X}(\mathcal{E},d) \arrow[r,"\operatorname{pr}_{d}"'] & X .
    \end{tikzcd}
\]

\begin{proposition}
    \label{prop:flag-projectivizations}
    The two projections give canonical equivalences
    \[
        \Flag_{X}(\mathcal{E};d,d + 1) \simeq \mathbb{P}_{\Gr_{X}(\mathcal{E},d)}(R_{d}) \simeq \mathbb{P}_{\Gr_{X}(\mathcal{E},d + 1)}(Q_{d + 1}^{\vee}),
    \]
    under which
    \[
        L_{d + 1} \simeq \mathcal{O}_{j_{d}}(1) \simeq \mathcal{O}_{\pi_{d}}(-1).
    \]
\end{proposition}

\begin{proof}
    It suffices to verify the equivalences of functors of points.
    We prove the first; the second is analogous.

    Fix a map \(T \to \Gr_{X}(\mathcal{E},d)\), corresponding to a quotient \(\mathcal{E}_{T} \twoheadrightarrow Q_{d,T}\) with fiber \(R_{d,T}\).
    A flag over this point is equivalent to a line bundle quotient \(R_{d,T} \twoheadrightarrow L\), obtained by taking the fiber of \(Q_{d + 1,T} \twoheadrightarrow Q_{d,T}\).
    Conversely, such a quotient determines the flag by the pushout
    \[
        Q_{d + 1,T} \simeq \mathcal{E}_{T} \sqcup_{R_{d,T}} L.
    \]
    These constructions are inverse and natural in \(T\), giving the first equivalence and identifying \(L_{d + 1}\) with \(\mathcal{O}_{j_{d}}(1)\).
    The second equivalence follows by dualizing the line subbundle \(L_{d + 1} \hookrightarrow Q_{d + 1}\), which identifies \(L_{d + 1}\) with \(\mathcal{O}_{\pi_{d}}(-1)\).
\end{proof}

For \(d \geq 0\) and \(m \in \mathbb{Z}\), the flag correspondence defines the Fourier--Mukai functors:
\[
    \begin{aligned}
        \Psi_{d,m} \colon \Perf(\Gr_{X}(\mathcal{E},d + 1)) & \longrightarrow \Perf(\Gr_{X}(\mathcal{E},d)), \\
        M                                               & \longmapsto j_{d,*}\bigl(L_{d + 1}^{m} \otimes \pi_{d}^{*}M\bigr).
    \end{aligned}
\]
By Grothendieck duality \cite[Theorem~4.18 and Corollary~4.47]{Jiang2023DerivedGrassmanniansSchurFunctors}, the left and right adjoints of \(\Psi_{d,m}\) are
\[
    \begin{aligned}
        \Psi_{d,m}^{L}(M) & \simeq \pi_{d,*}\bigl( L_{d + 1}^{d + 1 - m} \otimes \pi_{d}^{*}\det Q_{d + 1}^{\vee} \otimes j_{d}^{*}M \bigr)[d], \\
        \Psi_{d,m}^{R}(M) & \simeq \pi_{d,*}\bigl( L_{d + 1}^{d - r - m} \otimes j_{d}^{*}(\det R_{d} \otimes M) \bigr)[r - d - 1].
    \end{aligned}
\]
As discussed in the Introduction, we seek to compare the compositions
\[
    \Psi_{d - 1,b}^{L}\Psi_{d - 1,a}, \qquad \Psi_{d,a + 1}\Psi_{d,b + 1}^{L}
\]
as endofunctors of \(\Perf(\Gr_{X}(\mathcal{E},d))\), for \(d \geq 1\) and \(a,b \in \mathbb{Z}\).
Their Fourier--Mukai kernels are obtained by convolution along the fiber products
\[
    \begin{aligned}
        Z^{-} & = \Flag_{X}(\mathcal{E};d - 1,d) \times_{\Gr_{X}(\mathcal{E},d - 1)} \Flag_{X}(\mathcal{E};d - 1,d), \\
        Z^{+} & = \Flag_{X}(\mathcal{E};d,d + 1) \times_{\Gr_{X}(\mathcal{E},d + 1)} \Flag_{X}(\mathcal{E};d,d + 1).
    \end{aligned}
\]
Over a map \(T \to X\), the space \(Z^{-}\) parametrizes diagrams
\[
    Q_{d}^{(1)} \twoheadrightarrow Q_{d - 1} \twoheadleftarrow Q_{d}^{(2)},
\]
whereas \(Z^{+}\) parametrizes diagrams
\[
    Q_{d}^{(1)} \twoheadleftarrow Q_{d + 1} \twoheadrightarrow Q_{d}^{(2)},
\]
with all quotient maps compatible with the maps from \(\mathcal{E}_{T}\).

A natural way to relate these spaces is to retain the common quotients of ranks \(d - 1\) and \(d + 1\), giving diagrams of the form
\[
    \begin{tikzcd}[column sep=large,row sep=large]
        & & \mathcal{E}_{T} \arrow[d,two heads] & & \\
        & & Q_{d + 1} \arrow[dl,two heads] \arrow[dr,two heads] & & \\
        & Q_{d}^{(1)} \arrow[dr,two heads] & & Q_{d}^{(2)} \arrow[dl,two heads] & \\
        & & Q_{d - 1}. & &
    \end{tikzcd}
\]
These diagrams define the space
\[
    \widetilde{Z} = \Flag_{X}(\mathcal{E};d - 1,d,d + 1) \times_{\Flag_{X}(\mathcal{E};d - 1,d + 1)} \Flag_{X}(\mathcal{E};d - 1,d,d + 1).
\]
Forgetting the rank \(d + 1\) or rank \(d - 1\) quotient gives natural morphisms
\[
    Z^{-} \xleftarrow{b_{-}} \widetilde{Z} \xrightarrow{b_{+}} Z^{+}.
\]
We will use this correspondence to compare the two Fourier--Mukai kernels.

To understand the geometry of these forgetful maps, we first describe the underlying construction in a more explicit form.
Let \(X\) be a prestack and let \(\mathcal{E} \in \QCoh(X)^{\leq 0}\) be connective.
Put
\[
    P = \Gr_{X}(\mathcal{E},1), \qquad G = \Gr_{X}(\mathcal{E},2), \qquad F = \Flag_{X}(\mathcal{E};1,2).
\]
Forgetting the common quotient of rank two defines a natural morphism
\[
    b \colon F \times_{G} F \longrightarrow P \times_{X} P.
\]
A \(T\)-point of \(P \times_{X} P\) consists of two line quotients
\[
    \mathcal{E}_{T} \twoheadrightarrow L_{1}^{(1)}, \qquad \mathcal{E}_{T} \twoheadrightarrow L_{1}^{(2)}.
\]
Away from the relative diagonal \(\Delta_{P/X}\), the induced map
\[
    \mathcal{E}_{T} \longrightarrow L_{1}^{(1)} \oplus L_{1}^{(2)}
\]
is a quotient of rank two and determines a unique lift to \(F \times_{G} F\).
Thus \(b\) is an equivalence away from the diagonal.

Moreover, the relative diagonals fit into a commutative square
\[
    \begin{tikzcd}[column sep=large,row sep=large]
        F \arrow[r,hook,"\Delta_{F/G}"] \arrow[d,"j_{1}"']
        & F \times_{G} F \arrow[d,"b"] \\
        P \arrow[r,hook,"\Delta_{P/X}"']
        & P \times_{X} P,
    \end{tikzcd}
\]
which is Cartesian on classical truncations.
Over the diagonal, the map \(b\) records the additional choice of a quotient of rank two extending the given line quotient.
This suggests that \(b\) is the blow-up of \(P \times_{X} P\) along its diagonal.

\subsection{Blow-up geometry of flag correspondences}

Khan and Rydh introduced a construction of derived blow-ups of quasi-smooth closed immersions via virtual Cartier divisors \cite[\S4]{KhanRydh2025VirtualCartierDivisorsBlowups}.
This moduli-theoretic description was subsequently extended to arbitrary closed immersions by Hekking, Khan, and Rydh \cite[Definition~4.3]{HekkingEtAl2025BlowupsDerivedWeilRestrictions}.
We recall their definition and use it to prove that the morphism \(b\) introduced above is the blow-up of \(P \times_{X}P\) along its relative diagonal.
We then apply this result to the two flag correspondences, identifying \(\widetilde{Z}\) with the blow-ups of \(Z^{-}\) and \(Z^{+}\) along their respective diagonals.
This is analogous to the common diagonal blow-ups appearing in Zhao's study of nested quiver varieties \cite[Theorems~1.1 and~6.3]{Zhao2023DerivedBlowupsNestedQuivers}.

\begin{definition}
    \label{def:derived-blowup}
    Let \(i \colon Z \hookrightarrow M\) be a closed immersion.
    The derived blow-up \(b \colon \operatorname{Bl}_{Z}(M) \to M\) is the relative derived scheme representing the functor that assigns to each \(f \colon T \to M\) the space of commutative squares
    \[
        \begin{tikzcd}
            D \arrow[r,hook,"j"] \arrow[d,"g"'] & T \arrow[d,"f"] \\
            Z \arrow[r,hook,"i"'] & M
        \end{tikzcd}
    \]
    such that \(j\) is a virtual Cartier divisor, i.e., a quasi-smooth closed immersion of virtual codimension one, the square is Cartesian on classical truncations, and the induced morphism
    \[
        g^{*}\mathbb{L}_{Z/M}[-1] \longrightarrow \mathbb{L}_{D/T}[-1]
    \]
    is surjective on \(\pi_{0}\).
\end{definition}

\begin{lemma}
    \label{lem:rank-one-flag-blowup}
    Let \(X\) be a prestack and let \(\mathcal{E} \in \QCoh(X)^{\leq 0}\) be connective.
    Put
    \[
        P = \Gr_{X}(\mathcal{E},1), \qquad G = \Gr_{X}(\mathcal{E},2), \qquad F = \Flag_{X}(\mathcal{E};1,2).
    \]
    Then the forgetful morphism
    \[
        b \colon F \times_{G} F \longrightarrow P \times_{X} P
    \]
    identifies \(F \times_{G} F\) with the derived blow-up of \(P \times_{X} P\) along its relative diagonal:
    \[
        F \times_{G} F \simeq \operatorname{Bl}_{\Delta_{P/X}}(P \times_{X} P).
    \]
    The exceptional divisor is the relative diagonal
    \[
        \Delta_{F/G} \colon F \hookrightarrow F \times_{G} F.
    \]
    Writing \(L_{k}^{(i)}\) for the pullback of the universal line bundle \(L_{k}\) along the \(i\)-th projection \(F \times_{G} F \to F\), its ideal sheaf is
    \[
        \mathcal{I}_{\Delta_{F/G}} \simeq L_{2}^{(1)} \otimes (L_{1}^{(2)})^{-1} \simeq L_{2}^{(2)} \otimes (L_{1}^{(1)})^{-1}.
    \]
\end{lemma}

\begin{proof}
    Fix \(T \to P \times_{X} P\), corresponding to two line quotients
    \[
        q \colon \mathcal{E}_{T} \twoheadrightarrow L, \qquad q^{\prime} \colon \mathcal{E}_{T} \twoheadrightarrow L^{\prime},
    \]
    and put \(K = \fib(q)\).
    By the definition of \(F \times_{G} F\), a lift of \((q,q^{\prime})\) along \(b\) amounts to a common quotient \(\mathcal{E}_{T} \twoheadrightarrow Q\) of rank two, together with quotient maps \(Q \twoheadrightarrow L\) and \(Q \twoheadrightarrow L^{\prime}\) whose composites are \(q\) and \(q^{\prime}\), respectively.

    Given such a lift, let \(H = \fib(Q \to L)\).
    Restricting the second quotient \(Q \twoheadrightarrow L^{\prime}\) to \(H\) gives a morphism \(s \colon H \to Q \to L^{\prime}\), or equivalently a section of \(H^{\vee} \otimes L^{\prime}\).
    Let \(j \colon D \hookrightarrow T\) be the derived zero locus of \(s\).
    Then \(j\) is quasi-smooth of relative virtual dimension \(-1\), with ideal sheaf
    \[
        \mathcal{I}_{D/T} \simeq H \otimes (L^{\prime})^{-1}.
    \]
    Over \(D\), the quotient \(Q \to L^{\prime}\) factors through \(Q/H \simeq L\), inducing an isomorphism \(\varphi \colon L|_{D} \simeq L^{\prime}|_{D}\) compatible with \(q\) and \(q^{\prime}\).
    Moreover, \(q^{\prime}|_{D} \simeq \varphi \circ q|_{D}\), so the two quotient maps define a morphism \(g \colon D \to P\) and a commutative square
    \[
        \begin{tikzcd}
            D \arrow[r,hook,"j"] \arrow[d,"g"']
            & T \arrow[d,"{(q,q^{\prime})}"] \\
            P \arrow[r,hook,"\Delta_{P/X}"']
            & P \times_{X} P.
        \end{tikzcd}
    \]
    This square is Cartesian on classical truncations.
    Indeed, the two line quotients agree classically precisely where the restriction of \(Q \to L^{\prime}\) to \(H\) vanishes, which is the classical zero locus of \(s\).

    It remains to verify the condition on cotangent complexes in \zcref{def:derived-blowup}.
    The square above induces a canonical map
    \[
        g^{*}\mathbb{L}_{P/(P \times_{X} P)}[-1] \longrightarrow \mathbb{L}_{D/T}[-1].
    \]
    For the source, the cotangent complex formula for the relative diagonal and the derived Grassmannian \cite[Theorem~4.18]{Jiang2023DerivedGrassmanniansSchurFunctors} gives
    \[
        g^{*}\mathbb{L}_{P/(P \times_{X} P)}[-1] \simeq g^{*}\mathbb{L}_{P/X} \simeq (K \otimes L^{-1})|_{D}.
    \]
    On the other hand, since \(D\) is the derived zero locus of \(s \colon H \to L^{\prime}\), its relative cotangent complex satisfies
    \[
        \mathbb{L}_{D/T}[-1] \simeq \bigl(H \otimes (L^{\prime})^{-1}\bigr)|_{D}.
    \]
    Under these identifications, the canonical map becomes
    \[
        (K \otimes L^{-1})|_{D} \longrightarrow \bigl(H \otimes (L^{\prime})^{-1}\bigr)|_{D}.
    \]
    Using \(L|_{D} \simeq L^{\prime}|_{D}\), this map is induced by the natural morphism \(K \to H\) in the commutative diagram
    \[
        \begin{tikzcd}
            K \arrow[r] \arrow[d]
            & \mathcal{E}_{T} \arrow[r] \arrow[d]
            & L \arrow[d,equal] \\
            H \arrow[r]
            & Q \arrow[r]
            & L.
        \end{tikzcd}
    \]
    Since \(\mathcal{E}_{T} \twoheadrightarrow Q\) is a quotient, \(K \to H\) is surjective on \(\pi_{0}\).
    Hence the conormal condition is satisfied, and the construction defines a \(T\)-point of \(\operatorname{Bl}_{\Delta_{P/X}}(P \times_{X} P)\).

    Conversely, let
    \[
        \begin{tikzcd}
            D \arrow[r,hook,"j"] \arrow[d,"g"']
            & T \arrow[d,"{(q,q^{\prime})}"] \\
            P \arrow[r,hook,"\Delta_{P/X}"']
            & P \times_{X} P
        \end{tikzcd}
    \]
    be a square satisfying the conditions of \zcref{def:derived-blowup}.
    As before, write
    \[
        q \colon \mathcal{E}_{T} \twoheadrightarrow L, \qquad q^{\prime} \colon \mathcal{E}_{T} \twoheadrightarrow L^{\prime}, \qquad K = \fib(q).
    \]
    Since \(j \colon D \hookrightarrow T\) is a derived Cartier divisor, its ideal \(\ell = \fib(\mathcal{O}_{T} \to j_{*}\mathcal{O}_{D})\) is a line bundle.

    The morphism \(g \colon D \to P\) identifies the two quotient maps over \(D\).
    Hence the composite \(K|_{D} \to \mathcal{E}_{T}|_{D} \xrightarrow{q^{\prime}|_{D}} L^{\prime}|_{D}\) is null-homotopic.
    On the other hand, tensoring the defining fiber sequence of \(\ell\) with \(L^{\prime}\) gives
    \[
        \ell \otimes L^{\prime} \longrightarrow L^{\prime} \longrightarrow L^{\prime} \otimes j_{*}\mathcal{O}_{D}.
    \]
    By the projection formula and adjunction, the composite \(K \to L^{\prime} \to L^{\prime} \otimes j_{*}\mathcal{O}_{D}\) is therefore null-homotopic.
    We thus obtain a factorization
    \[
        K \xrightarrow{h} \ell \otimes L^{\prime} \longrightarrow L^{\prime}.
    \]

    We claim that \(h\) is surjective on \(\pi_{0}\).
    Along \(D\), the conormal map of the given square is
    \[
        g^{*}\mathbb{L}_{P/(P \times_{X} P)}[-1] \simeq (K \otimes L^{-1})|_{D} \longrightarrow \mathbb{L}_{D/T}[-1] \simeq \ell|_{D}.
    \]
    Under the identification \(L|_{D} \simeq L^{\prime}|_{D}\), this map is precisely the restriction of \(h \otimes L^{-1}\).
    The conormal condition in \zcref{def:derived-blowup} therefore implies that \(h\) is surjective on \(\pi_{0}\) along \(D\).

    Away from \(D\), we have \(\ell \simeq \mathcal{O}_{T}\), so \(h\) identifies with \(K \to L^{\prime}\).
    If this map were not surjective at a geometric point, the quotient \(q^{\prime}\) would factor through \(q\), and the two line quotients would coincide at that point.
    This contradicts the classical Cartesian condition, which identifies their coincidence locus with \(D\).
    Hence \(h\) is surjective on \(\pi_{0}\) everywhere.

    We now construct the desired quotient of rank two.
    Define \(Q\) by the cocartesian square
    \[
        \begin{tikzcd}[column sep=large,row sep=large]
            K \arrow[r] \arrow[d,"h"']
            & \mathcal{E}_{T} \arrow[d] \\
            \ell \otimes L^{\prime} \arrow[r]
            & Q.
        \end{tikzcd}
    \]
    Since the square is cocartesian, it induces a fiber sequence
    \[
        \ell \otimes L^{\prime} \longrightarrow Q \longrightarrow L.
    \]
    As both ends are line bundles, \(Q\) is a vector bundle of rank two.
    Moreover, the pushout identifies
    \[
        \cofib(\mathcal{E}_{T} \to Q) \simeq \cofib(h).
    \]
    Since \(h\) is surjective on \(\pi_{0}\), so is \(\mathcal{E}_{T} \to Q\).
    Thus we obtain a quotient of rank two.

    The maps \(\mathcal{E}_{T} \to L^{\prime}\) and \(\ell \otimes L^{\prime} \to L^{\prime}\) agree on \(K\), so the pushout induces a map \(Q \to L^{\prime}\).
    This map is surjective on \(\pi_{0}\), since its composite with \(\mathcal{E}_{T} \twoheadrightarrow Q\) is \(q^{\prime}\).
    Together with \(Q \twoheadrightarrow L\), it gives two compatible quotient maps
    \[
        Q \twoheadrightarrow L, \qquad Q \twoheadrightarrow L^{\prime},
    \]
    whose composites with \(\mathcal{E}_{T} \twoheadrightarrow Q\) recover \(q\) and \(q^{\prime}\).
    Hence we obtain a lift \(T \to F \times_{G} F\).

    Finally, we verify that the two constructions are inverse.
    Starting with a common quotient \(\mathcal{E}_{T} \twoheadrightarrow Q\) of rank two, let \(H = \fib(Q \to L)\).
    The derived zero locus of \(H \to L^{\prime}\) has ideal
    \[
        \ell \simeq H \otimes (L^{\prime})^{-1}.
    \]
    Under this identification, the map \(h \colon K \to \ell \otimes L^{\prime}\) is the natural map \(K \to H\), so the pushout recovers the original \(Q\).

    Conversely, starting with the blow-up data, the construction gives a fiber sequence
    \[
        \ell \otimes L^{\prime} \longrightarrow Q \longrightarrow L,
    \]
    and the composite
    \[
        \ell \otimes L^{\prime} \longrightarrow Q \longrightarrow L^{\prime}
    \]
    is the canonical map defining the divisor \(D\).
    Its derived zero locus is therefore \(D\), and the resulting identification of the two quotient maps over \(D\) agrees with the original one.
    The constructions are inverse and natural in \(T\), completing the proof.
\end{proof}

\begin{remark}
    The proof of \zcref{lem:rank-one-flag-blowup} also applies to pairs of quotients of ranks \(1\) and \(n\), for any \(n \geq 2\).
    More precisely, for a connective complex \(\mathcal{E}\) over \(X\), the natural forgetful morphism
    \[
        \begin{aligned}
             & \Flag_{X}(\mathcal{E};1,n + 1) \times_{\Gr_{X}(\mathcal{E},n + 1)} \Flag_{X}(\mathcal{E};n,n + 1) \\
             & \qquad \longrightarrow \Gr_{X}(\mathcal{E},1) \times_{X} \Gr_{X}(\mathcal{E},n)
        \end{aligned}
    \]
    identifies its source with the derived blow-up along the incidence correspondence
    \[
        \Flag_{X}(\mathcal{E};1,n) \hookrightarrow \Gr_{X}(\mathcal{E},1) \times_{X} \Gr_{X}(\mathcal{E},n).
    \]
    The exceptional divisor is naturally identified with
    \[
        \Flag_{X}(\mathcal{E};1,n,n + 1).
    \]
    The same construction can also be used to obtain analogous Clifford-type relations for the associated flag operators, which we do not pursue in this paper.
\end{remark}

We now apply \zcref{lem:rank-one-flag-blowup} to the two flag correspondences considered in \zcref{subsec:grassmannians-flags}.

\begin{lemma}
    \label{lem:common-flag-blowup}
    Let \(X\) be a prestack, let \(\mathcal{E} \in \QCoh(X)^{\leq 0}\) be connective, and let \(d \geq 1\).
    Put
    \[
        F_{-} = \Flag_{X}(\mathcal{E};d - 1,d), \qquad F_{+} = \Flag_{X}(\mathcal{E};d,d + 1),
    \]
    and
    \[
        \begin{aligned}
            Z^{-} & = F_{-} \times_{\Gr_{X}(\mathcal{E},d - 1)} F_{-}, \\
            Z^{+} & = F_{+} \times_{\Gr_{X}(\mathcal{E},d + 1)} F_{+}.
        \end{aligned}
    \]
    Let
    \[
        \widetilde{Z} = \Flag_{X}(\mathcal{E};d - 1,d,d + 1) \times_{\Flag_{X}(\mathcal{E};d - 1,d + 1)} \Flag_{X}(\mathcal{E};d - 1,d,d + 1),
    \]
    and let
    \[
        b_{-} \colon \widetilde{Z} \longrightarrow Z^{-}, \qquad b_{+} \colon \widetilde{Z} \longrightarrow Z^{+}
    \]
    be the morphisms forgetting the quotient of rank \(d + 1\) and rank \(d - 1\), respectively.

    Then \(b_{-}\) and \(b_{+}\) identify \(\widetilde{Z}\) with the derived blow-ups of \(Z^{-}\) and \(Z^{+}\) along their respective diagonals
    \[
        i_{-} \colon F_{-} \hookrightarrow Z^{-}, \qquad i_{+} \colon F_{+} \hookrightarrow Z^{+}.
    \]
    Under both identifications, the exceptional divisor is the relative diagonal
    \[
        D \simeq \Flag_{X}(\mathcal{E};d - 1,d,d + 1) \hookrightarrow \widetilde{Z}.
    \]
\end{lemma}

\begin{proof}
    By \zcref{prop:flag-projectivizations}, we have
    \[
        \begin{aligned}
            F_{-} & \simeq \mathbb{P}_{\Gr_{X}(\mathcal{E},d - 1)}(R_{d - 1}), \\
            F_{+} & \simeq \mathbb{P}_{\Gr_{X}(\mathcal{E},d + 1)}(Q_{d + 1}^{\vee}).
        \end{aligned}
    \]
    For the first identification, a common quotient of rank two of \(R_{d - 1}\) extends, by pushout, the fixed quotient of rank \(d - 1\) to a quotient of rank \(d + 1\) of \(\mathcal{E}\).
    Thus the forgetful morphism in \zcref{lem:rank-one-flag-blowup}, applied to \(R_{d - 1}\), is identified with \(b_{-}\).

    Dually, a common quotient of rank two of \(Q_{d + 1}^{\vee}\) corresponds to a subbundle of rank two of \(Q_{d + 1}\) containing the two kernel lines.
    Quotienting by this subbundle gives the common quotient of rank \(d - 1\), identifying the corresponding forgetful morphism with \(b_{+}\).

    In both cases, the relative diagonal corresponds to \(D\), which is therefore the exceptional divisor.
\end{proof}

\subsection{Comparison of compositions with left adjoints}
\label{subsec:rees-calculation}

By \zcref{lem:common-flag-blowup}, the two fiber products \(Z^{-}\) and \(Z^{+}\) admit the common blow-up \(\widetilde{Z}\).
Derived blow-ups also admit a description in terms of Rees algebras and derived Proj, developed by Hekking \cite[\S\S6--7]{Hekking2021GradedAlgebrasProjectiveSpectraBlowups} and identified with the moduli-theoretic construction via virtual Cartier divisors by Hekking--Khan--Rydh \cite[Theorem~4.7]{HekkingEtAl2025BlowupsDerivedWeilRestrictions}.
We use this description to compare the corresponding Fourier--Mukai kernels, and first recall the relevant properties of the extended Rees algebra.

Let \(i \colon Z \hookrightarrow M\) be a closed immersion, and let
\[
    \mathcal{R}^{\mathrm{ext}}_{Z/M} = \bigoplus_{n \in \mathbb{Z}}\mathcal{R}_{n}, \qquad \deg(\tau^{-1}) = -1,
\]
be its extended Rees algebra, a graded \(\mathcal{O}_{M}[\tau^{-1}]\)-algebra \cite[Definition~3.20]{HekkingEtAl2025BlowupsDerivedWeilRestrictions}.
Its nonnegative part
\[
    \mathcal{R}_{Z/M} \coloneqq \bigoplus_{n \geq 0}\mathcal{R}_{n}
\]
is the Rees algebra of \(i\).

By \cite[Proposition~3.21 and Corollary~3.22]{HekkingEtAl2025BlowupsDerivedWeilRestrictions},
\[
    \mathcal{R}_{n} \simeq \mathcal{O}_{M} \qquad (n \leq 0), \qquad \mathcal{R}_{1} \simeq \mathcal{I}_{Z/M} \coloneqq \fib(\mathcal{O}_{M} \to i_{*}\mathcal{O}_{Z}),
\]
and multiplication by \(\tau^{-1}\) gives the filtration
\[
    \cdots \longrightarrow \mathcal{R}_{n + 1} \longrightarrow \mathcal{R}_{n} \longrightarrow \cdots \longrightarrow \mathcal{R}_{0} = \mathcal{O}_{M}.
\]

The following result identifies the blow-up with the derived Proj of the Rees algebra and describes its exceptional divisor.
Here, \(\operatorname{Proj}_{M}(\mathcal{A})\) denotes the derived projective spectrum of a nonnegatively graded algebra \(\mathcal{A}\) \cite{Hekking2021GradedAlgebrasProjectiveSpectraBlowups}.

\begin{proposition}[{\cite[Theorem~4.7 and Equation~(4.7)]{HekkingEtAl2025BlowupsDerivedWeilRestrictions}}]
    \label{prop:rees-blowup-proj}
    Let \(i \colon Z \hookrightarrow M\) be a closed immersion, and put \(\mathcal{N} = \mathbb{L}_{Z/M}[-1]\).
    Let \(b \colon \operatorname{Bl}_{Z}(M) \to M\) be the derived blow-up, with exceptional divisor \(D\).
    There is a canonical equivalence
    \[
        \operatorname{Bl}_{Z}(M) \simeq \operatorname{Proj}_{M}(\mathcal{R}_{Z/M}),
    \]
    under which the tautological line bundle \(\mathcal{O}(1)\) is identified with \(\mathcal{O}(-D)\).
    Moreover, the exceptional divisor is canonically identified with
    \[
        D \simeq \mathbb{P}_{Z}(\mathcal{N}),
    \]
    under which
    \[
        \mathcal{O}(-D)|_{D} \simeq \mathcal{O}_{\mathbb{P}_{Z}(\mathcal{N})}(1).
    \]
\end{proposition}

In the following, write \(\widetilde{M} = \operatorname{Bl}_{Z}(M)\), with blow-up map \(b\) and exceptional divisor \(D\).

For each \(n \in \mathbb{Z}\), we construct two canonical maps
\[
    \mathcal{R}_{n} \longrightarrow \mathcal{O}_{M}, \qquad \theta_{n} \colon \mathcal{R}_{n} \longrightarrow b_{*}\mathcal{O}_{\widetilde{M}}(-nD).
\]

The first map comes from the Rees filtration.
For \(n \geq 0\), it is given by multiplication by \(\tau^{-n}\), while for \(n \leq 0\), it is the canonical equivalence \(\mathcal{R}_{n} \simeq \mathcal{O}_{M}\).

To construct the second map, recall that \zcref{prop:rees-blowup-proj} identifies \(\widetilde{M}\) with \(\operatorname{Proj}_{M}(\mathcal{R}_{Z/M})\).
The derived Proj construction provides a canonical morphism of graded algebras
\[
    b^{*}\mathcal{R}_{Z/M} \longrightarrow \bigoplus_{n \geq 0}\mathcal{O}_{\widetilde{M}}(n) \simeq \bigoplus_{n \geq 0}\mathcal{O}_{\widetilde{M}}(-nD).
\]
Taking the degree \(n\) component, for \(n \geq 0\), gives
\[
    b^{*}\mathcal{R}_{n} \longrightarrow \mathcal{O}_{\widetilde{M}}(-nD),
\]
whose adjoint defines \(\theta_{n}\).
For \(n < 0\), using \(\mathcal{R}_{n} \simeq \mathcal{O}_{M}\), we instead define \(\theta_{n}\) by adjunction from the natural map
\[
    \mathcal{O}_{\widetilde{M}} \longrightarrow \mathcal{O}_{\widetilde{M}}(-nD)
\]
associated with the effective divisor \((-n)D\).
The two descriptions agree at \(n = 0\).

We first compute the cofiber of \(\mathcal{R}_{n} \to \mathcal{O}_{M}\).
By \cite[Corollary~3.22]{HekkingEtAl2025BlowupsDerivedWeilRestrictions}, the extended Rees algebra satisfies
\[
    \mathcal{R}^{\mathrm{ext}}_{Z/M}/(\tau^{-1}) \simeq i_{*}\Sym^{*}\bigl(\mathbb{L}_{Z/M}[-1]\bigr).
\]
Taking the degree \(n\) component gives a canonical fiber sequence
\[
    \mathcal{R}_{n + 1} \longrightarrow \mathcal{R}_{n} \longrightarrow i_{*}\Sym^{n}\bigl(\mathbb{L}_{Z/M}[-1]\bigr), \qquad n \geq 0.
\]
Consequently, the cofiber of \(\mathcal{R}_{n} \to \mathcal{O}_{M}\) vanishes for \(n \leq 0\).
For \(n > 0\), it admits a finite filtration with graded pieces
\begin{equation}
    \label{eq:rees-augmentation-cofiber-filtration}
    \gr_{t}\cofib(\mathcal{R}_{n} \to \mathcal{O}_{M}) \simeq i_{*}\Sym^{t}\bigl(\mathbb{L}_{Z/M}[-1]\bigr), \qquad 0 \leq t < n,
\end{equation}
attached in the order \(t = n - 1,\ldots,0\).

We compute the cofiber of \(\theta_{n}\) by reducing locally to the blow-up of a zero section, where the generalized Serre theorem \cite[Theorem~5.6]{Jiang2025DerivedProjectivizations} gives the required calculation.
This follows the approach of Zhao \cite[Theorem~5.6]{Zhao2023GeneralizedVanishingBlowups}, who studied comparison maps between Rees algebras and pushforwards of line bundles for derived blow-ups of quasi-smooth derived Artin stacks over \(\mathbb{C}\).
Our argument carries out an analogous comparison for closed immersions admitting local zero-section models, and applies in particular to the diagonal blow-ups arising from flag correspondences over arbitrary bases.

\begin{lemma}
    \label{lem:zero-section-rees-comparison}
    Let \(i \colon Z \hookrightarrow M\) be a closed immersion equipped with a retraction \(p \colon M \to Z\), and let
    \[
        b \colon \widetilde{M} = \operatorname{Bl}_{Z}(M) \longrightarrow M
    \]
    be its derived blow-up, with exceptional divisor \(D\).
    Put
    \[
        \mathcal{N} = \mathbb{L}_{Z/M}[-1].
    \]
    Assume that \(\mathcal{N}\) is perfect of Tor-amplitude in \([0,1]\) and rank \(\rho\).
    Suppose that, for every map \(u \colon T \to Z\) with \(T\) affine, the induced section of \(M \times_{Z} T \to T\) admits, Zariski locally on \(T\), an open neighborhood identified over \(T\) with the zero section
    \[
        T \longrightarrow \Spec_{T}\Sym(u^{*}\mathcal{N}).
    \]

    Then \(\theta_{n}\) is an equivalence for \(n > -\rho\).
    For \(n \leq -\rho\), its cofiber admits a finite filtration with graded pieces
    \begin{equation}
        \label{eq:rees-tautological-cofiber-filtration}
        \gr_{t}\cofib(\theta_{n}) \simeq i_{*}\bigl( (\Sym^{t}\mathcal{N})^{\vee} \otimes \det \mathcal{N}^{\vee} \bigr)[1 - \rho], \qquad 0 \leq t \leq -\rho - n,
    \end{equation}
    attached in the order \(t = 0,1,\ldots,-\rho - n\).
\end{lemma}

\begin{proof}
    We first consider the zero section of
    \[
        M = \Spec_{Z}\Sym(\mathcal{N}).
    \]
    By the universal property of the extended Rees algebra, its pushforward along \(p\) is freely generated over \(\mathcal{O}_{Z}[\tau^{-1}]\) by \(\mathcal{N}\) in degree one, with the degree-zero generators of \(\Sym(\mathcal{N})\) given by multiplication by \(\tau^{-1}\).
    Consequently,
    \[
        p_{*}\mathcal{R}_{n} \simeq \bigoplus_{m \geq n}\Sym^{m}\mathcal{N},
    \]
    where \(\Sym^{m}\mathcal{N} = 0\) for \(m < 0\).

    The projective spectrum of the nonnegative part parametrizes a line quotient \(\mathcal{N}_{T} \twoheadrightarrow L\) together with a map \(L \to \mathcal{O}_{T}\).
    Thus \zcref{prop:rees-blowup-proj} gives
    \[
        \widetilde{M} \simeq \Spec_{\mathbb{P}_{Z}(\mathcal{N})} \left(\bigoplus_{m \geq 0}\mathcal{O}(m)\right).
    \]
    This is the total space of \(\mathcal{O}(-1)\), with \(D\) as its zero section.
    Writing \(e \colon \mathbb{P}_{Z}(\mathcal{N}) \to Z\), we obtain
    \[
        p_{*}b_{*}\mathcal{O}_{\widetilde{M}}(-nD) \simeq \bigoplus_{m \geq n}e_{*}\mathcal{O}_{\mathbb{P}_{Z}(\mathcal{N})}(m).
    \]
    Under these identifications, \(p_{*}\theta_{n}\) is the direct sum of the tautological maps
    \[
        \phi_{m} \colon \Sym^{m}\mathcal{N} \longrightarrow e_{*}\mathcal{O}_{\mathbb{P}_{Z}(\mathcal{N})}(m), \qquad m \geq n.
    \]

    By the generalized Serre theorem \cite[Theorem~5.6]{Jiang2025DerivedProjectivizations}, the maps \(\phi_{m}\) are equivalences for \(m > -\rho\).
    Since \(p\) is affine, \(p_{*}\) is conservative, so \(\theta_{n}\) is an equivalence for \(n > -\rho\).

    For a general immersion, let \(u \colon T \to Z\) be a map with \(T\) affine, and write \(v \colon M_{T} = M \times_{Z} T \to M\).
    Base change gives
    \[
        v^{*}\mathcal{R}_{n} \simeq \mathcal{R}_{n,T}, \qquad v^{*}b_{*}\mathcal{O}_{\widetilde{M}}(-nD) \simeq b_{T,*}\mathcal{O}_{\widetilde{M}_{T}}(-nD_{T}),
    \]
    compatibly with \(\theta_{n}\).
    For \(n > -\rho\), the zero-section calculation therefore shows that \(\theta_{n}\) is an equivalence Zariski locally near \(T\), while away from \(T\) it is the identity.
    Hence \(\theta_{n}\) is an equivalence in general.

    We now construct the filtration of its cofiber.
    By \zcref{prop:rees-blowup-proj}, the exceptional divisor is
    \[
        D \simeq \mathbb{P}_{Z}(\mathcal{N}), \qquad \mathcal{O}(-D)|_{D} \simeq \mathcal{O}_{D}(1).
    \]
    Write \(e \colon D \to Z\) for the exceptional projection.
    The Rees filtration and restriction to \(D\) give a canonical commutative diagram of fiber sequences
    \[
        \begin{tikzcd}[column sep=large,row sep=large]
            \mathcal{R}_{m + 1}
            \arrow[r] \arrow[d,"\theta_{m + 1}"']
            & \mathcal{R}_{m}
            \arrow[r] \arrow[d,"\theta_{m}"]
            & i_{*}\Sym^{m}\mathcal{N}
            \arrow[d,"i_{*}\phi_{m}"] \\
            b_{*}\mathcal{O}_{\widetilde{M}}(-(m + 1)D)
            \arrow[r]
            & b_{*}\mathcal{O}_{\widetilde{M}}(-mD)
            \arrow[r]
            & i_{*}e_{*}\mathcal{O}_{D}(m).
        \end{tikzcd}
    \]
    Here the right vertical map is induced by the universal quotient \(e^{*}\mathcal{N} \to \mathcal{O}_{D}(1)\).

    The generalized Serre theorem identifies its cofiber as
    \[
        \cofib(i_{*}\phi_{m}) \simeq i_{*}\bigl( (\Sym^{-\rho - m}\mathcal{N})^{\vee} \otimes \det \mathcal{N}^{\vee} \bigr)[1 - \rho]
    \]
    for \(m \leq -\rho\), and shows that it vanishes for \(m > -\rho\).
    Taking the cofibers of the vertical maps therefore gives canonical fiber sequences
    \begin{equation}
        \label{eq:rees-comparison-cofiber-recursion}
        \cofib(\theta_{m + 1}) \longrightarrow \cofib(\theta_{m}) \longrightarrow i_{*}\bigl( (\Sym^{-\rho - m}\mathcal{N})^{\vee} \otimes \det \mathcal{N}^{\vee} \bigr)[1 - \rho].
    \end{equation}

    Since \(\cofib(\theta_{1 - \rho}) = 0\), these sequences give a finite filtration of \(\cofib(\theta_{n})\) for \(n \leq -\rho\).
    Reindexing by \(t = -\rho - m\) gives the graded pieces in \eqref{eq:rees-tautological-cofiber-filtration}, attached in the order \(t = 0,1,\ldots,-\rho - n\).
\end{proof}

We now apply \zcref{lem:zero-section-rees-comparison} to the two diagonal immersions
\[
    i_{-} \colon F_{-} \hookrightarrow Z^{-}, \qquad i_{+} \colon F_{+} \hookrightarrow Z^{+}.
\]
Our goal is to compare the Fourier--Mukai kernels through their common blow-up \(\widetilde{Z}\).
For \(i_{+}\), we need to show that \(\mathcal{O}_{Z^{+}} \to b_{+,*}\mathcal{O}_{\widetilde{Z}}\) is an equivalence.
For \(i_{-}\), we will use the Rees comparison maps to compute their cofibers.
Both calculations require the conormal complexes of the diagonals and the local zero-section models in \zcref{lem:zero-section-rees-comparison}.

By \zcref{prop:flag-projectivizations}, we have
\[
    \begin{aligned}
        F_{-} & \simeq \mathbb{P}_{\Gr_{X}(\mathcal{E},d - 1)}(R_{d - 1}), \\
        F_{+} & \simeq \mathbb{P}_{\Gr_{X}(\mathcal{E},d + 1)}(Q_{d + 1}^{\vee}).
    \end{aligned}
\]
Since \(Z^{-}\) and \(Z^{+}\) are the corresponding relative fiber products, each diagonal becomes a section of a derived projectivization after base change along the first projection.
We therefore recall an affine chart around a given line quotient.

Let \(E \twoheadrightarrow L\) be a line bundle quotient over a derived scheme \(T\), with fiber \(K\).
Zariski locally on \(T\), the quotient splits, so that \(E \simeq K \oplus L\).
Consider the open subspace of \(\mathbb{P}_{T}(E)\) consisting of line quotients \(E \twoheadrightarrow L^{\prime}\) whose restriction to \(L\) is an isomorphism.
Identifying \(L^{\prime}\) with \(L\), such a quotient is uniquely determined by a map \(K \to L\).
Hence this open subspace is
\begin{equation}
    \label{eq:projectivization-affine-chart}
    \Spec_{T}\Sym(K \otimes L^{-1}),
\end{equation}
with the original quotient corresponding to the zero section.

\begin{lemma}
    \label{lem:lower-diagonal-local-model}
    Assume that \(\mathcal{E}\) is perfect of Tor-amplitude in \([0,1]\) and rank \(r\).
    The conormal complexes of the two diagonals are
    \[
        \begin{aligned}
            N_{-} & \coloneqq \mathbb{L}_{F_{-}/Z^{-}}[-1] \simeq \pi_{d - 1}^{*}R_{d} \otimes L_{d}^{-1}, \\
            N_{+} & \coloneqq \mathbb{L}_{F_{+}/Z^{+}}[-1] \simeq j_{d}^{*}Q_{d}^{\vee} \otimes L_{d + 1}.
        \end{aligned}
    \]
    Here \(N_{-}\) is perfect of Tor-amplitude in \([0,1]\) and rank \(r - d\), while \(N_{+}\) is a vector bundle of rank \(d\).

    Regard each diagonal as a section of the first projection \(Z^{\pm} \to F_{\pm}\).
    For any map \(u \colon T \to F_{\pm}\) with \(T\) affine, the induced section admits, Zariski locally on \(T\), an open neighborhood identified with the zero section
    \[
        T \longrightarrow \Spec_{T}\Sym(u^{*}N_{\pm}).
    \]
\end{lemma}

\begin{proof}
    For a relative diagonal \(F \hookrightarrow F \times_{B} F\), its conormal complex is canonically identified with \(\mathbb{L}_{F/B}\).
    In our case, this gives
    \[
        \begin{aligned}
            N_{-} & \simeq \mathbb{L}_{F_{-}/\Gr_{X}(\mathcal{E},d - 1)}, \\
            N_{+} & \simeq \mathbb{L}_{F_{+}/\Gr_{X}(\mathcal{E},d + 1)}.
        \end{aligned}
    \]
    The two projectivization descriptions above have universal fiber sequences
    \[
        \pi_{d - 1}^{*}R_{d} \longrightarrow j_{d - 1}^{*}R_{d - 1} \longrightarrow L_{d}
    \]
    and
    \[
        j_{d}^{*}Q_{d}^{\vee} \longrightarrow \pi_{d}^{*}Q_{d + 1}^{\vee} \longrightarrow L_{d + 1}^{\vee}.
    \]
    The cotangent formula for derived projectivizations \cite[Theorem~4.18]{Jiang2023DerivedGrassmanniansSchurFunctors} gives the asserted identifications of \(N_{-}\) and \(N_{+}\), together with their Tor-amplitudes and ranks.

    It remains to verify the local zero-section descriptions.
    Fix a map \(u \colon T \to F_{-}\) with \(T\) affine.
    The first projection identifies
    \[
        Z^{-} \times_{F_{-}}T \simeq T \times_{\Gr_{X}(\mathcal{E},d - 1)}F_{-}.
    \]
    Under the projectivization description of \(F_{-}\), the diagonal section corresponds to the quotient
    \[
        u^{*}j_{d - 1}^{*}R_{d - 1} \twoheadrightarrow u^{*}L_{d},
    \]
    whose fiber is \(u^{*}\pi_{d - 1}^{*}R_{d}\).
    Applying \eqref{eq:projectivization-affine-chart}, we obtain an open neighborhood
    \[
        \Spec_{T}\Sym\bigl( u^{*}\pi_{d - 1}^{*}R_{d} \otimes (u^{*}L_{d})^{-1} \bigr) \simeq \Spec_{T}\Sym(u^{*}N_{-}),
    \]
    with the diagonal as the zero section.

    Similarly, for a map \(u \colon T \to F_{+}\) with \(T\) affine, the diagonal section corresponds to
    \[
        u^{*}\pi_{d}^{*}Q_{d + 1}^{\vee} \twoheadrightarrow u^{*}L_{d + 1}^{\vee},
    \]
    with fiber \(u^{*}j_{d}^{*}Q_{d}^{\vee}\).
    The same affine chart gives an open neighborhood
    \[
        \Spec_{T}\Sym\bigl( u^{*}j_{d}^{*}Q_{d}^{\vee} \otimes u^{*}L_{d + 1} \bigr) \simeq \Spec_{T}\Sym(u^{*}N_{+}).
    \]
    This proves the local descriptions for both diagonals.
\end{proof}

We can now apply \zcref{lem:zero-section-rees-comparison}.
Since \(N_{+}\) has rank \(d > 0\), the Rees comparison map for \(i_{+}\) at \(n = 0\) is an equivalence.
As this map is the adjunction unit, we obtain
\begin{equation}
    \label{eq:upper-diagonal-structure-sheaf}
    \mathcal{O}_{Z^{+}} \xrightarrow{\sim} b_{+,*}\mathcal{O}_{\widetilde{Z}}.
\end{equation}
For \(i_{-}\), the conormal complex \(N_{-}\) has rank \(r - d\).
Thus the same lemma computes the cofibers of its Rees comparison maps.

We now apply the preceding calculations to compare compositions of flag functors.
For \(d \geq 1\) and \(a,b \in \mathbb{Z}\), consider
\begin{equation}
    \label{eq:mixed-flag-composites}
    \Psi_{d - 1,b}^{L}\Psi_{d - 1,a}, \qquad \Psi_{d,a + 1}\Psi_{d,b + 1}^{L}
\end{equation}
as endofunctors of \(\Perf(\Gr_{X}(\mathcal{E},d))\).
We first compute their Fourier--Mukai kernels and then compare them through the common blow-up \(\widetilde{Z}\).

Recall the Cartesian squares
\[
    \begin{tikzcd}[
            column sep=1.2cm,
            row sep=0.9cm,
            baseline=(minus-source.base)
        ]
        |[alias=minus-source]| Z^{-}
        \arrow[r,"\operatorname{pr}_{2}^{-}"]
        \arrow[d,"\operatorname{pr}_{1}^{-}"']
        & F_{-} \arrow[d,"j_{d - 1}"] \\
        F_{-} \arrow[r,"j_{d - 1}"']
        & \Gr_{X}(\mathcal{E},d - 1)
    \end{tikzcd}
    \hspace{1cm}
    \begin{tikzcd}[
            column sep=1.2cm,
            row sep=0.9cm,
            baseline=(plus-source.base)
        ]
        |[alias=plus-source]| Z^{+}
        \arrow[r,"\operatorname{pr}_{2}^{+}"]
        \arrow[d,"\operatorname{pr}_{1}^{+}"']
        & F_{+} \arrow[d,"\pi_{d}"] \\
        F_{+} \arrow[r,"\pi_{d}"']
        & \Gr_{X}(\mathcal{E},d + 1).
    \end{tikzcd}
\]
For \(i = 1,2\), put
\[
    \begin{aligned}
        q_{-,i} & = \pi_{d - 1} \circ \operatorname{pr}_{i}^{-} \colon Z^{-} \longrightarrow \Gr_{X}(\mathcal{E},d), \\
        q_{+,i} & = j_{d} \circ \operatorname{pr}_{i}^{+} \colon Z^{+} \longrightarrow \Gr_{X}(\mathcal{E},d).
    \end{aligned}
\]
The first and second factors correspond to the target and source of each composition, respectively.
Write \(q_{\pm} = (q_{\pm,1},q_{\pm,2})\) for the induced maps to \(\Gr_{X}(\mathcal{E},d) \times_{X}\Gr_{X}(\mathcal{E},d)\).
We denote the universal flag lines pulled back from the two factors by \(L_{d}^{(i)}\) on \(Z^{-}\) and \(L_{d + 1}^{(i)}\) on \(Z^{+}\).

The two Cartesian squares give the base-change equivalences
\[
    j_{d - 1}^{*}j_{d - 1,*} \simeq (\operatorname{pr}_{1}^{-})_{*}(\operatorname{pr}_{2}^{-})^{*}, \qquad \pi_{d}^{*}\pi_{d,*} \simeq (\operatorname{pr}_{1}^{+})_{*}(\operatorname{pr}_{2}^{+})^{*}.
\]
Applying these equivalences and the projection formula to the expressions for \(\Psi_{k,m}\) and \(\Psi_{k,m}^{L}\) obtained above, we find that
\begin{equation}
    \label{eq:mixed-flag-convolution-kernels}
    \begin{aligned}
        \Psi_{d - 1,b}^{L}\Psi_{d - 1,a}(T) & \simeq q_{-,1,*}\bigl( W^{-}_{a,b} \otimes q_{-,2}^{*}T \bigr), \\
        \Psi_{d,a + 1}\Psi_{d,b + 1}^{L}(T) & \simeq q_{+,1,*}\bigl( W^{+}_{a,b} \otimes q_{+,2}^{*}T \bigr),
    \end{aligned}
\end{equation}
where
\[
    \begin{aligned}
        W^{-}_{a,b} & = q_{-,1}^{*}\det Q_{d}^{\vee} \otimes (L_{d}^{(1)})^{d - b} \otimes (L_{d}^{(2)})^{a}[d - 1], \\
        W^{+}_{a,b} & = \det Q_{d + 1}^{\vee} \otimes (L_{d + 1}^{(1)})^{a + 1} \otimes (L_{d + 1}^{(2)})^{d - b}[d].
    \end{aligned}
\]
The determinant factor in \(W^{-}_{a,b}\) comes from the first flag factor, while \(\det Q_{d + 1}^{\vee}\) in \(W^{+}_{a,b}\) is pulled back from the common quotient of rank \(d + 1\).
Formula~\eqref{eq:mixed-flag-convolution-kernels} shows that the Fourier--Mukai kernels of the two compositions are
\[
    q_{-,*}W^{-}_{a,b}, \qquad q_{+,*}W^{+}_{a,b}.
\]

We next compare these kernels on the common blow-up \(\widetilde{Z}\).
By \zcref{lem:rank-one-flag-blowup}, the ideal sheaf of the exceptional divisor satisfies
\begin{equation}
    \label{eq:common-flag-line-identities}
    \mathcal{O}_{\widetilde{Z}}(-D) \simeq L_{d + 1}^{(1)} \otimes (L_{d}^{(2)})^{-1} \simeq L_{d + 1}^{(2)} \otimes (L_{d}^{(1)})^{-1}.
\end{equation}
Moreover, the universal quotient sequence gives
\[
    \det Q_{d + 1}^{\vee} \simeq b_{-}^{*}q_{-,1}^{*}\det Q_{d}^{\vee} \otimes (L_{d + 1}^{(1)})^{-1}.
\]
Putting \(s = a - b\) and \(n = d + s\), these identities yield
\begin{equation}
    \label{eq:common-flag-kernel-weights}
    \begin{aligned}
        b_{-}^{*}W^{-}_{a,b} \otimes \mathcal{O}_{\widetilde{Z}}(-nD) & \simeq b_{+}^{*}W^{+}_{a,b}[-1], \\
        i_{-}^{*}W^{-}_{a,b}                                         & \simeq \pi_{d - 1}^{*}\det Q_{d}^{\vee} \otimes L_{d}^{n}[d - 1].
    \end{aligned}
\end{equation}
Thus the comparison naturally involves the second composition shifted by \([-1]\).

Recall from \zcref{lem:lower-diagonal-local-model} that
\begin{equation}
    \label{eq:lower-diagonal-conormal}
    N_{-} \simeq \pi_{d - 1}^{*}R_{d} \otimes L_{d}^{-1}, \qquad \rho = r - d, \qquad L_{d} \simeq \mathcal{O}_{\pi_{d - 1}}(-1).
\end{equation}
Let \(\mathcal{R}_{n}\) be the degree-\(n\) component of the extended Rees algebra of \(i_{-}\).
The two canonical maps
\[
    \mathcal{R}_{n} \longrightarrow \mathcal{O}_{Z^{-}}, \qquad \theta_{n} \colon \mathcal{R}_{n} \longrightarrow b_{-,*}\mathcal{O}_{\widetilde{Z}}(-nD)
\]
provide a common source for the comparison.

Define the following Fourier--Mukai functors on \(\QCoh(\Gr_{X}(\mathcal{E},d))\):
\[
    \begin{aligned}
        \mathcal{M}_{a,b} \colon T     & \longmapsto q_{-,1,*}\bigl(W^{-}_{a,b} \otimes \mathcal{R}_{n} \otimes q_{-,2}^{*}T\bigr), \\
        \mathcal{H}^{+}_{a,b} \colon T & \longmapsto q_{-,1,*}\bigl(W^{-}_{a,b} \otimes \cofib(\mathcal{R}_{n} \to \mathcal{O}_{Z^{-}}) \otimes q_{-,2}^{*}T\bigr), \\
        \mathcal{H}^{-}_{a,b} \colon T & \longmapsto q_{-,1,*}\bigl(W^{-}_{a,b} \otimes \cofib(\theta_{n}) \otimes q_{-,2}^{*}T\bigr).
    \end{aligned}
\]

\begin{theorem}
    \label{thm:mixed-flag-rees-comparison}
    For \(d \geq 1\) and \(a,b \in \mathbb{Z}\), put \(s = a - b\) and \(u = -r - s\).
    There are canonical fiber sequences of endofunctors of \(\Perf(\Gr_{X}(\mathcal{E},d))\):
    \begin{equation}
        \label{eq:mixed-flag-comparison-triangles}
        \begin{aligned}
            \mathcal{M}_{a,b} & \longrightarrow \Psi_{d - 1,b}^{L}\Psi_{d - 1,a} \longrightarrow \mathcal{H}^{+}_{a,b}, \\
            \mathcal{M}_{a,b} & \longrightarrow \Psi_{d,a + 1}\Psi_{d,b + 1}^{L}[-1] \longrightarrow \mathcal{H}^{-}_{a,b}.
        \end{aligned}
    \end{equation}
    The functors \(\mathcal{H}^{+}_{a,b}\) and \(\mathcal{H}^{-}_{a,b}\) vanish for \(s < 0\) and \(u < 0\), respectively.
    For \(s \geq 0\) and \(u \geq 0\), respectively, they admit finite filtrations with graded pieces
    \begin{align}
        \label{eq:positive-flag-correction-graded}
        \gr_{t}\mathcal{H}^{+}_{a,b} & \simeq - \otimes \Sym^{t}R_{d} \otimes \Gamma^{s - t}Q_{d},                                                         &  & 0 \leq t \leq s, \\
        \label{eq:negative-flag-correction-graded}
        \gr_{t}\mathcal{H}^{-}_{a,b} & \simeq - \otimes \bigl(\Sym^{t}R_{d} \otimes \Gamma^{u - t}Q_{d}\bigr)^{\vee} \otimes \det \mathcal{E}^{\vee}[2d - r], &  & 0 \leq t \leq u.
    \end{align}
    The graded pieces are attached in the orders \(t = s,s - 1,\ldots,0\) and \(t = 0,1,\ldots,u\), respectively.
\end{theorem}

\begin{proof}
    Put \(n = d + s\).
    The two Rees maps
    \[
        \mathcal{R}_{n} \longrightarrow \mathcal{O}_{Z^{-}}, \qquad \theta_{n} \colon \mathcal{R}_{n} \longrightarrow b_{-,*}\mathcal{O}_{\widetilde{Z}}(-nD)
    \]
    induce morphisms of Fourier--Mukai kernels with common source \(q_{-,*}(W^{-}_{a,b} \otimes \mathcal{R}_{n})\).
    By \eqref{eq:mixed-flag-convolution-kernels}, the first target is the kernel of \(\Psi_{d - 1,b}^{L}\Psi_{d - 1,a}\).
    For the second target, the projection formula, \eqref{eq:common-flag-kernel-weights}, and \eqref{eq:upper-diagonal-structure-sheaf} give
    \[
        \begin{aligned}
            q_{-,*}\bigl(W^{-}_{a,b} \otimes b_{-,*}\mathcal{O}_{\widetilde{Z}}(-nD)\bigr) & \simeq (q_{-}b_{-})_{*}\bigl(b_{-}^{*}W^{-}_{a,b} \otimes \mathcal{O}_{\widetilde{Z}}(-nD)\bigr) \\
                                                                                          & \simeq (q_{+}b_{+})_{*}b_{+}^{*}W^{+}_{a,b}[-1] \simeq q_{+,*}W^{+}_{a,b}[-1].
        \end{aligned}
    \]
    This is the kernel of \(\Psi_{d,a + 1}\Psi_{d,b + 1}^{L}[-1]\).
    Taking cofibers proves \eqref{eq:mixed-flag-comparison-triangles}.

    We now compute \(\mathcal{H}^{+}_{a,b}\) and \(\mathcal{H}^{-}_{a,b}\).
    Let \(\Delta\) denote the relative diagonal of \(\Gr_{X}(\mathcal{E},d)\).
    Since \(q_{-}i_{-} = \Delta\pi_{d - 1}\), the graded pieces of the kernels of \(\mathcal{H}^{+}_{a,b}\) and \(\mathcal{H}^{-}_{a,b}\) are supported on \(\Delta\).
    By \eqref{eq:common-flag-kernel-weights} and \eqref{eq:lower-diagonal-conormal},
    \[
        i_{-}^{*}W^{-}_{a,b} \simeq \omega_{\pi_{d - 1}} \otimes L_{d}^{s}, \qquad N_{-} \simeq \pi_{d - 1}^{*}R_{d} \otimes L_{d}^{-1}, \qquad \rho = r - d.
    \]

    For \(\mathcal{H}^{+}_{a,b}\), the defining cofiber vanishes when \(n \leq 0\).
    If \(n > 0\), the \(t\)-th graded piece of \eqref{eq:rees-augmentation-cofiber-filtration}, for \(0 \leq t < n\), gives the Fourier--Mukai kernel
    \[
        q_{-,*}i_{-,*}\bigl(i_{-}^{*}W^{-}_{a,b} \otimes \Sym^{t}N_{-}\bigr) \simeq \Delta_{*}\bigl(\Sym^{t}R_{d} \otimes \pi_{d - 1,!}L_{d}^{s - t}\bigr).
    \]
    Since \(\pi_{d - 1}\) is a projective bundle of relative dimension \(d - 1\), we have
    \[
        \pi_{d - 1,!}L_{d}^{m} \simeq
        \begin{cases}
            \Gamma^{m}Q_{d}, & m \geq 0, \\
            0,               & 1 - d \leq m < 0.
        \end{cases}
    \]
    If \(s < 0\), every exponent \(s - t\) occurring in the filtration lies in the vanishing range, so \(\mathcal{H}^{+}_{a,b} = 0\).
    If \(s \geq 0\), precisely the indices \(0 \leq t \leq s\) survive, giving \eqref{eq:positive-flag-correction-graded}.
    Their attachment order is \(t = s,s - 1,\ldots,0\), as in \eqref{eq:rees-augmentation-cofiber-filtration}.

    For \(\mathcal{H}^{-}_{a,b}\), \zcref{lem:zero-section-rees-comparison} applies to \(i_{-}\) with conormal rank \(\rho = r - d\).
    Since
    \[
        -\rho - n = -r - s = u,
    \]
    it gives \(\mathcal{H}^{-}_{a,b} = 0\) for \(u < 0\).
    Suppose \(u \geq 0\).
    Using
    \[
        \det N_{-}^{\vee} \simeq \pi_{d - 1}^{*} \bigl(\det \mathcal{E}^{\vee} \otimes \det Q_{d}\bigr) \otimes L_{d}^{\rho},
    \]
    the \(t\)-th graded piece of \eqref{eq:rees-tautological-cofiber-filtration} satisfies
    \[
        i_{-}^{*}W^{-}_{a,b} \otimes (\Sym^{t}N_{-})^{\vee} \otimes \det N_{-}^{\vee}[1 - \rho] \simeq \pi_{d - 1}^{*}\bigl((\Sym^{t}R_{d})^{\vee} \otimes \det \mathcal{E}^{\vee}\bigr) \otimes L_{d}^{t - u}[2d - r].
    \]
    Here we used \(n + t + \rho = t - u\) and \(d - \rho = 2d - r\).
    Pushing forward along \(\pi_{d - 1}\) gives
    \[
        \Delta_{*}\bigl( (\Sym^{t}R_{d})^{\vee} \otimes \det \mathcal{E}^{\vee} \otimes \pi_{d - 1,*}L_{d}^{t - u} \bigr)[2d - r].
    \]
    For \(0 \leq t \leq u\), projective-bundle cohomology gives
    \[
        \pi_{d - 1,*}L_{d}^{t - u} \simeq (\Gamma^{u - t}Q_{d})^{\vee}.
    \]
    These kernels yield \eqref{eq:negative-flag-correction-graded}, with the pieces attached in the order \(t = 0,1,\ldots,u\), by \eqref{eq:rees-comparison-cofiber-recursion}.

    Since the graded pieces of \(\mathcal{H}^{\pm}_{a,b}\) are tensor functors by perfect complexes, both functors preserve perfect complexes.
    As the compositions of flag functors also preserve perfect complexes, the fiber sequences imply the same for \(\mathcal{M}_{a,b}\).
\end{proof}

\begin{remark}
    \label{rem:hsu-ordinary-grassmannians}
    For \(X = \Spec \mathbb{C}\) and \(\mathcal{E} = \mathbb{C}^{N}\), put \(s = a - b\).
    Up to determinant twists and shifts, Hsu \cite[Theorem~5.7 and Proposition~5.10]{Hsu2025ShiftedZeroAffineAlgebra} established the equivalence for \(-N < s < 0\) and the two boundary triangles at \(s = 0,-N\).
    He also computed the first off-wall correction triangles at \(s = 1\) and \(s = -N - 1\), including the two-step filtration of the positive correction kernel with graded pieces \(\Delta_{*}R_{d}\) and \(\Delta_{*}Q_{d}\) \cite[Section~6]{Hsu2025ShiftedZeroAffineAlgebra}.
\end{remark}

We now identify the natural transformations for \(s = 0\) and \(u = 0\).

\begin{lemma}
    \label{lem:flag-left-counit-identification}
    For \(a \in \mathbb{Z}\), the natural transformation
    \begin{equation}
        \label{eq:flag-left-counit}
        \Psi_{d - 1,a}^{L}\Psi_{d - 1,a} \longrightarrow \mathcal{H}^{+}_{a,a} \simeq \id
    \end{equation}
    is the counit of \(\Psi_{d - 1,a}^{L} \dashv \Psi_{d - 1,a}\).
\end{lemma}

\begin{proof}
    At \(s = 0\), we have \(n = d\) and \(i_{-}^{*}W^{-}_{a,a} \simeq \omega_{\pi_{d - 1}}\).
    By \eqref{eq:rees-augmentation-cofiber-filtration}, the \(t\)-th graded piece of \(\mathcal{H}^{+}_{a,a}\) has kernel
    \[
        \Delta_{*}\bigl(\Sym^{t}R_{d} \otimes \pi_{d - 1,!}L_{d}^{-t}\bigr), \qquad 0 \leq t < d.
    \]
    Since \(\pi_{d - 1,!}L_{d}^{-t} = 0\) for \(1 \leq t < d\), only the term \(t = 0\) survives, giving \(\mathcal{H}^{+}_{a,a} \simeq \id\).
    The map to this term is induced by restriction to the diagonal, followed by the counit for \(\pi_{d - 1}\).

    More explicitly, for \(T \in \Perf(\Gr_{X}(\mathcal{E},d))\), the resulting natural transformation is
    \[
        \begin{aligned}
            \Psi_{d - 1,a}^{L}\Psi_{d - 1,a}(T) & \simeq \pi_{d - 1,!}\bigl(L_{d}^{-a} \otimes j_{d - 1}^{*}j_{d - 1,*} (L_{d}^{a} \otimes \pi_{d - 1}^{*}T)\bigr) \\
             & \longrightarrow \pi_{d - 1,!}\pi_{d - 1}^{*}T \longrightarrow T.
        \end{aligned}
    \]
    These are the counits of \(j_{d - 1}^{*} \dashv j_{d - 1,*}\) and \(\pi_{d - 1,!} \dashv \pi_{d - 1}^{*}\), respectively.
    Their composite is the counit of \(\Psi_{d - 1,a}^{L} \dashv \Psi_{d - 1,a}\).
\end{proof}

\begin{lemma}
    \label{lem:flag-right-counit-identification}
    Suppose \(a - b = -r\).
    Relative Serre duality gives a canonical equivalence
    \[
        \Psi_{d,b + 1}^{L} \simeq \bigl(- \otimes \det \mathcal{E}^{\vee}\bigr) [2d + 1 - r]\Psi_{d,a + 1}^{R}.
    \]
    Under this equivalence, the natural transformation in the second fiber sequence of \zcref{thm:mixed-flag-rees-comparison} is the counit of \(\Psi_{d,a + 1} \dashv \Psi_{d,a + 1}^{R}\), tensored by \(\det \mathcal{E}^{\vee}\) and shifted by \([2d - r]\):
    \begin{equation}
        \label{eq:flag-serre-boundary}
        \begin{aligned}
            \Psi_{d,a + 1}\Psi_{d,b + 1}^{L}[-1] & \simeq \bigl(- \otimes \det \mathcal{E}^{\vee}\bigr) [2d - r]\Psi_{d,a + 1}\Psi_{d,a + 1}^{R} \\
             & \longrightarrow \bigl(- \otimes \det \mathcal{E}^{\vee}\bigr) [2d - r] \simeq \mathcal{H}^{-}_{a,b}.
        \end{aligned}
    \end{equation}
\end{lemma}

\begin{proof}
    The first equivalence follows from the adjoint formulas in \zcref{subsec:grassmannians-flags}, using \(a = b - r\) and
    \[
        \det Q_{d + 1}^{\vee} \simeq j_{d}^{*}\det \mathcal{E}^{\vee} \otimes j_{d}^{*}\det R_{d} \otimes L_{d + 1}^{-1}.
    \]

    We identify the natural transformation.
    The counit of \(\Psi_{d,a + 1} \dashv \Psi_{d,a + 1}^{R}\) is the composite
    \[
        \begin{aligned}
            \Psi_{d,a + 1}\Psi_{d,a + 1}^{R}(M) & \simeq j_{d,*}\bigl(L_{d + 1}^{a + 1}
                \otimes \pi_{d}^{*}\pi_{d,*} (L_{d + 1}^{-a - 1} \otimes j_{d}^{!}M)\bigr) \\
             & \longrightarrow j_{d,*}j_{d}^{!}M \longrightarrow M,
        \end{aligned}
    \]
    where the two arrows are the counits of \(\pi_{d}^{*} \dashv \pi_{d,*}\) and \(j_{d,*} \dashv j_{d}^{!}\).

    We compare this composite with the Rees construction.
    Put \(\rho = r - d\), so that \(n = -\rho\).
    The two projections from the exceptional divisor fit into the Cartesian square
    \[
        \begin{tikzcd}
            D \arrow[r,"e_{+}"] \arrow[d,"e_{-}"']
            & F_{+} \arrow[d,"j_{d}"]\\
            F_{-} \arrow[r,"\pi_{d - 1}"']
            & \Gr_{X}(\mathcal{E},d).
        \end{tikzcd}
    \]
    Since \(\theta_{1 - \rho}\) is an equivalence, the diagram in the proof of \zcref{lem:zero-section-rees-comparison} identifies the map to \(\cofib(\theta_{-\rho})\) with the composite
    \begin{equation}
        \label{eq:rees-boundary-map}
        b_{-,*}\mathcal{O}_{\widetilde{Z}}(\rho D) \longrightarrow i_{-,*}e_{-,*}\mathcal{O}_{D}(-\rho) \longrightarrow i_{-,*}\det N_{-}^{\vee}[1 - \rho].
    \end{equation}

    The projection \(e_{+}\) is a projective bundle, so \(e_{+,*}\mathcal{O}_{D} \simeq \mathcal{O}_{F_{+}}\).
    Using \eqref{eq:common-flag-kernel-weights} and \eqref{eq:upper-diagonal-structure-sheaf}, the first arrow of \eqref{eq:rees-boundary-map}, after tensoring by \(W^{-}_{a,b}\) and pushing forward, becomes the map induced by
    \[
        \pi_{d}^{*}\pi_{d,*} \longrightarrow \id.
    \]

    For the second arrow, the relative dualizing complex of \(e_{-}\) is
    \[
        e_{-}^{!}\mathcal{O}_{F_{-}} \simeq e_{-}^{*}\det N_{-} \otimes \mathcal{O}_{D}(-\rho)[\rho - 1].
    \]
    Thus this arrow is obtained from the counit \(e_{-,*}e_{-}^{!} \to \id\).
    Since \(e_{-}\) is the base change of \(j_{d}\), its pushforward under the Cartesian square above identifies this map with the counit
    \[
        j_{d,*}j_{d}^{!} \longrightarrow \id.
    \]
    The determinant twist and shift are given by
    \[
        i_{+}^{*}W^{+}_{a,b}[-1] \simeq j_{d}^{*}\det \mathcal{E}^{\vee} \otimes \omega_{j_{d}}[2d - r].
    \]

    Consequently, the map obtained from the Rees comparison is precisely the composite of the two adjunction counits displayed above, with the stated determinant twist and shift.
\end{proof}

%% file: src/clifford-relations.tex
\section{Clifford relations and Fock decompositions}
\label{sec:clifford-flags}

In this section, we establish the Clifford relations induced by the flag functors on \(K_{0}\).
When \(r > 0\), these relations yield the standard Fock decomposition of the resulting Clifford module.
Motivated by this representation-theoretic structure, we construct a categorical counterpart in the form of a semiorthogonal decomposition of derived Grassmannians using the flag functors.

\subsection{Clifford relations}
Put
\[
    \mathcal{C}_{d} = \Perf(\Gr_{X}(\mathcal{E},d)), \qquad \mathcal{C} = \bigoplus_{d \geq 0}\mathcal{C}_{d},
\]
where objects of \(\mathcal{C}\) have finite support.
For each \(i \in \mathbb{Z}\), define functors \(F_{i}\) and \(F_{i}^{L}\) on \(\mathcal{C}\) by
\begin{equation}
    \label{eq:normalized-flag-functors}
    \begin{aligned}
        F_{i}(M)     & = \Psi_{d,i - r}(M) \otimes (\det Q_{d})^{-1}, &  & M \in \mathcal{C}_{d + 1}, \\
        F_{i}^{L}(N) & = \Psi_{d,i - r}^{L}(N \otimes \det Q_{d}),    &  & N \in \mathcal{C}_{d}.
    \end{aligned}
\end{equation}
We set \(F_{i}|_{\mathcal{C}_{0}} = 0\).
These functors form an adjoint pair \(F_{i}^{L} \dashv F_{i}\).

We first establish the square-zero relations.

\begin{lemma}
    \label{lem:flag-squares}
    For every \(i \in \mathbb{Z}\), there are equivalences of functors
    \begin{equation}
        \label{eq:flag-square-zero}
        F_{i}^{2} \simeq 0, \qquad (F_{i}^{L})^{2} \simeq 0.
    \end{equation}
\end{lemma}

\begin{proof}
    Fix \(d \geq 0\), and consider the composition \(F_{i}^{2} \colon \mathcal{C}_{d + 2} \to \mathcal{C}_{d}\).
    Let
    \[
        p \colon \Flag_{X}(\mathcal{E};d,d + 1,d + 2) \longrightarrow \Flag_{X}(\mathcal{E};d,d + 2)
    \]
    be the morphism forgetting the quotient \(Q_{d + 1}\).
    On the target, put \(K = \ker(Q_{d + 2} \to Q_{d})\), a vector bundle of rank two.
    By the universal property of projectivization, we identify \(p\) with \(\mathbb{P}(K)\), under which \(L_{d + 1} \simeq \mathcal{O}_{p}(1)\).

    By convolution, \(F_{i}^{2}\) is represented by the Fourier--Mukai kernel \(\mathcal{L}\) on \(\Flag_{X}(\mathcal{E};d,d + 1,d + 2)\), where
    \[
        \mathcal{L} = (\det Q_{d})^{-1} \otimes (\det Q_{d + 1})^{-1} \otimes L_{d + 1}^{i - r} \otimes L_{d + 2}^{i - r}.
    \]
    Using
    \[
        \det Q_{d + 1} \simeq \det Q_{d} \otimes L_{d + 1}, \qquad \det K \simeq L_{d + 1} \otimes L_{d + 2},
    \]
    we obtain
    \[
        \begin{aligned}
            \mathcal{L} & \simeq (\det Q_{d})^{-2} \otimes L_{d + 1}^{i - r - 1} \otimes L_{d + 2}^{i - r} \\
                        & \simeq p^{*}\bigl( (\det Q_{d})^{-2} \otimes (\det K)^{i - r} \bigr) \otimes \mathcal{O}_{p}(-1).
        \end{aligned}
    \]
    Since \(p_{*}\mathcal{O}_{p}(-1) \simeq 0\), the projection formula gives \(p_{*}\mathcal{L} \simeq 0\).
    Hence \(F_{i}^{2} \simeq 0\).
    Taking left adjoints gives \((F_{i}^{L})^{2} \simeq 0\).
\end{proof}

The definitions of \(F_{i}\) and \(F_{i}^{L}\) and the projection formula give, on \(\mathcal{C}_{d}\) for \(d \geq 1\),
\begin{equation}
    \label{eq:normalized-mixed-flag-composites}
    \begin{aligned}
        F_{i}^{L}F_{j} & \simeq \Psi_{d - 1,i - r}^{L}\Psi_{d - 1,j - r}, \\
        F_{j}F_{i}^{L} & \simeq \Psi_{d,j - r + 1}\Psi_{d,i - r + 1}^{L}.
    \end{aligned}
\end{equation}

We now consider the induced operators \(f_{i} = [F_{i}]\) and \(f_{i}^{*} = [F_{i}^{L}]\) on \(K_{0}(\mathcal{C})\).
Let \(p \colon \mathbb{P}_{X}(\mathcal{E}) \to X\) be the projection, and put
\[
    c_{s} = [p_{*}\mathcal{O}_{\mathbb{P}_{X}(\mathcal{E})}(s)] \in K_{0}(\Perf(X)), \qquad s \in \mathbb{Z}.
\]
These coefficients act on \(K_{0}(\mathcal{C})\) by tensoring with pullbacks from \(X\).

\begin{theorem}
    \label{thm:k-theory-clifford-relations}
    The operators \(f_{i}\) and \(f_{i}^{*}\) on \(K_{0}(\mathcal{C})\) are square-zero and satisfy, for all \(i,j \in \mathbb{Z}\),
    \begin{equation}
        \label{eq:k-theory-flag-relations}
        \begin{aligned}
            f_{i}f_{j} + f_{j}f_{i} & = f_{i}^{*}f_{j}^{*} + f_{j}^{*}f_{i}^{*} = 0, \\
            f_{i}^{*}f_{j} + f_{j}f_{i}^{*} & = c_{j - i}\,\id.
        \end{aligned}
    \end{equation}
\end{theorem}

\begin{proof}
    We first consider the relations between the \(f_{i}\).
    For \(i = j\), the same-type relations follow from \zcref{lem:flag-squares}.
    By symmetry, it remains to consider \(i < j\).
    Put \(m = j - i - 1\).
    Let
    \[
        q \colon \Flag_{X}(\mathcal{E};d,d + 1,d + 2) \longrightarrow \Flag_{X}(\mathcal{E};d,d + 2)
    \]
    be the morphism forgetting the intermediate quotient, and put \(K = \ker(Q_{d + 2} \to Q_{d})\).
    As in the proof of \zcref{lem:flag-squares}, \(q\) is the projection of \(\mathbb{P}(K)\), and
    \[
        L_{d + 1} \simeq \mathcal{O}_{q}(1), \qquad L_{d + 2} \simeq q^{*}\det K \otimes \mathcal{O}_{q}(-1).
    \]

    The compositions \(F_{i}F_{j}\) and \(F_{j}F_{i}\) are represented by the kernels \(W_{ij}\) and \(W_{ji}\) on \(\Flag_{X}(\mathcal{E};d,d + 1,d + 2)\), respectively, where
    \[
        \begin{aligned}
            W_{ij} & = (\det Q_{d})^{-2} \otimes L_{d + 1}^{i - r - 1} \otimes L_{d + 2}^{j - r}, \\
            W_{ji} & = (\det Q_{d})^{-2} \otimes L_{d + 1}^{j - r - 1} \otimes L_{d + 2}^{i - r}.
        \end{aligned}
    \]
    Projective-bundle cohomology gives
    \[
        \begin{aligned}
            q_{*}W_{ij} & \simeq (\det Q_{d})^{-2} \otimes (\det K)^{i - r} \otimes \Gamma^{m}K[-1], \\
            q_{*}W_{ji} & \simeq (\det Q_{d})^{-2} \otimes (\det K)^{i - r} \otimes \Sym^{m}K.
        \end{aligned}
    \]
    For any vector bundle \(K\), we have \([\Gamma^{m}K] = [\Sym^{m}K]\) in \(K_{0}\).
    Indeed, after pullback to the complete flag bundle of \(K\), both admit filtrations with the same line-bundle factors, and the pullback is injective on \(K_{0}\).
    The two pushforwards therefore have opposite classes.
    Further pushforward gives
    \[
        f_{i}f_{j} + f_{j}f_{i} = 0.
    \]
    The kernels of the left-adjoint compositions are obtained from these kernels by Grothendieck duality.
    Their classes cancel as well, giving \(f_{i}^{*}f_{j}^{*} + f_{j}^{*}f_{i}^{*} = 0\).

    For the mixed relation, put \(s = j - i\) and first suppose \(d \geq 1\).
    Apply \zcref{thm:mixed-flag-rees-comparison} with
    \[
        a = j - r, \qquad b = i - r,
    \]
    using \eqref{eq:normalized-mixed-flag-composites}.
    Taking classes in the two fiber sequences eliminates the common term \(\mathcal{M}_{a,b}\) and gives
    \[
        f_{i}^{*}f_{j} + f_{j}f_{i}^{*} = [\mathcal{H}^{+}_{a,b}] - [\mathcal{H}^{-}_{a,b}].
    \]

    The symmetric-power filtration associated with
    \[
        R_{d} \longrightarrow \mathcal{E} \longrightarrow Q_{d}
    \]
    gives
    \[
        \sum_{t = 0}^{k} [\Sym^{t}R_{d}][\Sym^{k - t}Q_{d}] = [\Sym^{k}\mathcal{E}], \qquad k \geq 0.
    \]
    Using the equality of symmetric- and divided-power classes, the filtrations in \zcref{thm:mixed-flag-rees-comparison} therefore yield
    \[
        \begin{aligned}
            [\mathcal{H}^{+}_{a,b}] & = [\Sym^{s}\mathcal{E}]\,\id, \\
            [\mathcal{H}^{-}_{a,b}] & = (-1)^{r}[\det \mathcal{E}^{\vee}] [(\Sym^{-r - s}\mathcal{E})^{\vee}]\,\id.
        \end{aligned}
    \]

    By the generalized Serre theorem \cite[Theorem~5.6]{Jiang2025DerivedProjectivizations},
    \begin{equation}
        \label{eq:projective-euler-coefficient}
        c_{s} = [\Sym^{s}\mathcal{E}] - (-1)^{r}[\det \mathcal{E}^{\vee}] [(\Sym^{-r - s}\mathcal{E})^{\vee}].
    \end{equation}
    This proves the mixed relation on \(\mathcal{C}_{d}\) for \(d \geq 1\).

    Finally, on \(\mathcal{C}_{0} = \Perf(X)\), we have \(F_{i}^{L}F_{j} = 0\).
    Since the flag correspondence at quotient rank zero is
    \[
        X\xleftarrow{p}\mathbb{P}_{X}(\mathcal{E}) \xrightarrow{\id}\mathbb{P}_{X}(\mathcal{E}),
    \]
    the definitions give
    \[
        F_{j}F_{i}^{L} \simeq - \otimes p_{*}\mathcal{O}_{\mathbb{P}_{X}(\mathcal{E})}(j - i).
    \]
    Taking classes gives the required relation on \(\mathcal{C}_{0}\).
\end{proof}

The relations above define an action of the Clifford algebra associated with the \(K_{0}(\Perf(X))\)-valued pairing \((c_{j - i})_{i,j \in \mathbb{Z}}\).
When \(r > 0\), its restriction to \(0 \leq i,j < r\) is nondegenerate, since
\[
    c_{j - i} =
    \begin{cases}
        [\Sym^{j - i}\mathcal{E}], & j \geq i, \\
        0,                       & j < i,
    \end{cases}
    \qquad 0 \leq i,j < r.
\]
Thus the pairing matrix is upper unitriangular, and an explicit change of annihilation operators gives the standard Clifford relations.

\begin{corollary}
    \label{cor:finite-clifford-action}
    Suppose \(r > 0\).
    For \(0 \leq i < r\), define
    \begin{equation}
        \label{eq:clifford-annihilation-generators}
        \bar{f}_{i}^{*} = \sum_{a = 0}^{r - 1 - i} (-1)^{a}[\bigwedge^{a}\mathcal{E}]\,f_{i + a}^{*}.
    \end{equation}
    Then the operators \(f_{i}\) and \(\bar{f}_{i}^{*}\) are square-zero and satisfy, for all \(0 \leq i,j < r\),
    \[
        \begin{aligned}
            f_{i} f_{j} + f_{j} f_{i} & = \bar{f}_{i}^{*}\bar{f}_{j}^{*} + \bar{f}_{j}^{*}\bar{f}_{i}^{*} = 0, \\
            \bar{f}_{i}^{*}f_{j} + f_{j}\bar{f}_{i}^{*} & = \delta_{ij}\,\id.
        \end{aligned}
    \]
\end{corollary}

\begin{proof}
    The identity
    \[
        \left(\sum_{a \geq 0}(-1)^{a}[\bigwedge^{a}\mathcal{E}]z^{a}\right) \left(\sum_{m \geq 0}[\Sym^{m}\mathcal{E}]z^{m}\right) = 1
    \]
    gives
    \[
        \begin{aligned}
            \bar{f}_{i}^{*}f_{j} + f_{j}\bar{f}_{i}^{*} & = \sum_{a = 0}^{r - 1 - i} (-1)^{a}[\bigwedge^{a}\mathcal{E}]c_{j - i - a} \\
             & = \delta_{ij}\,\id.
        \end{aligned}
    \]
    The same-type anticommutation relations follow from \zcref{thm:k-theory-clifford-relations} by \(K_{0}(\Perf(X))\)-linearity.
    The square-zero property of the \(f_{i}\) follows from \zcref{lem:flag-squares}, while that of the \(\bar{f}_{i}^{*}\) follows because they are \(K_{0}(\Perf(X))\)-linear combinations of pairwise anticommuting square-zero operators.
\end{proof}

In contrast, suppose \(r = -m \leq 0\).
At any geometric point \(\bar{x} \to X\), the generating series of derived symmetric powers satisfies
\[
    \sum_{a \geq 0} [\Sym^{a}\mathcal{E}_{\bar{x}}]z^{a} = (1 - z)^{m} \quad \text{in }K_{0}(k(\bar{x}))[[z]].
\]
Since \([\det\mathcal{E}_{\bar{x}}^{\vee}] = 1\), the generalized Serre formula \eqref{eq:projective-euler-coefficient} gives, for \(0 \leq s \leq m\),
\[
    \begin{aligned}
        c_{s}|_{\bar{x}} & = (-1)^{s}\binom{m}{s} - (-1)^{m}(-1)^{m - s}\binom{m}{m - s} \\
         & = 0.
    \end{aligned}
\]
Both terms vanish outside this range, so \(c_{s}|_{\bar{x}} = 0\) for every \(s \in \mathbb{Z}\).
Thus the Clifford pairing becomes zero at every geometric point, and no finite restriction can be nondegenerate as in the positive-rank case.
\subsection{Categorical Fock decompositions}
\label{subsec:categorical-fock-decompositions}

Suppose that \(r > 0\), and put \(M = K_{0}(\mathcal{C})\).
By \zcref{cor:finite-clifford-action}, \(M\) is a module over the standard Clifford algebra with \(r\) pairs of generators.
For \(-1 \leq i < r\), put
\[
    M_{i} = \bigcap_{j = i + 1}^{r - 1}\ker\bar{f}_{j}^{*} = \bigcap_{j = i + 1}^{r - 1}\ker f_{j}^{*}.
\]
This defines an increasing filtration of \(M\):
\[
    V = M_{-1} \subseteq M_{0} \subseteq \cdots \subseteq M_{r - 1} = M,
\]
where \(V = M_{-1}\) is the vacuum module.
For \(0 \leq i < r\), the Clifford relations give
\[
    \begin{aligned}
        f_{i}M_{i} & = f_{i}M_{i - 1} = \ker(f_{i}|_{M_{i}}), \\
        f_{i}^{*}M_{i} & = M_{i - 1} = \im(f_{i}^{*}\cdots f_{r - 1}^{*}).
    \end{aligned}
\]
Consequently, we have a direct sum decomposition
\begin{equation}
    \label{eq:fock-module-recursion}
    M_{i} = f_{i}M_{i - 1} \oplus M_{i - 1}.
\end{equation}
Iterating these decompositions yields the Fock decomposition
\begin{equation}
    \label{eq:fock-module-decomposition}
    M \simeq \bigoplus_{I \subseteq [r]}f_{I}V,
\end{equation}
where \([r] = \{0,\ldots,r - 1\}\) and \(f_{I} = f_{i_{m}}\cdots f_{i_{1}}\) for \(I = \{i_{1} < \cdots < i_{m}\}\).

Jiang, Li, and Zhao \cite{JiangEtAl2025InstantonsAffineVertexAlgebra} observed that the Fock decomposition \eqref{eq:fock-module-decomposition} is the \(K_{0}\)-shadow of the semiorthogonal decompositions of derived Grassmannians, and used these decompositions to construct the Clifford action.
Such decompositions were established by Toda \cite{Toda2023QuotSchemesHallProducts} for Quot schemes of locally free quotients over smooth quasi-projective complex varieties using categorified Hall products, and in \cite{Jiang2023DerivedCategoriesGrassmannians} for derived Grassmannians over arbitrary prestacks over \(\mathbb{Q}\) using incidence correspondences.
Reversing the direction of \cite{JiangEtAl2025InstantonsAffineVertexAlgebra}, we deduce these decompositions from the categorical Clifford relations, which gives a new proof over arbitrary bases.

\begin{definition}
    For \(-1 \leq i < r\), define the full subcategory
    \[
        \mathcal{C}^{(i)} = \bigcap_{j = i + 1}^{r - 1}\ker F_{j}^{L}.
    \]
    These subcategories form an increasing filtration
    \[
        \mathcal{V} = \mathcal{C}^{(-1)} \subseteq \mathcal{C}^{(0)} \subseteq \cdots \subseteq \mathcal{C}^{(r - 1)} = \mathcal{C},
    \]
    where \(\mathcal{V}\) is called the vacuum category.
    We write \(\mathcal{V}_{d} = \mathcal{V} \cap \mathcal{C}_{d}\).
\end{definition}

We now study the action of the flag functors and their left adjoints on this filtration.
To describe the subcategories \(\mathcal{C}^{(i)}\) as images of compositions of left adjoints, we first prove the following vanishing relation.

\begin{lemma}
    \label{lem:flag-word-vanishing}
    For every \(0 \leq i \leq j < r\), we have
    \[
        F_{r - 1}\cdots F_{i}F_{j} \simeq 0.
    \]
\end{lemma}

\begin{proof}
    Fix \(d \geq 0\) and put \(m = r - i\).
    By \eqref{eq:normalized-flag-functors} and convolution, \(F_{r - 1}\cdots F_{i} \colon \mathcal{C}_{d + m} \to \mathcal{C}_{d}\) is represented by the following kernel on \(\Flag_{X}(\mathcal{E};d,d + 1,\ldots,d + m)\):
    \[
        \bigotimes_{v = 1}^{m} \left((\det Q_{d + v - 1})^{-1} \otimes L_{d + v}^{-v}\right) \simeq (\det Q_{d + m})^{-m},
    \]
    where we use \(\det Q_{d + v} \simeq \det Q_{d + v - 1} \otimes L_{d + v}\).

    Consider the forgetful morphism
    \[
        p \colon \Flag_{X}(\mathcal{E};d,d + 1,\ldots,d + m) \longrightarrow \Flag_{X}(\mathcal{E};d,d + m).
    \]
    This is a complete flag bundle, so \(p_{*}\mathcal{O} \simeq \mathcal{O}\).
    Thus the same functor is represented by the kernel \((\det Q_{d + m})^{-m}\) on \(\Flag_{X}(\mathcal{E};d,d + m)\).

    Composing with \(F_{j}\), we obtain the following kernel on \(\Flag_{X}(\mathcal{E};d,d + m,d + m + 1)\):
    \[
        (\det Q_{d + m})^{-m - 1} \otimes L_{d + m + 1}^{j - r} \simeq (\det Q_{d + m + 1})^{-m - 1} \otimes L_{d + m + 1}^{j - i + 1}.
    \]
    Consider the forgetful morphism
    \[
        q \colon \Flag_{X}(\mathcal{E};d,d + m,d + m + 1) \longrightarrow \Flag_{X}(\mathcal{E};d,d + m + 1).
    \]
    Put \(K = \ker(Q_{d + m + 1} \to Q_{d})\), a vector bundle of rank \(m + 1\) on the target.
    By the universal property of projectivization, we identify \(q\) with \(\mathbb{P}(K^{\vee})\), under which \(L_{d + m + 1} \simeq \mathcal{O}_{q}(-1)\).
    Since \(1 \leq j - i + 1 \leq m\), the projective-bundle formula gives
    \[
        q_{*}L_{d + m + 1}^{j - i + 1} \simeq 0.
    \]
    The desired vanishing follows from the projection formula.
\end{proof}

We can now identify the images of the flag functors and their left adjoints on the filtration.

\begin{lemma}
    \label{lem:fock-filtration-compatibility}
    For every \(0 \leq i < r\), we have
    \[
        \begin{aligned}
            F_{i}(\mathcal{C}^{(i)}) & = F_{i}(\mathcal{C}^{(i - 1)}) = \mathcal{C}^{(i)} \cap \ker F_{i}, \\
            F_{i}^{L}(\mathcal{C}^{(i)}) & = \mathcal{C}^{(i - 1)} = \im(F_{i}^{L}\cdots F_{r - 1}^{L}).
        \end{aligned}
    \]
\end{lemma}

\begin{proof}
    The mixed equivalence and fiber sequence below follow from \zcref{thm:mixed-flag-rees-comparison,lem:flag-left-counit-identification}, together with the generalized Serre theorem on \(\mathcal{C}_{0}\).
    The vanishing relation follows by taking left adjoints in \zcref{lem:flag-word-vanishing}:
    \begin{gather}
        F_{j}^{L}F_{i} \simeq F_{i}F_{j}^{L}[-1], \qquad 0 \leq i < j < r,
        \label{eq:categorical-clifford-mixed}\\
        F_{j}^{L}F_{i}^{L}\cdots F_{r - 1}^{L} \simeq 0, \qquad 0 \leq i \leq j < r,
        \label{eq:categorical-clifford-adjoint-vanishing}\\
        F_{i}^{L}F_{i}\xrightarrow{\epsilon_{i}}\id \longrightarrow F_{i}F_{i}^{L}, \qquad 0 \leq i < r.
        \label{eq:categorical-clifford-projection}
    \end{gather}
    Here the last line is a fiber sequence, and \(\epsilon_{i}\) is the adjunction counit.

    We first prove the second assertion by descending induction on \(i\).
    For \(i = r - 1\), we have \(\mathcal{C}^{(r - 1)} = \mathcal{C}\), which is the essential image of the empty composition.
    Suppose that
    \[
        \mathcal{C}^{(i)} = \im(F_{i + 1}^{L}\cdots F_{r - 1}^{L}).
    \]
    We show that
    \[
        F_{i}^{L}(\mathcal{C}^{(i)}) = \mathcal{C}^{(i - 1)} = \im(F_{i}^{L}\cdots F_{r - 1}^{L}).
    \]
    By the induction hypothesis and \eqref{eq:categorical-clifford-adjoint-vanishing}, we have
    \[
        F_{i}^{L}(\mathcal{C}^{(i)}) = \im(F_{i}^{L}\cdots F_{r - 1}^{L}) \subseteq \bigcap_{j = i}^{r - 1}\ker F_{j}^{L} = \mathcal{C}^{(i - 1)}.
    \]
    Conversely, let \(N \in \mathcal{C}^{(i - 1)}\).
    By \eqref{eq:categorical-clifford-mixed}, we have \(F_{i}N \in \mathcal{C}^{(i)}\).
    Since \(F_{i}^{L}N = 0\), the fiber sequence \eqref{eq:categorical-clifford-projection} gives
    \[
        N \simeq F_{i}^{L}F_{i}N \in F_{i}^{L}(\mathcal{C}^{(i)}).
    \]
    This proves the reverse inclusion and completes the induction.

    For the first assertion, the mixed equivalence and \(F_{i}^{2} \simeq 0\) give
    \[
        F_{i}(\mathcal{C}^{(i - 1)}) \subseteq F_{i}(\mathcal{C}^{(i)}) \subseteq \mathcal{C}^{(i)} \cap \ker F_{i}.
    \]
    Conversely, for \(N \in \mathcal{C}^{(i)} \cap \ker F_{i}\), the fiber sequence \eqref{eq:categorical-clifford-projection} and the second assertion give
    \[
        N \simeq F_{i}F_{i}^{L}N \in F_{i}(\mathcal{C}^{(i - 1)}).
    \]
    Hence both inclusions are equalities.
\end{proof}

\begin{theorem}
    \label{thm:categorical-fock-decomposition}
    Suppose that \(r > 0\).
    For every \(0 \leq i < r\), the restriction
    \[
        F_{i} \colon \mathcal{C}^{(i - 1)} \longrightarrow \mathcal{C}^{(i)}
    \]
    is fully faithful and gives a semiorthogonal decomposition
    \begin{equation}
        \label{eq:categorical-fock-recursion}
        \mathcal{C}^{(i)} = \left\langle F_{i}\mathcal{C}^{(i - 1)}, \mathcal{C}^{(i - 1)} \right\rangle.
    \end{equation}
    For \(I = \{i_{1} < \cdots < i_{m}\} \subseteq [r]\), put \(F_{I} = F_{i_{m}}\cdots F_{i_{1}}\) and \(F_{\varnothing} = \id\).
    Consequently, for every \(d \geq 0\), the functors
    \[
        F_{I} \colon \mathcal{V}_{d + |I|} \longrightarrow \mathcal{C}_{d}, \qquad I \subseteq [r],
    \]
    are fully faithful and give a semiorthogonal decomposition
    \begin{equation}
        \label{eq:categorical-fock-decomposition}
        \mathcal{C}_{d} = \left\langle F_{I}\mathcal{V}_{d + |I|} \ \middle|\ I \subseteq [r] \right\rangle.
    \end{equation}
    The components are ordered in reverse colexicographic order: at the largest index where two subsets differ, the subset containing that index comes first.
\end{theorem}
\begin{proof}
    Fix \(0 \leq i < r\).
    By the fiber sequence \eqref{eq:categorical-clifford-projection} and \((F_{i}^{L})^{2} \simeq 0\), we obtain a semiorthogonal decomposition
    \[
        \mathcal{C} = \left\langle F_{i}(\ker F_{i}^{L}), \ker F_{i}^{L} \right\rangle.
    \]
    Indeed, for \(N \in \ker F_{i}^{L}\), the counit \(F_{i}^{L}F_{i}N \to N\) is an equivalence, so \(F_{i}\) is fully faithful on \(\ker F_{i}^{L}\).
    Semiorthogonality follows from adjunction, while the fiber sequence gives generation, since \(F_{i}^{L}F_{i}M \in \ker F_{i}^{L}\) and \(F_{i}^{L}M \in \ker F_{i}^{L}\) for every \(M \in \mathcal{C}\).

    We now restrict this decomposition to \(\mathcal{C}^{(i)}\).
    By \zcref{lem:fock-filtration-compatibility}, for every \(M \in \mathcal{C}^{(i)}\), we have
    \[
        F_{i}^{L}F_{i}M \in \mathcal{C}^{(i - 1)}, \qquad F_{i}F_{i}^{L}M \in F_{i}\mathcal{C}^{(i - 1)}.
    \]
    Moreover, \(\mathcal{C}^{(i)} \cap \ker F_{i}^{L} = \mathcal{C}^{(i - 1)}\).
    Thus the same fiber sequence gives \eqref{eq:categorical-fock-recursion}, with \(F_{i}\) fully faithful on \(\mathcal{C}^{(i - 1)}\).

    Iterating these decompositions from \(i = r - 1\) down to \(i = 0\) gives the fully faithful components \(F_{I}\mathcal{V}\).
    At each step, the components indexed by subsets containing \(i\) precede those indexed by subsets not containing \(i\), yielding the reverse colexicographic order.
    Taking the degree-\(d\) components gives \eqref{eq:categorical-fock-decomposition}.
\end{proof}

By \zcref{lem:fock-filtration-compatibility}, we have
\[
    \mathcal{V} = \im F_{[r]}^{L}.
\]
In particular, \(\mathcal{V}_{d} = 0\) for \(d < r\).
Passing to \(K_{0}\), \zcref{thm:categorical-fock-decomposition} identifies \(K_{0}(\mathcal{V})\) with the vacuum module \(V\) and recovers the Fock decomposition \eqref{eq:fock-module-decomposition}.

In \zcref{sec:integral-grassmannian-decompositions}, we establish categorical analogues of the extended Clifford relations of Jiang, Li, and Zhao \cite{JiangEtAl2025InstantonsAffineVertexAlgebra}, recovering their spin representation on \(K_{0}\).
In particular, we identify the vacuum category geometrically as
\[
    \mathcal{V}_{e + r} \simeq \Perf\bigl(\Gr_{X}(\mathcal{E}^{\vee}[1],e)\bigr), \qquad e \geq 0,
\]
and show that our semiorthogonal decompositions agree with those in \cite{Jiang2023DerivedCategoriesGrassmannians}, including the embedding functors and the ordering of the components.

%% file: src/spherical-flags.tex
\section{Rank zero: spherical twists and flops}
\label{sec:spherical-flags}

Throughout this section, assume that \(r = 0\).
Applying \zcref{thm:mixed-flag-rees-comparison} with \(d\) replaced by \(d + 1\) and \(a = b\) gives two fiber sequences with a common first term.
By \zcref{lem:flag-left-counit-identification,lem:flag-right-counit-identification}, their maps to the last terms are the counits of \(\Psi_{d,a}^{L} \dashv \Psi_{d,a}\) and \(\Psi_{d + 1,a + 1} \dashv \Psi_{d + 1,a + 1}^{R}\), respectively, after the determinant twist and shift in the second case.
We will use this relation to prove sphericality and derive the flop--flop formulas from the incidence equivalences.

\subsection{Spherical flag functors}
\label{subsec:spherical-functors}

Spherical twists were introduced by Seidel and Thomas in the setting of spherical objects \cite{SeidelThomas2001BraidGroupActions}.
We recall the corresponding notion for functors.

\begin{definition}
    \label{def:spherical-functor}
    Let \(F \colon \mathcal{C} \to \mathcal{D}\) be a functor of stable \(\infty\)-categories with left and right adjoints \(L = F^{L}\) and \(R = F^{R}\).
    We define its cotwist and twist by
    \[
        C_{F} = \fib(LF \to \id_{\mathcal{C}}), \qquad T_{F} = \cofib(FR \to \id_{\mathcal{D}}).
    \]
    We call \(F\) spherical if \(C_{F}\) and \(T_{F}\) are autoequivalences.
\end{definition}

For convenience, our definitions of twists and cotwists may differ from the usual conventions by shifts.
We use the following sphericality criterion of Anno--Logvinenko \cite[Theorem~5.1 and Lemmas~5.6--5.7]{AnnoLogvinenko2017SphericalDGFunctors}, which also holds for stable \(\infty\)-categories by the results of Christ \cite[Lemma~2.13, Remark~2.14, and Theorem~2.15]{Christ2023SphericalMonadicAdjunctions}.

\begin{theorem}
    \label{thm:spherical-criterion}
    Suppose that \(T_{F}\) is an autoequivalence and that the canonical map induced by adjunction
    \[
        \alpha_{F} \colon C_{F}R \longrightarrow LFR \longrightarrow L
    \]
    is an equivalence.
    Then \(F\) is spherical.
\end{theorem}

The criterion requires the canonical map \(\alpha_{F}\) to be an equivalence, whereas our geometric calculations only give a natural equivalence \(C_{F}R \simeq L\).
Inspired by Meachan's argument for Fourier--Mukai functors \cite[Proposition~3.2]{Meachan2021SphericalFunctors}, we use the following formal lemma to verify this condition.

\begin{lemma}
    \label{lem:spherical-rigidity}
    Suppose that \(T_{F}\) is an autoequivalence and there exists a natural equivalence \(C_{F}R \simeq L\).
    Then the canonical map \(\alpha_{F}\) is an equivalence.
\end{lemma}

\begin{proof}
    The left adjoint of the twist admits the dual description (cf. \cite[Equation~(5.2) and Proposition~5.3]{AnnoLogvinenko2017SphericalDGFunctors})
    \[
        T_{F}^{L} \simeq \fib(\id_{\mathcal{D}} \to FL),
    \]
    where the map is the unit of \(L \dashv F\).
    Since \(T_{F}\) is an autoequivalence, so is \(T_{F}^{L}\).

    Let \(\alpha_{F} \colon C_{F}R \to L\) be the canonical map, and put \(K = \cofib(\alpha_{F})\).
    We first show that \(FK \simeq 0\) by identifying it with the cofiber of the counit
    \[
        T_{F}^{L}T_{F} \longrightarrow \id_{\mathcal{D}}
    \]
    of the adjunction \(T_{F}^{L} \dashv T_{F}\).

    Let \(\eta \colon \id_{\mathcal{D}} \to FL\) and \(\varepsilon \colon FR \to \id_{\mathcal{D}}\) be the unit and counit of \(L \dashv F\) and \(F \dashv R\), respectively.
    Naturality gives the following diagram, whose rows and columns are cofiber sequences:
    \[
        \begin{tikzcd}[column sep=huge, row sep=large]
            FR \arrow[r,"\eta FR"] \arrow[d,"\varepsilon"'] & FLFR \arrow[r] \arrow[d,"FL\varepsilon"] & FC_{F}R \arrow[d] \\
            \id_{\mathcal{D}} \arrow[r,"\eta"'] \arrow[d] & FL \arrow[r] \arrow[d] & T_{F}^{L}[1] \arrow[d] \\
            T_{F} \arrow[r,"\eta T_{F}"'] & FLT_{F} \arrow[r] & T_{F}^{L}T_{F}[1].
        \end{tikzcd}
    \]
    Here the first row follows from the triangle identity for \(L \dashv F\).
    The map \(FC_{F}R \to T_{F}^{L}[1]\) is the composite
    \[
        FC_{F}R \xrightarrow{F\alpha_{F}} FL \longrightarrow T_{F}^{L}[1].
    \]

    Filling in the corresponding pushouts gives
    \[
        \begin{tikzcd}[column sep=huge, row sep=large]
            FC_{F}R \arrow[r,"F\alpha_{F}"] \arrow[d] & FL \arrow[r] \arrow[d] & T_{F}^{L}[1] \arrow[d] \\
            0 \arrow[r] & FK \arrow[r] \arrow[d] & T_{F}^{L}T_{F}[1] \arrow[d,"{\varepsilon_{T_{F}^{L},T_{F}}[1]}"] \\
            & 0 \arrow[r] & \id_{\mathcal{D}}[1].
        \end{tikzcd}
    \]
    All three squares are pushouts.
    The rightmost map is the shifted counit of \(T_{F}^{L} \dashv T_{F}\).
    Hence
    \[
        FK \simeq \cofib\bigl(T_{F}^{L}T_{F} \xrightarrow{\varepsilon_{T_{F}^{L},T_{F}}} \id_{\mathcal{D}}\bigr).
    \]
    Since \(T_{F}^{L}\) is an autoequivalence, this counit is an equivalence, and therefore \(FK = 0\).

    Now let \(M \in \mathcal{D}\) and \(Z \in \ker F\).
    By adjunction,
    \[
        \Map_{\mathcal{C}}(LM,Z) \simeq \Map_{\mathcal{D}}(M,FZ) \simeq *.
    \]
    The assumed natural equivalence \(C_{F}R \simeq L\) also gives \(\Map_{\mathcal{C}}(C_{F}RM,Z) \simeq *\).
    Applying \(\Map_{\mathcal{C}}(-,Z)\) to the cofiber sequence
    \[
        C_{F}RM \xrightarrow{\alpha_{F}} LM \longrightarrow KM
    \]
    shows that \(\Map_{\mathcal{C}}(KM,Z) \simeq *\).
    Since \(FKM = 0\), we may take \(Z = KM\), which implies \(KM = 0\).
    Thus \(\alpha_{F}\) is an equivalence.
\end{proof}

For the flag functors, use \zcref{def:spherical-functor} to write
\[
    C_{d,a} = C_{\Psi_{d,a}}, \qquad T_{d,a} = T_{\Psi_{d,a}}.
\]

The counit identifications in \zcref{lem:flag-left-counit-identification,lem:flag-right-counit-identification} give the following relation.

\begin{lemma}
    \label{lem:flag-cotwist}
    For every \(d \geq 0\) and \(a \in \mathbb{Z}\), there is a canonical equivalence
    \begin{equation}
        \label{eq:flag-cotwist-twist-comparison}
        C_{d,a} \simeq \bigl(- \otimes \det \mathcal{E}^{\vee}\bigr) T_{d + 1,a + 1}[2d + 1].
    \end{equation}
\end{lemma}

\begin{proof}
    Apply \zcref{thm:mixed-flag-rees-comparison} with \(d\) replaced by \(d + 1\) and both parameters equal to \(a\).
    By \zcref{lem:flag-left-counit-identification}, the first fiber sequence identifies its common first term with
    \[
        \mathcal{M}_{a,a} \simeq \fib(\Psi_{d,a}^{L}\Psi_{d,a} \to \id) = C_{d,a}.
    \]
    By \zcref{lem:flag-right-counit-identification}, the second sequence becomes
    \[
        C_{d,a} \longrightarrow \bigl(- \otimes \det \mathcal{E}^{\vee}\bigr) \Psi_{d + 1,a + 1}\Psi_{d + 1,a + 1}^{R}[2d + 2] \longrightarrow \bigl(- \otimes \det \mathcal{E}^{\vee}\bigr)[2d + 2],
    \]
    where the last arrow is the counit of \(\Psi_{d + 1,a + 1} \dashv \Psi_{d + 1,a + 1}^{R}\) with the indicated twist and shift.
    Taking its fiber gives \eqref{eq:flag-cotwist-twist-comparison}.
\end{proof}

\begin{theorem}
    \label{thm:spherical-flags}
    For every \(d \geq 0\) and \(a \in \mathbb{Z}\), the functor
    \[
        \Psi_{d,a} \colon \Perf(\Gr_{X}(\mathcal{E},d + 1)) \longrightarrow \Perf(\Gr_{X}(\mathcal{E},d))
    \]
    is spherical.
\end{theorem}

\begin{proof}
    As in the proof of \zcref{lem:flag-squares}, the projective-line calculation gives \(\Psi_{d,a}\Psi_{d + 1,a + 1} \simeq 0\).
    Taking right adjoints yields \(\Psi_{d + 1,a + 1}^{R}\Psi_{d,a}^{R} \simeq 0\), and hence \(T_{d + 1,a + 1}\Psi_{d,a}^{R} \simeq \Psi_{d,a}^{R}\).
    By \zcref{lem:flag-cotwist} and relative Serre duality,
    \[
        \begin{aligned}
            C_{d,a}\Psi_{d,a}^{R} & \simeq \bigl(- \otimes \det \mathcal{E}^{\vee}\bigr) T_{d + 1,a + 1}\Psi_{d,a}^{R}[2d + 1] \\
                                  & \simeq \bigl(- \otimes \det \mathcal{E}^{\vee}\bigr) \Psi_{d,a}^{R}[2d + 1] \simeq \Psi_{d,a}^{L}.
        \end{aligned}
    \]
    Thus the natural equivalence required by \zcref{lem:spherical-rigidity} holds for every \(d\) and \(a\).

    We prove by induction on \(d\) that \(\Psi_{d,a}\) is spherical for every \(a \in \mathbb{Z}\).

    For \(d = 0\), let \(p \colon \mathbb{P}_{X}(\mathcal{E}) \to X\) be the projection.
    Then
    \[
        \Psi_{0,a}\Psi_{0,a}^{R} \simeq p_{*}p^{!}.
    \]
    The generalized Serre theorem \cite[Theorem~5.6(2iii)]{Jiang2025DerivedProjectivizations} gives the fiber sequence
    \[
        \det \mathcal{E}^{\vee} \longrightarrow \mathcal{O}_{X} \longrightarrow p_{*}\mathcal{O}_{\mathbb{P}_{X}(\mathcal{E})},
    \]
    where the second arrow is the adjunction unit.
    By relative Grothendieck duality and the projection formula, the dual of this unit induces the counit \(p_{*}p^{!} \to \id\).
    Dualizing the sequence therefore gives
    \[
        T_{0,a} \simeq - \otimes \det \mathcal{E}.
    \]
    In particular, \(T_{0,a}\) is an autoequivalence.
    Applying \zcref{lem:spherical-rigidity,thm:spherical-criterion}, we conclude that \(\Psi_{0,a}\) is spherical for every \(a\).

    Now let \(d \geq 1\) and suppose that \(\Psi_{d - 1,b}\) is spherical for every \(b \in \mathbb{Z}\).
    Fix \(a \in \mathbb{Z}\).
    By the induction hypothesis, \(C_{d - 1,a - 1}\) is an autoequivalence.
    The equivalence
    \[
        C_{d - 1,a - 1} \simeq \bigl(- \otimes \det \mathcal{E}^{\vee}\bigr) T_{d,a}[2d - 1]
    \]
    from \zcref{lem:flag-cotwist} shows that \(T_{d,a}\) is an autoequivalence.
    Together with the natural equivalence established above, \zcref{lem:spherical-rigidity,thm:spherical-criterion} imply that \(\Psi_{d,a}\) is spherical.
\end{proof}

\subsection{Incidence equivalences and flop--flop}
\label{subsec:incidence-twists}

We now compare the spherical twists of the flag functors with the incidence equivalences between the derived Grassmannians of \(\mathcal{E}\) and \(\mathcal{E}^{\vee}[1]\).
We recall the functors needed for this comparison, leaving the geometry of the incidence correspondences and the proof of their compatibility with the flag functors to \zcref{sec:integral-grassmannian-decompositions}.

For \(d \geq 0\), consider the incidence correspondence
\[
    \Gr_{X}(\mathcal{E},d) \xleftarrow{r_{+}} I_{d,d} \xrightarrow{r_{-}} \Gr_{X}(\mathcal{E}^{\vee}[1],d)
\]
of \zcref{def:incidence-correspondence}.
Write \(Q_{d}^{-}\) for the universal quotient on \(\Gr_{X}(\mathcal{E}^{\vee}[1],d)\), and put
\[
    \mathcal{O}(k,k) = r_{+}^{*}(\det Q_{d})^{k} \otimes r_{-}^{*}(\det Q_{d}^{-})^{k}.
\]
For \(k \in \mathbb{Z}\), define
\begin{equation}
    \label{eq:rank-zero-incidence-functors}
    \begin{aligned}
        \Phi_{d,k} \colon \Perf(\Gr_{X}(\mathcal{E}^{\vee}[1],d)) & \to \Perf(\Gr_{X}(\mathcal{E},d)),           & M & \mapsto r_{+,*}\bigl(r_{-}^{*}M \otimes \mathcal{O}(k,k)\bigr), \\
        \Phi_{d,k}^{\vee} \colon \Perf(\Gr_{X}(\mathcal{E},d))    & \to \Perf(\Gr_{X}(\mathcal{E}^{\vee}[1],d)), & N & \mapsto r_{-,*}\bigl(r_{+}^{*}N \otimes \mathcal{O}(k,k)\bigr).
    \end{aligned}
\end{equation}
The equivalence of \(\Phi_{d,0}\) follows from \cite[Lemma~3.11 and Corollary~3.12]{Jiang2023DerivedCategoriesGrassmannians}.
Applying the same result to \(\mathcal{E}^{\vee}[1]\), together with the projection formula, shows that all the functors in \eqref{eq:rank-zero-incidence-functors} are equivalences.

The following compatibility with the flag functors is proved in \zcref{subsec:rank-zero-incidence-comparison}.

\begin{proposition}
    \label{prop:flag-unit-fiber-equivalence}
    There is a natural equivalence
    \begin{equation}
        \label{eq:incidence-unit-fiber}
        T_{d,0}^{L}\Phi_{d,0} \simeq \bigl(- \otimes \det \mathcal{E}^{\vee}\bigr)\Phi_{d,1}.
    \end{equation}
\end{proposition}

\begin{theorem}
    \label{thm:flop-product}
    For every \(d \geq 0\) and \(k \in \mathbb{Z}\), there are natural equivalences of \(X\)-linear functors
    \begin{align}
        \label{eq:incidence-ratio}
        \Phi_{d,k}\Phi_{d,k + 1}^{-1}                & \simeq \bigl(- \otimes \det \mathcal{E}^{\vee}\bigr)T_{d,-k}, \\
        \label{eq:flop-product}
        \bigl(\Phi_{d,0}\Phi_{d,0}^{\vee}\bigr)^{-1} & \simeq T_{d,d}\cdots T_{d,2}T_{d,1}.
    \end{align}
    For \(d = 0\), the product in the second formula is empty.
\end{theorem}

\begin{proof}
    Since \(\Psi_{d,0}\) is spherical, \(T_{d,0}^{L} \simeq T_{d,0}^{-1}\).
    Taking inverses in \eqref{eq:incidence-unit-fiber} gives
    \[
        \Phi_{d,0}\Phi_{d,1}^{-1} \simeq \bigl(- \otimes \det \mathcal{E}^{\vee}\bigr)T_{d,0}.
    \]
    The projection formula gives
    \[
        \Phi_{d,k} \simeq \bigl(- \otimes (\det Q_{d})^{k}\bigr) \Phi_{d,0} \bigl(- \otimes (\det Q_{d}^{-})^{k}\bigr).
    \]
    Similarly,
    \[
        T_{d,-k} \simeq \bigl(- \otimes (\det Q_{d})^{k}\bigr) T_{d,0} \bigl(- \otimes (\det Q_{d})^{-k}\bigr).
    \]
    Conjugating the formula for \(k = 0\) proves \eqref{eq:incidence-ratio}.

    We next compute the inverse of the incidence equivalence.
    By \cite[Lemma~2.10]{Jiang2023DerivedCategoriesGrassmannians}, \(r_{+}\) is the Grassmannian of rank \(d\) quotients of \(R_{d}^{\vee}[1]\).
    Since
    \[
        \det (R_{d}^{\vee}[1]) \simeq \det R_{d} \simeq \det \mathcal{E} \otimes (\det Q_{d})^{-1},
    \]
    Grothendieck duality \cite[Theorem~4.18 and Corollary~4.47]{Jiang2023DerivedGrassmanniansSchurFunctors} gives
    \[
        \omega_{r_{+}} \simeq (\det \mathcal{E})^{d} \otimes \mathcal{O}(-d,-d).
    \]
    Computing the right adjoint of \(\Phi_{d,k}\) therefore yields
    \begin{equation}
        \label{eq:inverse-incidence}
        \Phi_{d,k}^{-1} \simeq \bigl(- \otimes (\det \mathcal{E})^{d}\bigr) \Phi_{d,-d - k}^{\vee}.
    \end{equation}
    In particular,
    \[
        \bigl(\Phi_{d,0}\Phi_{d,0}^{\vee}\bigr)^{-1} \simeq \bigl(- \otimes (\det \mathcal{E})^{d}\bigr) \Phi_{d,-d}\Phi_{d,0}^{-1}.
    \]
    Applying \eqref{eq:incidence-ratio} successively for \(k = -d,\ldots,-1\) gives
    \[
        \Phi_{d,-d}\Phi_{d,0}^{-1} \simeq \bigl(- \otimes (\det \mathcal{E}^{\vee})^{d}\bigr) T_{d,d}\cdots T_{d,1}.
    \]
    Since the twists are \(X\)-linear, the determinant factors cancel, proving \eqref{eq:flop-product}.
    The case \(d = 0\) follows directly from the identity incidence correspondence.
\end{proof}

\begin{remark}
    \label{rem:d=1}
    For \(d = 1\), our formula specializes to the flop--flop formula for derived projectivizations studied by Jiang and Leung \cite[Theorem~4.1(iii) and Example~4.3]{JiangLeung2022ProjectivizationsFlops}.
    Let \(p \colon \mathbb{P}_{X}(\mathcal{E}) \to X\) be the projection.
    Since \(\Psi_{0,0} = p_{*}\) and \(\Psi_{0,0}^{L} = p^{*}\), we have \(C_{0,0} \simeq T_{p^{*}}[-1]\), where \(T_{p^{*}}\) denotes the spherical twist of \(p^{*}\).
    Applying \eqref{eq:flag-cotwist-twist-comparison} with \(d = a = 0\), we obtain
    \[
        \bigl(\Phi_{1,0}\Phi_{1,0}^{\vee}\bigr)^{-1} \simeq T_{1,1} \simeq \bigl(- \otimes \det \mathcal{E}\bigr)T_{p^{*}}[-2].
    \]
    Equivalently,
    \[
        \Phi_{1,0}\Phi_{1,0}^{\vee} \simeq \bigl(- \otimes \det \mathcal{E}^{\vee}\bigr)T_{p^{*}}^{-1}[2].
    \]

    In the classical Atiyah flop, \(\Psi_{1,1}\) is the spherical functor associated with \(i_{*}\mathcal{O}_{\mathbb{P}^{1}}(-1)\), where \(i \colon \mathbb{P}^{1} \hookrightarrow Y\) is the exceptional curve.
    Thus the inverse flop--flop is the spherical twist \(T_{i_{*}\mathcal{O}_{\mathbb{P}^{1}}(-1)}\).
\end{remark}

\begin{example}
    \label{ex:DS}
    For \(X = \operatorname{Hom}(W,V)\), where \(W\) and \(V\) are \(n\)-dimensional vector spaces over \(\mathbb{C}\), and \(\mathcal{E}\) the tautological complex \(W \otimes \mathcal{O}_{X} \to V \otimes \mathcal{O}_{X}\), the spaces \(\Gr_{X}(\mathcal{E},d)\) and \(\Gr_{X}(\mathcal{E}^{\vee}[1],d)\) are the two sides of the local Grassmannian flop of \cite{DonovanSegal2014WindowShiftsGrassmannianTwists}.
    By \cite[Theorems~3 and~3.2.21]{BallardEtAl2021KernelsGrassmannFlops}, the inverse window equivalences of \cite{DonovanSegal2014WindowShiftsGrassmannianTwists} are given by the incidence transforms \(\Phi_{d,k}\), since in this case \(I_{d,d}\) coincides with the classical fiber product over \(X\).
    Hence \eqref{eq:incidence-ratio} recovers the identification of window shifts with spherical twists in \cite{DonovanSegal2014WindowShiftsGrassmannianTwists}, and \eqref{eq:flop-product} factors the inverse of the flop--flop autoequivalence into \(d\) twists.
\end{example}

%% file: src/abel-applications.tex
\section{Applications to linear series on curves}
\label{sec:abel-applications}

We apply our semiorthogonal decompositions and flop--flop formulas to generalized Abel maps for families of integral curves.
Using the derived Grassmannian descriptions of \cite{Jiang2025AbelMapsIntegralCurves}, we obtain semiorthogonal decompositions in arbitrary characteristic for positive-rank complexes and nontrivial autoequivalences in rank zero.

\subsection{Generalized linear series and categorical decompositions}

Let \(\pi:\mathcal{C} \to S\) be a flat, locally projective, finitely presented family of geometrically integral curves of arithmetic genus \(g\).
Let \(X \to S\) be a derived stack equipped with a family \(\mathcal{J}\) of rank one torsion-free sheaves of degree \(d\) on \(\mathcal{C}_{X} = \mathcal{C} \times_{S} X\), perfect relative to \(X\).
Write \(f:\mathcal{C}_{X} \to X\), and put
\[
    \mathcal{E} = (f_{*}\mathcal{J})^{\vee}.
\]
By \cite[Corollary~2.10 and Section~4.1]{Jiang2025AbelMapsIntegralCurves}, the complex \(\mathcal{E}\) is perfect of Tor-amplitude in \([0,1]\), of rank
\[
    r = d + 1 - g.
\]
For \(s \geq 0\), write
\begin{equation}
    \label{eq:abel-grassmannians}
    \mathbb{G}_{s} = \Gr_{X}(\mathcal{E},s), \qquad \mathbb{G}_{s}^{\vee} = \Gr_{X}(\mathcal{E}^{\vee}[1],s).
\end{equation}

The spaces \(\mathbb{G}_{s}\) parametrize \(s\)-dimensional linear systems of sections, while \(\mathbb{G}_{s}^{\vee}\) describe their Serre-dual counterparts.
Their structural morphisms to \(X\) are generalized Abel maps in the sense of \cite[Section~4]{Jiang2025AbelMapsIntegralCurves}.

For example, if \(C\) is a smooth projective curve, taking \(X = \Pic^{d}(C)\) with its Poincar\'e line bundle recovers the classical Abel--Jacobi map for \(s = 1\):
\[
    C^{(d)} \longrightarrow \Pic^{d}(C), \qquad D \longmapsto \mathcal{O}_{C}(D).
\]

Applying our Fock decompositions and flop--flop formulas from \zcref{sec:clifford-flags,sec:spherical-flags,sec:integral-grassmannian-decompositions} gives the following result.

\begin{theorem}
    \label{thm:abel-categorical-applications}
    With the notation above, the following hold.
    \begin{enumerate}

        \item If \(r > 0\), then for every \(s \geq 0\) there is an \(X\)-linear semiorthogonal decomposition
              \[
                  \Perf(\mathbb{G}_{s}) = \left\langle \operatorname{Im}(F_{I}\Phi_{\varnothing}) \ \middle|\ I \subseteq [r],\ s + |I| \geq r \right\rangle,
              \]
              where the fully faithful embedding functors are
              \[
                  F_{I}\Phi_{\varnothing}: \Perf(\mathbb{G}^{\vee}_{s + |I| - r}) \longrightarrow \Perf(\mathbb{G}_{s}).
              \]
        \item If \(r = 0\), then every flag functor \(\Psi_{s,m}\) is spherical.
              The incidence functors \(\Phi_{s,k}\) and \(\Phi_{s,k}^{\vee}\) of \eqref{eq:rank-zero-incidence-functors} are equivalences, and
              \begin{equation}
                  \label{eq:critical-abel-flop}
                  \bigl(\Phi_{s,0}\Phi_{s,0}^{\vee}\bigr)^{-1} \simeq T_{s,s}\cdots T_{s,1},
              \end{equation}
              with the empty product understood as the identity for \(s = 0\).
    \end{enumerate}
\end{theorem}

When \(X = \overline{\Pic}^{\,d}_{\mathcal{C}/S}\) is equipped with a universal rank one torsion-free sheaf \(\mathcal{J}\), for example after rigidification along a section through the smooth locus, the derived Grassmannians \(\mathbb{G}_{s}\) parametrize \(s\)-dimensional linear systems on these sheaves \cite[Sections~4.2--4.3]{Jiang2025AbelMapsIntegralCurves}.
If \(\mathcal{C}/S\) is Gorenstein, \(\mathbb{G}_{1}\) is the derived Hilbert scheme of \(d\) points.
For locally planar families over a classical base, this derived Hilbert scheme is itself classical \cite[Corollary~4.6 and Lemma~4.7]{Jiang2025AbelMapsIntegralCurves}.
In particular, for a smooth curve \(C\), taking \(X = \Pic^{d}(C)\) recovers the symmetric power \(C^{(d)}\) when \(s = 1\), while the higher \(\mathbb{G}_{s}\) are derived enhancements of the varieties of linear series \(G_{d}^{s - 1}(C)\).

For \(d \geq g\), the case \(s = 1\) is covered by the semiorthogonal decomposition for derived projectivizations, which holds in arbitrary characteristic \cite[Theorem~4.4]{Jiang2025AbelMapsIntegralCurves}.
For \(s > 1\), our theorem removes the characteristic-zero assumption in the decompositions of higher linear series \cite[Corollary~4.9(2)]{Jiang2025AbelMapsIntegralCurves}.
At the critical degree \(d = g - 1\), the case \(s = 1\) for smooth curves recovers the classical \(\Theta\)-flop associated with the Abel--Jacobi maps, studied by Toda and Jiang--Leung \cite{Toda2021StablePairModuliSpaces,JiangLeung2022ProjectivizationsFlops}.
Our rank-zero results extend this picture to higher linear-series spaces, providing spherical flag functors and flop--flop formulas.

\subsection{Brill--Noether loci and autoequivalences}

We now study the connection between our rank-zero flop--flop formulas and classical Brill--Noether geometry.
In this subsection, we further assume that \(S\) and \(X\) are classical schemes and that \(\mathcal{J}\) is an \(X\)-flat coherent sheaf.
We specialize to \(d = g - 1\), so that \(r = 0\).

For \(s \geq 0\), let \(a_{s}:\mathbb{G}_{s} \to X\) be the projection, and consider the Brill--Noether loci
\[
    B_{s} = \{x \in X:h^{0}(\mathcal{C}_{x},\mathcal{J}_{x}) \geq s\}.
\]
These are closed determinantal loci, with
\[
    X = B_{0} \supseteq B_{1} \supseteq B_{2} \supseteq \cdots,
\]
and the image of \(\mathbb{G}_{s} \to X\) is \(B_{s}\) on geometric points.

The spherical twists in the flop--flop formula \eqref{eq:critical-abel-flop} arise from the flag functors \(\Psi_{s,m}\), which relate linear systems of dimensions \(s + 1\) and \(s\).
We show that their nontriviality is controlled by the next Brill--Noether locus \(B_{s + 1}\).

\begin{proposition}
    \label{prop:abel-window-changes}
    For every \(s \geq 0\) and \(m \in \mathbb{Z}\), there is a fiber sequence of endofunctors of \(\Perf(\mathbb{G}_{s})\)
    \[
        \Psi_{s,m}\Psi_{s,m}^{R} \longrightarrow \id \longrightarrow T_{s,m},
    \]
    where the first term takes values in perfect complexes supported over \(B_{s + 1}\).
    Consequently, the spherical twists in \eqref{eq:critical-abel-flop}, and hence the flop--flop autoequivalence, become the identity over \(X \setminus B_{s + 1}\).
    If \(s \geq 1\) and \(B_{s + 1} \neq \varnothing\), then \(T_{s,m}\) has infinite order on \(\Perf(\mathbb{G}_{s})\).
\end{proposition}

\begin{proof}
    The fiber sequence is the defining triangle of the spherical twist.
    Since the structural map \(\mathbb{G}_{s + 1} \to X\) has image \(B_{s + 1}\), every object in the image of \(\Psi_{s,m}\) is supported over this locus.
    This proves the first assertion.

    For the second, the spherical adjunction identities and \zcref{lem:flag-cotwist} give
    \begin{equation}
        \label{eq:abel-window-shift}
        T_{s,m}^{-1}\Psi_{s,m} \simeq \Psi_{s,m}(- \otimes a_{s + 1}^{*}\det\mathcal{E}^{\vee})[2s],
    \end{equation}
    using \(\Psi_{s,m}\Psi_{s + 1,m + 1} \simeq 0\).
    If \(B_{s + 1} \neq \varnothing\), the projective-bundle formula implies \(\Psi_{s,m} \neq 0\).
    Thus, for some nonzero \(A \in \operatorname{Im}\Psi_{s,m}\), we have
    \[
        T_{s,m}^{-q}(A) \simeq A \otimes a_{s}^{*}(\det\mathcal{E}^{\vee})^{\otimes q}[2sq]
    \]
    for every \(q \geq 1\).
    Taking a nonzero geometric fiber of \(A\) shows that \(T_{s,m}^{q}\not\simeq \id\) for \(s \geq 1\).
\end{proof}

The proposition shows that the spherical corrections in the flop--flop formula are supported over the deeper Brill--Noether loci.
We now consider a case where the last nonempty locus consists of finitely many points, allowing us to describe the corresponding spherical functors explicitly.

Let \(C/\mathbb{C}\) be a smooth projective Brill--Noether--Petri general curve of genus \(g = n^{2}\), with \(n \geq 2\), and take \(X = \Pic^{g - 1}(C)\) with its Poincar\'e line bundle.
For \(0 \leq s \leq n\), the derived Grassmannians \(\mathbb{G}_{s}\) are classical, smooth projective varieties of dimension
\[
    \dim\mathbb{G}_{s} = g - s^{2},
\]
while \(\mathbb{G}_{s} = \varnothing\) for \(s > n\), by the Brill--Noether and Gieseker--Petri theorems \cite[V~(1.5) and VII~(2.3)]{ArbarelloEtAl1985GeometryAlgebraicCurvesI}\cite[Theorem~1.1]{Gieseker1982StableCurvesSpecialDivisors} and \cite[Lemma~4.8]{Jiang2025AbelMapsIntegralCurves}.

In particular, \(\mathbb{G}_{n}\) consists of \(N_{n}\) reduced points corresponding to line bundles \(L_{\alpha}\), where \(N_{n} = (n^{2})!\,\prod_{i = 0}^{n - 1}\bigl(i!/(n + i)!\bigr)\) is the Castelnuovo number \cite[Main Theorem~II(b)]{GriffithsHarris1980SpecialLinearSystems}.
Their fibers in \(Y = \mathbb{G}_{n - 1}\) give pairwise disjoint closed embeddings
\[
    i_{\alpha}:P_{\alpha} \simeq \mathbb{P}^{n - 1} \hookrightarrow Y.
\]

\begin{theorem}
    \label{cor:square-genus-flop-product}
    For \(m \in \mathbb{Z}\), put
    \[
        S_{\alpha,m} = i_{\alpha,*}\mathcal{O}_{P_{\alpha}}(-m) \in \Perf(Y).
    \]
    For each fixed \(m\), these are mutually orthogonal \((2n - 1)\)-spherical objects, and
    \begin{equation}
        \label{eq:square-genus-weighted-twists}
        T_{n - 1,m} \simeq \prod_{\alpha = 1}^{N_{n}}T_{S_{\alpha,m}}.
    \end{equation}
    In particular,
    \begin{equation}
        \label{eq:square-genus-flop-factors}
        \bigl(\Phi_{n - 1,0}\Phi_{n - 1,0}^{\vee}\bigr)^{-1} \simeq \left(\prod_{\alpha} T_{S_{\alpha,n - 1}}\right)\cdots \left(\prod_{\alpha} T_{S_{\alpha,1}}\right).
    \end{equation}
    Moreover, writing \(S_{\alpha} = S_{\alpha,1}\), the twists \(T_{S_{\alpha}}\) generate a free abelian subgroup
    \[
        \mathbb{Z}^{N_{n}} \subseteq \Aut(\Perf(Y)).
    \]
\end{theorem}

\begin{proof}
    The flag correspondence between \(\mathbb{G}_{n}\) and \(Y\) is \(\coprod_{\alpha} P_{\alpha}\), with universal line restricting to \(\mathcal{O}_{P_{\alpha}}(-1)\).
    Thus \(\Psi_{n - 1,m}\) decomposes into the functors defined by the objects \(S_{\alpha,m}\).
    Since \(\mathbb{G}_{n + 1} = \varnothing\), \zcref{lem:flag-cotwist} gives
    \[
        C_{n - 1,m} \simeq (- \otimes \det\mathcal{E}^{\vee})[2n - 1].
    \]
    By \zcref{thm:spherical-flags} and the disjointness of the supports, each summand is spherical with cotwist \([2n - 1]\).
    Hence the objects \(S_{\alpha,m}\) are mutually orthogonal and \((2n - 1)\)-spherical.

    The counit of \(\Psi_{n - 1,m}\) is the direct sum of the evaluation maps for these objects.
    Orthogonality therefore identifies its spherical twist with the product of the individual twists, proving \eqref{eq:square-genus-weighted-twists}.
    The flop--flop factorization follows from \eqref{eq:critical-abel-flop}.

    Finally,
    \[
        T_{S_{\alpha}}(S_{\alpha}) \simeq S_{\alpha}[2 - 2n], \qquad T_{S_{\beta}}(S_{\alpha}) \simeq S_{\alpha} \quad (\beta \neq \alpha).
    \]
    The twists commute by orthogonality.
    Any relation \(\prod_{\alpha} T_{S_{\alpha}}^{a_{\alpha}} \simeq \id\) would therefore imply
    \[
        S_{\alpha}[(2 - 2n)a_{\alpha}] \simeq S_{\alpha}
    \]
    for every \(\alpha\).
    Since \(n \geq 2\), all \(a_{\alpha}\) must vanish.
\end{proof}

The spherical twists constructed above are independent as autoequivalences, but their differences are invisible to rational cohomology and numerical \(K\)-theory.

\begin{proposition}
    \label{prop:square-genus-cohomological-action}
    Fix an index \(\alpha_{0}\).
    The autoequivalences
    \[
        T_{S_{\alpha}}T_{S_{\alpha_{0}}}^{-1}, \qquad \alpha \neq \alpha_{0},
    \]
    freely generate a subgroup \(\mathbb{Z}^{N_{n} - 1}\) of \(\Aut(\Perf(Y))\) acting trivially on rational cohomology and numerical \(K\)-theory.
\end{proposition}

\begin{proof}
    We first show that \(\operatorname{ch}(S_{\alpha})\) is independent of \(\alpha\).

    Write \(J = \Pic^{g - 1}(C)\), and let \(p:C \times J \to J\) be the projection, with Poincar\'e bundle \(\mathcal{P}\).
    Choose a sufficiently positive effective divisor \(D\) on \(C\).
    The restriction sequence for \(\mathcal{P}(D)\) gives a two-term presentation
    \[
        p_{*}\mathcal{P} \simeq \left[p_{*}\mathcal{P}(D) \longrightarrow p_{*}(\mathcal{P}(D)|_{D \times J})\right] = :[E_{0} \to E_{1}].
    \]
    Both bundles have the same rank, and \(\underline{\Hom}(E_{0},E_{1})\) is ample by \cite[Lemma~2.2]{FultonLazarsfeld1981ConnectednessDegeneracyLoci}.
    The dual complex realizes \(Y\) as a zero locus in \(G = \Gr_{J}(E_{0}^{\vee},n - 1)\).
    By \cite[Proposition~1.2 and Remark~1.9]{FultonLazarsfeld1981ConnectednessDegeneracyLoci} and Lefschetz duality, the inclusion induces isomorphisms
    \begin{equation}
        \label{eq:square-genus-lefschetz-comparison}
        H_{i}(Y,\mathbb{Q})\xrightarrow{\sim}H_{i}(G,\mathbb{Q}), \qquad i < 2n - 1.
    \end{equation}

    The relative flag variety \(\Flag_{J}(E_{0}^{\vee};n - 1,n)\) is a projective bundle over the connected Grassmannian \(\Gr_{J}(E_{0}^{\vee},n)\).
    Its fibers over the points corresponding to \(L_{\alpha}\) map isomorphically onto \(P_{\alpha}\), with the relative hyperplane class restricting to \(h_{\alpha} = c_{1}(\mathcal{O}_{P_{\alpha}}(1))\).
    Hence the classes
    \[
        i_{\alpha,*}(h_{\alpha}^{a} \cap [P_{\alpha}]), \qquad 0 \leq a \leq n - 1,
    \]
    have the same images in \(H_{*}(G,\mathbb{Q})\).
    Since their degrees are at most \(2n - 2\), \eqref{eq:square-genus-lefschetz-comparison} shows that they agree in \(H_{*}(Y,\mathbb{Q})\).

    Grothendieck--Riemann--Roch gives
    \[
        \operatorname{ch}(S_{\alpha})\operatorname{td}(Y) = i_{\alpha,*}\bigl(e^{-h_{\alpha}}\operatorname{td}(P_{\alpha})\bigr).
    \]
    Since \(\operatorname{td}(P_{\alpha})\) is a polynomial in \(h_{\alpha}\), the preceding comparison proves that \(\operatorname{ch}(S_{\alpha})\) is independent of \(\alpha\).

    A spherical twist acts on \(K_{0}(Y)\) by
    \[
        [T_{S_{\alpha}}(M)] = [M] - \chi(S_{\alpha},M)[S_{\alpha}].
    \]
    By Riemann--Roch, the equality of Chern characters implies that the twists \(T_{S_{\alpha}}\) induce the same action on numerical \(K\)-theory.
    Their actions on rational cohomology also agree by Grothendieck--Riemann--Roch.
    Thus the autoequivalences \(T_{S_{\alpha}}T_{S_{\alpha_{0}}}^{-1}\) act trivially on both groups, while their independence follows from the preceding theorem.
\end{proof}

%% file: src/incidence-and-extended-clifford-relations.tex
\section{Incidence correspondences and extended Clifford relations}
\label{sec:integral-grassmannian-decompositions}

In this appendix, we review the incidence geometry of derived Grassmannians and use it to establish the incidence--flag comparison needed for the rank-zero flop--flop formula.
In positive rank, we establish categorical analogues of the extended Clifford relations of Jiang, Li, and Zhao \cite{JiangEtAl2025InstantonsAffineVertexAlgebra}.
Using these relations, we identify the vacuum category and compare our semiorthogonal decompositions with those of \cite{Jiang2023DerivedCategoriesGrassmannians}.

\subsection{Incidence geometry and the flop--flop comparison}
\label{subsec:incidence-flag-geometry}
\label{subsec:rank-zero-incidence-comparison}

Let \(X\) be a prestack and let \(\mathcal{E}\) be a perfect complex of rank \(r\) and Tor-amplitude in \([0,1]\) over \(X\).
For \(d,e \geq 0\), consider the derived Grassmannians \(\Gr_{X}(\mathcal{E},d)\) and \(\Gr_{X}(\mathcal{E}^{\vee}[1],e)\), with universal fiber sequences
\[
    R_{d}^{+} \longrightarrow \mathcal{E} \longrightarrow Q_{d}^{+}, \qquad R_{e}^{-} \longrightarrow \mathcal{E}^{\vee}[1] \longrightarrow Q_{e}^{-}.
\]

\begin{definition}
    \label{def:incidence-correspondence}
    For \(d,e \geq 0\), define
    \[
        I_{d,e} = \Incid_{X}(\mathcal{E};d,e) : = \Gr_{\Gr_{X}(\mathcal{E},d)} \bigl((R_{d}^{+})^{\vee}[1],e\bigr).
    \]
    The universal quotient determines the incidence correspondence
    \[
        \Gr_{X}(\mathcal{E},d) \xleftarrow{r_{+}}I_{d,e} \xrightarrow{r_{-}}\Gr_{X}(\mathcal{E}^{\vee}[1],e).
    \]
    Equivalently, there is a canonical identification
    \[
        I_{d,e} \simeq \Gr_{\Gr_{X}(\mathcal{E}^{\vee}[1],e)} \bigl((R_{e}^{-})^{\vee}[1],d\bigr)
    \]
    by \cite[Lemma~2.10]{Jiang2023DerivedCategoriesGrassmannians}.

\end{definition}
On \(I_{d,e}\), let \(\mathcal{U}_{d,e}\) be the perfect complex defined by the fiber sequence
\begin{equation}
    \label{eq:incidence-universal-complex}
    (Q_{e}^{-})^{\vee}[1] \longrightarrow R_{d}^{+} \longrightarrow \mathcal{U}_{d,e}.
\end{equation}
It has Tor-amplitude in \([0,1]\) and rank \(r - d + e\).

\begin{proposition}
    \label{prop:incidence-flag-geometry}
    For \(d,e \geq 0\) and \(m \geq 1\), there is a canonical equivalence over \(I_{d,e}\)
    \begin{equation}
        \label{eq:incidence-flag-geometry}
        \begin{aligned}
             & \Flag_{X}(\mathcal{E};d,d + 1,\ldots,d + m) \times_{\Gr_{X}(\mathcal{E},d + m)}I_{d + m,e} \\
             & \qquad \simeq \Flag_{I_{d,e}}(\mathcal{U}_{d,e};1,\ldots,m).
        \end{aligned}
    \end{equation}
\end{proposition}

\begin{proof}
    Fix a point of \(I_{d,e}\).
    A flag on the left determines a chain of quotients of \(\mathcal{U}_{d,e}\), obtained from the fibers of \(Q_{d + v} \twoheadrightarrow Q_{d}\) via \eqref{eq:incidence-universal-complex}.
    Conversely, such a flag of quotients of \(\mathcal{U}_{d,e}\) reconstructs the original flag of \(\mathcal{E}\) by pushouts.
    These constructions are mutually inverse and natural, as in \cite[Proposition~2.12]{Jiang2023DerivedCategoriesGrassmannians}.
\end{proof}

We now assume \(r = 0\).
The preceding incidence--flag equivalence, together with the generalized Serre theorem for derived projectivizations, gives a fiber sequence relating the incidence and flag functors.

\begin{lemma}
    \label{lem:incidence-projectivization}
    For every \(d \geq 0\), there is a canonical fiber sequence of functors
    \begin{equation}
        \label{eq:incidence-functor-triangle}
        (- \otimes \det\mathcal{E}^{\vee})\Phi_{d,1} \longrightarrow \Phi_{d,0} \xrightarrow{v_{d}}\Psi_{d,0}\Phi_{(d + 1,d)},
    \end{equation}
    where \(\Phi_{(d + 1,d)} = r_{+,*}r_{-}^{*}\) is defined by the incidence correspondence \(I_{d + 1,d}\).
\end{lemma}

\begin{proof}
    By \zcref{prop:incidence-flag-geometry} with \(e = d\) and \(m = 1\), we have
    \[
        \mathbb{P}_{I_{d,d}}(\mathcal{U}_{d,d}) \simeq \Flag_{X}(\mathcal{E};d,d + 1) \times_{\Gr_{X}(\mathcal{E},d + 1)}I_{d + 1,d}.
    \]
    Since \(\mathcal{U}_{d,d}\) has rank zero, the generalized Serre theorem gives a canonical fiber sequence
    \[
        \det\mathcal{U}_{d,d}^{\vee} \longrightarrow \mathcal{O}_{I_{d,d}} \longrightarrow q_{*}\mathcal{O}_{\mathbb{P}_{I_{d,d}}(\mathcal{U}_{d,d})},
    \]
    where \(q:\mathbb{P}_{I_{d,d}}(\mathcal{U}_{d,d}) \to I_{d,d}\) is the projection.
    The universal fiber sequence \eqref{eq:incidence-universal-complex} gives
    \[
        \det\mathcal{U}_{d,d}^{\vee} \simeq r_{+}^{*}\det\mathcal{E}^{\vee} \otimes r_{+}^{*}\det Q_{d} \otimes r_{-}^{*}\det Q_{d}^{-}.
    \]
    Tensoring the fiber sequence by pullbacks from \(\Gr_{X}(\mathcal{E}^{\vee}[1],d)\) and applying \(r_{+,*}\), the projection formula and base change identify its three terms with the functors in \eqref{eq:incidence-functor-triangle}.
    The map \(v_{d}\) is induced by the adjunction unit \(\id \to q_{*}q^{*}\).
\end{proof}

To identify the fiber sequence \eqref{eq:incidence-functor-triangle} with the one defining the spherical twist, we need to identify its natural transformation \(v_{d}\).
Under the adjunction \(\Psi_{d,0}^{L}\dashv\Psi_{d,0}\), this transformation corresponds to
\begin{equation}
    \label{eq:incidence-adjunction}
    b_{d}:\Psi_{d,0}^{L}\Phi_{d,0} \longrightarrow \Phi_{(d + 1,d)}.
\end{equation}

\begin{proposition}
    \label{prop:incidence-comparison}
    For every \(d \geq 0\), the natural transformation \(b_{d}\) is an equivalence.
\end{proposition}

\begin{proof}
    The natural transformation \(b_{d}\) is an equivalence if and only if its cofiber is the zero functor.
    Since this vanishing is fpqc-local on \(X\), we may assume that \(X\) is affine and write \(\mathcal{E} \simeq [A \to B]\), where \(A\) and \(B\) are vector bundles with \(B\) in degree \(0\).

    We first identify the Fourier--Mukai kernels of the two functors in \eqref{eq:incidence-adjunction} and the natural map between them.
    Let \(p\) be the base change of
    \[
        \pi_{d}:\Flag_{X}(\mathcal{E};d,d + 1) \longrightarrow \Gr_{X}(\mathcal{E},d + 1)
    \]
    along \(\Gr_{X}(A^{\vee},d) \to X\).
    Then \(p\) is a projective bundle of relative dimension \(d\).
    Let \(Q\) be the universal quotient on \(\Gr_{X}(A^{\vee},d)\), and write \(K_{j}\) for the kernel of the universal rank-\(j\) quotient of \(B\).
    The flag correspondence gives an exact sequence
    \begin{equation}
        \label{eq:incidence-koszul-flag-sequence}
        0 \longrightarrow p^{*}K_{d + 1} \longrightarrow K_{d} \longrightarrow L_{d + 1} \longrightarrow 0, \qquad L_{d + 1} \simeq \mathcal{O}_{p}(-1).
    \end{equation}

    The derived Grassmannian condition gives a map \(A \to K_{d + 1}\).
    Composing with \(Q^{\vee} \hookrightarrow A\) yields a map
    \[
        \sigma:Q^{\vee} \longrightarrow K_{d + 1},
    \]
    whose pullback and composition with \(p^{*}K_{d + 1} \to K_{d}\) give
    \[
        \sigma^{\prime}:Q^{\vee} \longrightarrow K_{d}.
    \]
    The incidence correspondence \(I_{d + 1,d}\) is the derived zero locus \(Z(\sigma)\), while the pullback of \(I_{d,d}\) to the flag space is \(Z(\sigma^{\prime})\).
    Consequently, the kernels of \(\Phi_{(d + 1,d)}\) and \(\Psi_{d,0}^{L}\Phi_{d,0}\) are represented by
    \[
        \mathcal{O}_{Z(\sigma)} \qquad \text{and} \qquad p_{!}\mathcal{O}_{Z(\sigma^{\prime})},
    \]
    respectively.

    Since \(\sigma^{\prime}\) factors through \(p^{*}K_{d + 1}\), there is a canonical morphism
    \[
        Z(p^{*}\sigma) \longrightarrow Z(\sigma^{\prime}),
    \]
    inducing a morphism of structure sheaves
    \[
        \mathcal{O}_{Z(\sigma^{\prime})} \longrightarrow p^{*}\mathcal{O}_{Z(\sigma)}.
    \]
    By \zcref{prop:incidence-flag-geometry}, the source \(Z(p^{*}\sigma)\) is identified with \(\mathbb{P}_{I_{d,d}}(\mathcal{U}_{d,d})\).
    Under this identification, the preceding morphism of derived zero loci is the map used to define \(v_{d}\) in \zcref{lem:incidence-projectivization}.
    Taking its adjoint along the flag projection therefore identifies \(b_{d}\) with
    \begin{equation}
        \label{eq:incidence-koszul-mate}
        p_{!}\mathcal{O}_{Z(\sigma^{\prime})} \longrightarrow p_{!}p^{*}\mathcal{O}_{Z(\sigma)} \longrightarrow \mathcal{O}_{Z(\sigma)},
    \end{equation}
    where the second arrow is the adjunction counit.

    It remains to prove that this composite is an equivalence.
    By \eqref{eq:incidence-koszul-flag-sequence}, the fiber of
    \[
        \mathcal{O}_{Z(\sigma^{\prime})} \longrightarrow p^{*}\mathcal{O}_{Z(\sigma)}
    \]
    admits a finite filtration with graded pieces
    \[
        p^{*}\mathcal{O}_{Z(\sigma)} \otimes \bigwedge^{a}Q^{\vee} \otimes L_{d + 1}^{-a}[a], \qquad 1 \leq a \leq d.
    \]
    Since \(p\) is a projective bundle of relative dimension \(d\), we have
    \[
        p_{!}L_{d + 1}^{-a} \simeq 0 \quad (1 \leq a \leq d), \qquad p_{!}p^{*} \simeq \id.
    \]
    The projection formula shows that the first map in \eqref{eq:incidence-koszul-mate} is an equivalence, while the second is an equivalence by adjunction.
    Hence \(b_{d}\) is an equivalence.
\end{proof}

The equivalence in \zcref{prop:incidence-comparison} allows us to identify the natural transformation \(v_{d}\), and hence the fiber sequence \eqref{eq:incidence-functor-triangle}, with the adjunction unit defining the left spherical twist.

\begin{proof}[Proof of \zcref{prop:flag-unit-fiber-equivalence}]
    By the definition of \(b_{d}\) under the adjunction \(\Psi_{d,0}^{L}\dashv\Psi_{d,0}\), the morphism \(v_{d}\) factors as
    \[
        \Phi_{d,0} \xrightarrow{\eta_{d,0}^{L}\Phi_{d,0}} \Psi_{d,0}\Psi_{d,0}^{L}\Phi_{d,0} \xrightarrow{\Psi_{d,0}b_{d}} \Psi_{d,0}\Phi_{(d + 1,d)}.
    \]
    Since \(b_{d}\) is an equivalence by \zcref{prop:incidence-comparison}, the fiber of \(v_{d}\) is canonically equivalent to \(T_{d,0}^{L}\Phi_{d,0}\).
    Comparing with \eqref{eq:incidence-functor-triangle}, we obtain
    \[
        T_{d,0}^{L}\Phi_{d,0} \simeq (- \otimes \det\mathcal{E}^{\vee})\Phi_{d,1},
    \]
    proving the assertion.
\end{proof}

This completes the incidence--flag comparison in rank zero, which gives the relation between incidence equivalences and spherical twists used in the proof of \zcref{thm:flop-product}.

\subsection{The vacuum category and extended Clifford relations}
\label{subsec:fock-incidence-identification}

Let \(r > 0\), and let \(\mathrm{Cl}_{r}\) be the standard Clifford algebra over a commutative ring \(R\), with \(r\) pairs of annihilation and creation operators \(p_{i},q_{i}\), respectively.
For any \(\mathrm{Cl}_{r}\)-module \(M\), the Fock decomposition gives an isomorphism
\[
    M \simeq V \otimes_{R}\bigwedge\nolimits_{R}^{\bullet}R^{r}, \qquad V = \bigcap_{i = 0}^{r - 1}\ker p_{i},
\]
where \(V\) is the vacuum module.

A natural question is how to identify this vacuum module with a more concrete module \(N\).
As observed by Jiang, Li, and Zhao \cite[Proposition~3.3 and \S5]{JiangEtAl2025InstantonsAffineVertexAlgebra}, such an identification is equivalent to giving maps
\[
    e:N \longrightarrow M, \qquad f:M \longrightarrow N
\]
satisfying the extended Clifford relations
\begin{equation}
    \label{eq:extended-clifford-relations}
    \begin{gathered}
        fe = \id_{N}, \qquad p_{i}e = 0, \qquad fq_{i} = 0, \qquad 0 \leq i < r,\\
        ef = p_{[r]}q_{[r]},
    \end{gathered}
\end{equation}
where \(p_{[r]} = p_{0}\cdots p_{r - 1}\) and \(q_{[r]} = q_{r - 1}\cdots q_{0}\).
In view of \zcref{thm:categorical-fock-decomposition}, we seek to identify the vacuum category \(\mathcal{V}\) by constructing categorical analogues of this extended data.
For this purpose, we recall the incidence semiorthogonal decomposition \cite[Theorem~3.2]{Jiang2023DerivedCategoriesGrassmannians}, established for derived Grassmannians over prestacks over \(\mathbb{Q}\).

Write \(B_{a,b}\) for the set of partitions of height at most \(a\) and width at most \(b\), and put
\[
    \mathcal{C}_{e}^{\vee} = \Perf\bigl(\Gr_{X}(\mathcal{E}^{\vee}[1],e)\bigr).
\]
For \(d \geq 1\), \(0 \leq i \leq \min\{r,d\}\), and \(\lambda \in B_{r - i,i}\), consider the incidence functors
\[
    \begin{aligned}
        \Phi_{d}^{(i,\lambda)} \colon \mathcal{C}_{d - i}^{\vee} & \longrightarrow \mathcal{C}_{d}, \\
        M & \longmapsto r_{+,*}\left(r_{-}^{*}M \otimes \mathbb{S}^{\lambda}(\mathcal{U}_{d,d - i})\right) \otimes (\det Q_{d})^{i}.
    \end{aligned}
\]
By \cite[Theorem 3.2]{Jiang2023DerivedCategoriesGrassmannians}, these functors are fully faithful and give a semiorthogonal decomposition
\[
    \mathcal{C}_{d} = \left\langle \im\Phi_{d}^{(i,\lambda)} \ \middle|\ 0 \leq i \leq \min\{r,d\}, \ \lambda \in B_{r - i,i} \right\rangle,
\]
with the ordering of the components given in \cite[Notation~3.1]{Jiang2023DerivedCategoriesGrassmannians}.

For \(d \geq r\), the case \(i = r\) consists of a single component corresponding to the empty partition.
Its embedding is
\[
    \Phi_{d}^{(r,\varnothing)} = \left(- \otimes (\det Q_{d})^{r}\right) r_{+,*}r_{-}^{*} \colon \mathcal{C}_{d - r}^{\vee} \longrightarrow \mathcal{C}_{d}.
\]

This suggests that the \(i = r\) component should realize the vacuum category in \zcref{thm:categorical-fock-decomposition}.
By the extended Clifford relations recalled above, we are led to consider its embedding and right adjoint as categorical analogues of the maps \(e\) and \(f\).
Put
\[
    \mathcal{C}^{\vee} = \bigoplus_{e \geq 0}\mathcal{C}_{e}^{\vee},
    \qquad \Phi_{\varnothing} = \bigoplus_{e \geq 0}\Phi_{e + r}^{(r,\varnothing)} \colon \mathcal{C}^{\vee} \longrightarrow \mathcal{C},
\]
and denote its right adjoint by \(\Phi_{\varnothing}^{R}\).

\begin{proposition}
    \label{prop:clifford-vacuum-incidence}
    The functors \(\Phi_{\varnothing}\) and \(\Phi_{\varnothing}^{R}\) satisfy the categorical extended Clifford relations
    \begin{equation}
        \label{eq:categorical-extended-clifford-relations}
        \begin{gathered}
            \Phi_{\varnothing}^{R}\Phi_{\varnothing} \simeq \id_{\mathcal{C}^{\vee}},\\
            F_{i}^{L}\Phi_{\varnothing} \simeq 0, \qquad \Phi_{\varnothing}^{R}F_{i} \simeq 0, \qquad 0 \leq i < r,\\
            \Phi_{\varnothing}\Phi_{\varnothing}^{R} \simeq F_{[r]}^{L}F_{[r]}.
        \end{gathered}
    \end{equation}
    In particular, \(\Phi_{\varnothing}\) induces an equivalence
    \[
        \mathcal{C}_{e}^{\vee} \xrightarrow{\sim}\mathcal{V}_{e + r}
    \]
    for every \(e \geq 0\).
\end{proposition}

\begin{proof}
    We recall local computations of \cite[Section~3]{Jiang2023DerivedCategoriesGrassmannians}.
    Working Zariski locally on \(X\), we may assume that \(X = \Spec R\) is affine and that \(\mathcal{E}\) admits a two-term presentation
    \[
        \mathcal{E} \simeq [R^{n - r}\xrightarrow{\delta}R^{n}],
    \]
    with \(R^{n}\) in degree \(0\).
    Writing \(B = R^{n}\), this presentation induces a natural affine morphism
    \[
        \iota_{d}:\Gr_{X}(\mathcal{E},d) \longrightarrow \Gr_{X}(B,d).
    \]
    Let \(R_{B,d}\) denote the universal kernel bundle on \(\Gr_{X}(B,d)\).
    By \cite[Lemma~3.10]{Jiang2023DerivedCategoriesGrassmannians}, the objects
    \[
        \iota_{d}^{*}\mathbb{S}^{\lambda}(R_{B,d}^{\vee}), \qquad \lambda \in B_{n - d,d},
    \]
    thickly generate \(\mathcal{C}_{d}\).
    In the following computations, we omit the pullback notation \(\iota_{d}^{*}\).

    By \cite[Lemma~3.11 and Corollary~3.12]{Jiang2023DerivedCategoriesGrassmannians}, the functor \(\Phi_{e + r}^{(r,\varnothing)}\) is fully faithful, with essential image thickly generated by
    \[
        (\det Q_{e + r})^{r} \otimes \mathbb{S}^{\lambda}(R_{B,e + r}^{\vee}), \qquad \lambda \in B_{n - e - r,e}.
    \]
    This gives the first relation.

    For \(0 \leq a < r\), the definition of the normalized flag functors and the projection formula give
    \[
        F_{r - 1 - a}^{L} \bigl((\det Q_{k})^{a} \otimes M\bigr) \simeq (\det Q_{k + 1})^{a + 1} \otimes \Psi_{k,0}^{L}(M).
    \]
    Applying \cite[Lemma~3.13(1)]{Jiang2023DerivedCategoriesGrassmannians} successively, we obtain
    \[
        F_{[r]}^{L} \bigl(\mathbb{S}^{\lambda}(R_{B,e}^{\vee})\bigr) \simeq
        \begin{cases}
            (\det Q_{e + r})^{r} \otimes \mathbb{S}^{\lambda}(R_{B,e + r}^{\vee}), &\lambda \in B_{n - e - r,e},\\
            0, &\text{otherwise}.
        \end{cases}
    \]
    Thus every generator of \(\im\Phi_{\varnothing}\) lies in the image of \(F_{[r]}^{L}\).
    By \eqref{eq:categorical-clifford-adjoint-vanishing}, we have
    \[
        F_{i}^{L}F_{[r]}^{L} \simeq 0, \qquad 0 \leq i < r,
    \]
    and hence \(F_{i}^{L}\Phi_{\varnothing} \simeq 0\).
    Taking right adjoints gives \(\Phi_{\varnothing}^{R}F_{i} \simeq 0\), proving the second and third relations.

    Moreover, the same calculation shows that the image of \(F_{[r]}^{L}\) is generated by objects in \(\im\Phi_{\varnothing}\).
    Since \(\Phi_{\varnothing}\) is fully faithful, the counit
    \[
        \Phi_{\varnothing}\Phi_{\varnothing}^{R}F_{[r]}^{L} \longrightarrow F_{[r]}^{L}
    \]
    is an equivalence.
    The equivalence of the adjunction unit, the vanishing relations above, and this counit equivalence are fpqc-local on \(X\), since the functors and adjunction maps commute with base change.
    Thus these conclusions hold over arbitrary \(X\).

    By \zcref{lem:fock-filtration-compatibility}, we have
    \[
        \mathcal{V} = \im F_{[r]}^{L}.
    \]
    Consequently, the counit \(\Phi_{\varnothing}\Phi_{\varnothing}^{R} \to \id\) is an equivalence on \(\mathcal{V}\).
    On the other hand, \(\Phi_{\varnothing}^{R}F_{i} \simeq 0\) implies that \(\Phi_{\varnothing}\Phi_{\varnothing}^{R}\) vanishes on every component \(F_{I}\mathcal{V}\) with \(I \neq \varnothing\).
    Hence, by \zcref{thm:categorical-fock-decomposition}, it is the right projection onto \(\mathcal{V}\).

    By the same theorem and \zcref{lem:flag-word-vanishing}, the functor \(F_{[r]}^{L}F_{[r]}\) is also the right projection onto \(\mathcal{V}\).
    Therefore
    \[
        \Phi_{\varnothing}\Phi_{\varnothing}^{R} \simeq F_{[r]}^{L}F_{[r]},
    \]
    proving the last relation.

    Finally, taking essential images gives
    \[
        \im\Phi_{\varnothing} = \im(\Phi_{\varnothing}\Phi_{\varnothing}^{R}) = \im(F_{[r]}^{L}F_{[r]}) = \mathcal{V}.
    \]
    Since \(\Phi_{\varnothing}\) is fully faithful, it induces the claimed equivalences
    \[
        \mathcal{C}_{e}^{\vee} \xrightarrow{\sim}\mathcal{V}_{e + r}.
    \]
\end{proof}
The identification of the vacuum category also allows us to compare our semiorthogonal decompositions with those of \cite[Theorem 3.2]{Jiang2023DerivedCategoriesGrassmannians}.
\begin{proposition}
    \label{prop:subset-flag-words}
    For every \(e \geq 0\) and \(\lambda \in B_{r - i,i}\), \(0 \leq i \leq r\), there is a natural equivalence
    \begin{equation}
        \label{eq:subset-flag-word}
        \Phi_{e + i}^{(i,\lambda)} \simeq F_{I}\Phi_{\varnothing}|_{\mathcal{C}_{e}^{\vee}},
    \end{equation}
    where \(I \subseteq [r]\) corresponds to \((i,\lambda)\) under the subset--partition bijection.
    In particular, our semiorthogonal decompositions agree with those of \cite[Theorem 3.2]{Jiang2023DerivedCategoriesGrassmannians}, including the embedding functors and the ordering of the components.
\end{proposition}

\begin{proof}
    The case \(I = \varnothing\) is immediate.
    Fix \(I = \{i_{1} < \cdots < i_{m}\} \subseteq [r]\) and \(e \geq 0\).
    Put \(d = e + r - m\), and define the partition
    \[
        \lambda = (\lambda_{1},\ldots,\lambda_{m}), \qquad \lambda_{j} = i_{m + 1 - j} - m + j, \quad 1 \leq j \leq m.
    \]
    Then \(\lambda \in B_{m,r - m}\).

    By the definition of the normalized flag functors and the projection formula,
    \[
        F_{t}\bigl(M \otimes (\det Q_{k + 1})^{a}\bigr) \simeq \Psi_{k,t - r + a}(M) \otimes (\det Q_{k})^{a - 1}.
    \]
    Applying this identity successively gives
    \[
        (F_{I}\Phi_{\varnothing})(M) \simeq (\det Q_{d})^{r - m} \otimes \Psi_{d,\lambda_{1}}\cdots \Psi_{d + m - 1,\lambda_{m}} (r_{+,*}r_{-}^{*}M),
    \]
    where the last push-pull functor is taken along \(I_{d + m,e}\).
    By \cite[Corollary~2.18]{Jiang2023DerivedCategoriesGrassmannians}, this is precisely \(\Phi_{d}^{(r - m,\lambda)}(M)\).

    As \(I\) ranges over the subsets of \([r]\), the pairs \((r - m,\lambda)\) exhaust the index set of \cite[Theorem 3.2]{Jiang2023DerivedCategoriesGrassmannians}, and the reverse colexicographic ordering agrees with the ordering of \cite[Notation~3.1]{Jiang2023DerivedCategoriesGrassmannians}.
    The last assertion follows from \zcref{prop:clifford-vacuum-incidence,thm:categorical-fock-decomposition}.
\end{proof}

Finally, passing to \(K_{0}\), we obtain an isomorphism
\begin{equation}
    \label{eq:jlz-spin-module}
    K_{0}(\mathcal{C}) \simeq K_{0}(\mathcal{C}^{\vee}) \otimes_{R}\bigwedge\nolimits_{R}^{\bullet}R^{r}, \qquad R = K_{0}(\Perf(X)).
\end{equation}
Under this identification, the operators \(f_{i}\) and \(\bar{f}_{i}^{*}\) of \zcref{cor:finite-clifford-action} act by exterior multiplication and contraction, respectively, recovering the spin representation of Jiang, Li, and Zhao \cite[Section~5.3 and Theorem~5.5]{JiangEtAl2025InstantonsAffineVertexAlgebra}.